# Input-Output-to-State Stability

Mikhail Krichman[*] and Eduardo D. Sontag[†]
Dep. of Mathematics, Rutgers University, NJ
{krichman,sontag}@math.rutgers.edu

Yuan Wang[‡]
Dept. of Mathematics, Florida Atlantic University, FL
ywang@control.math.fau.edu

**Abstract**

This work explores Lyapunov characterizations of the input-output-to-state stability (IOSS) property for nonlinear systems. The notion of IOSS is a natural generalization of the standard zero-detectability property used in the linear case. The main contribution of this work is to establish a complete equivalence between the input-output-to-state stability property and the existence of a certain type of smooth Lyapunov function. As corollaries, one shows the existence of "norm-estimators", and obtains characterizations of nonlinear detectability in terms of relative stability and of finite-energy estimates.

## 1 Introduction

This paper concerns itself with the following question, for dynamical systems: *is it possible to estimate, on the basis of external information provided by past input and output signals, the magnitude of the internal state $x(t)$ at time $t$?* The rest of this introduction will explain, in very informal and intuitive terms, the motivation for this question, closely related to the "zero detectability" problem, sketching the issues that arise and the main results. Precise definitions are provided in the next section.

State estimation is central to control theory. It arises in signal processing applications (Kalman filters), as well as in stabilization based on partial information (observers). By and large, the theory of state estimation is well-understood for linear systems, but it is still poorly developed for more general classes of systems, such as finite-dimensional deterministic systems, with which this paper is concerned. An outstanding open question is the derivation of useful necessary and sufficient conditions for the existence of observers, i.e., "algorithms" (dynamical systems) which converge to an estimate $\hat{x}(t)$ of the state $x(t)$ of the system of interest, using the information provided by $\{u(s), s \leq t\}$, the set of past input values, and by $\{y(s), s \leq t\}$, the set of past output measurements. In the context of stabilization to an equilibrium, let us say to the zero state $x = 0$ if we are working in an Euclidean space, a weaker type of estimate is sometimes enough: it may suffice to have a norm-estimate, that is to say, an upper bound $\hat{x}(t)$ on the *magnitude* (norm) $|x(t)|$ of the state $x(t)$. Indeed, it is often the case (cf. [33] and

---

[*]Supported in part by US Air Force Grant F49620-98-1-0421
[†]Supported in part by US Air Force Grant F49620-98-1-0242
[‡]Supported in part by NSF Grant DMS-9457826

Assumption UEC (73) in [19]) that norm-estimates suffice for control applications. To be more precise, one wishes that $\hat{x}(t)$ becomes eventually an upper bound on $|x(t)|$ as $t \to \infty$. We are thus interested in *norm-estimators* which, when driven by the i/o data generated by the system, produce such an upper bound $\hat{x}(t)$, cf. Figure 1.

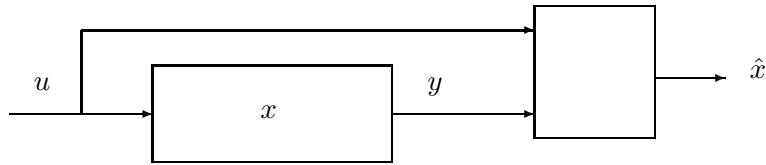

Figure 1: Norm-Estimator

In order to understand the issues that arise, let us start by considering the very special case when the external data (inputs $u$ and outputs $y$) vanish identically. The obvious estimate (assuming, as we will, that everything is normalized so that the zero state is an equilibrium for the unforced system, and the output is zero when $x = 0$) is $\hat{x}(t) \equiv 0$. However, the only way that this estimate fulfills the goal of upper bounding the norm of the true state as $t \to \infty$ is if $x(t) \to 0$. In other words, one obvious necessary property for the possibility of norm-estimation is that the origin must be a globally asymptotically stable state with respect to the "subsystem" consisting of those states for which the input $u \equiv 0$ produces the output $y \equiv 0$. One says in this case that the original system is *zero-detectable*. For *linear systems*, zero detectability is equivalent to detectability, that is to say, the property that if any two trajectories produce the same output, then they approach each other. Zero-detectability is a central property in the general theory of nonlinear stabilization on the basis of output measurements; see for instance, among many other references, [35, 18, 50, 11, 17]. Our work can be seen as a contribution towards the better characterization and understanding of this fundamental concept.

However, zero-detectability by itself is far from being sufficient for our purposes, since it fails to be "well-posed" enough. One easily sees that, at the least, one should ask that, when inputs and outputs are small, states should also be small, and if inputs and outputs converge to zero as $t \to \infty$, states do too, cf. Figure 2. Moreover, when defining formally the notion

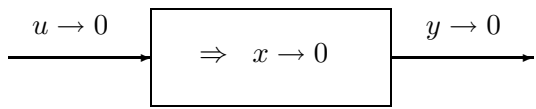

Figure 2: State Converges to Zero if External Data does of norm-estimator and the natural necessary and sufficient conditions for its existence, other requirements appear: the existence of asymptotic bounds on states, as a function of bounds on input/output data, and the need to describe the "overshoot" (transient behavior) of the state.

One way to approach the formal definition, so as to incorporate all the above characteristics in a simple manner, is to look at the analogous questions for the stability problem, which, for linear systems, is known to be technically dual to detectability. This leads one to the area which deals precisely with this circle of ideas: *input-to-state stability* (ISS).

Input-to-state stability was introduced in [37], and has proved to be a very useful paradigm in the study of nonlinear stability; see for instance the textbooks [17, 21, 23, 24], and the papers [10, 15, 16, 20, 29, 12, 33, 34, 44, 42, 49, 48], as well as its variants such as integral ISS (cf. [2, 4, 25, 39]) and input/output stability (cf. [37, 45, 46]). The notion of ISS takes into account the effect of initial states in a manner fully compatible with Lyapunov stability, and incorporates naturally the idea of "nonlinear gain" functions; the reader may wish to consult [40] for an exposition — as well as [44] for several new characterizations obtained after that exposition was written. Roughly speaking, a system is ISS provided that, no matter what is the initial state, if the inputs are small, then the state must eventually be small. Dualizing this definition one arrives at the notion of detectability which is the main subject of study of this paper: *input/output to state stability* (IOSS). (The terminology "IOSS" is not to be confused with the totally different concept called input/output stability (IOS), cf. [37, 45, 46], which refers instead to stability of outputs, rather than to detectability.)

A system $\dot{x} = f(x, u)$ with measurement ("output") map $y = h(x)$ is IOSS if there are some functions $\beta \in \mathcal{KL}$ and $\gamma_1, \gamma_2 \in \mathcal{K}_\infty$ such that the estimate:

$$|x(t)| \leq \max\left\{\beta(|x(0)|, t), \gamma_1\left(\|u|_{[0,t]}\|\right), \gamma_2\left(\|y|_{[0,t]}\|\right)\right\}$$

holds for any initial state $x(0)$ and any input $u(\cdot)$, where $x(\cdot)$ is the ensuing trajectory and $y(t) = h(x(t))$ the respective output function. (States $x(t)$, input values $u(t)$, and output values $y(t)$ lie in appropriate Euclidean spaces. We use $|\cdot|$ to denote Euclidean norm and $\|\cdot\|$ for supremum norm. Precise definitions and technical assumptions are discussed later.) The terminology IOSS is self-explanatory: formally there is "stability from the i/o data to the state". The term was introduced in the paper [47], but the same notion had appeared before: it represents a natural combination of the notions of "strong" observability (cf. [37]) and ISS, and was called simply "detectability" in [41] (where it is phrased in input/output, as opposed to state space, terms, and applied to questions of parameterization of controllers) and was called "strong unboundedness observability" in [20] (more precisely, this last notion allows also an additive nonnegative constant in the right-hand side of the estimate). In [47], two of the authors described relationships between the existence of full state observers and the IOSS property, or more precisely, a property which we called "incremental IOSS". The use of ISS-like formalism for studying observers, and hence implicitly the IOSS property, has also appeared several times in other authors' work, such as the papers [32, 27].

One of the main results of this paper is that a system is IOSS if and only if it admits a norm-estimator (in a sense also to be made precise). This result is in turn a consequence of a necessary and sufficient characterization of the IOSS property in terms of smooth dissipation functions, namely, there is a proper (radially unbounded) and positive definite smooth function $V$ of states (a "storage function" in the language of dissipative systems introduced by Willems [52] and further developed by Hill and Moylan [14, 13] and others) such that a *dissipation inequality*

$$\frac{d}{dt}V(x(t)) \leq -\sigma_1(|x(t)|) + \sigma_2(|y(t)|) + \sigma_3(|u(t)|) \tag{1}$$

holds along all trajectories, with the functions $\sigma_i$ of class $\mathcal{K}_\infty$. This provides an "infinitesimal" description of IOSS, and a norm-observer is easily built from $V$. Such a characterization in dissipation terms was conjectured in [47], and we provide here a complete solution to the problem. (The paper [47] also explains how the existence of $V$ links the IOSS property to "passivity" of systems.)

It is worth pointing out that several authors had independently suggested that one should *define* "detectability" in dissipation terms. For example, in [28], Equation 15, one finds detectability defined by the requirement that there should exist a differentiable storage function $V$ satisfying our dissipation inequality but with the special choice $\sigma_2(r) := r^2$ (there were no inputs in the class of systems considered there). A variation of this is to weaken the dissipation inequality, to require merely

$$x \neq 0 \;\Rightarrow\; \frac{d}{dt}V(x(t)) < \sigma_2(|y(t)|)$$

(again, with no inputs), as done for instance in the definition of detectability given in [31]. Observe that this represents a slight weakening of our property, in so far as there is no "margin" of stability $-\sigma_1(|x(t)|)$. One of our contributions is to show that such alternative definitions (when posed in the right generality) are in fact equivalent to IOSS.

A key preliminary step in the construction of $V$, just as it was for the analogous result for the ISS property obtained in [42], is the characterization of the IOSS property in robustness terms, by means of a "small gain" argument. The IOSS property is shown to be equivalent to the existence of a "robustness margin" $\rho \in \mathcal{K}_\infty$. This means that every system obtained by closing the loop with a feedback law $\Delta$ (even dynamic and/or time-varying) for which $|\Delta(t)| \leq \rho(|x(t)|)$ for all $t$, cf. Figure 3, is OSS (i.e., is IOSS as a system with no inputs). In order to formulate precisely this notion of robust detectability, we need to consider auxiliary "systems with disturbances". Since such systems must be introduced anyhow, we decided to present all our results (and definitions, even of IOSS) for systems with disturbances, in the process gaining extra generality in our results.

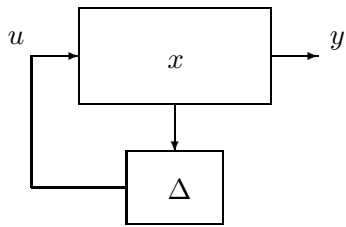

Figure 3: Robust Detectability

The core of the paper is, thus, the construction of $V$ for "robustly detectable" (more precisely, "robust IOSS") systems $\dot{x} = g(x,d)$ which are obtained by substituting $u = d\rho(|x|)$ in the original system, and letting $d = d(\cdot)$ be an arbitrary measurable function taking values in a unit ball. The function $V$ must satisfy a differential inequality of the form $\dot{V}(x(t)) \leq -\sigma_1(|x(t)|) + \sigma_2(|y(t)|)$ along all trajectories, that is to say, the following partial differential inequality:

$$\nabla V(x) \cdot g(x,d) \;\leq\; -\sigma_1(|x|) + \sigma_2(|y|),$$

for some functions $\sigma_1$ and $\sigma_2$ of class $\mathcal{K}_\infty$. But one last reduction consists of turning this problem into one of building Lyapunov functions for "relatively asymptotically stable" systems. Indeed, one observes that the main property needed for $V$ is that it should *decrease along trajectories as long as $y(t)$ is sufficiently smaller than $x(t)$*. This leads us to the notion of "global asymptotic stability modulo outputs" and its Lyapunov-theoretic characterization.

The construction of $V$ relies upon the solution of an appropriate optimal control problem, for which $V$ is the value function. This problem is obtained by "fuzzifying" the dynamics near the set where $y \ll x$, so as to obtain a problem whose value function is continuous. Several elementary facts about relaxed controls are used in deriving the conclusions. The last major ingredient is the use of techniques from nonsmooth analysis, and in particular inf-convolutions, in order to obtain a Lipschitz, and from there by a standard regularization argument, a smooth, function $V$, starting from the continuous $V$ that was obtained from the optimal control problem.

Finally, we will also discuss a version of detectability which relies upon "energy" estimates instead of uniform estimates. Such versions of detectability are fairly standard in control theory; see for instance [11], which defined "$L^2$-detectability" by a requirement that the state trajectory should be in $L^2$ if the observations are. The corresponding "integral to integral" notion uses a very interesting concept introduced in [30], that of "unboundedness observability" (UO), which amounts to a "relative (modulo outputs) forward completeness" property. It is shown that, for systems with no controls, the integral variant of OSS is equivalent to the conjunction of OSS and UO.

It is worth remarking that the main result in this paper amounts to providing necessary and sufficient conditions for the existence of a smooth (and proper and positive definite) solution $V$ to a partial differential inequality which is equivalent to asking that (1) holds along all trajectories, namely:

$$\max_{u \in \mathbb{R}^m} \{\nabla V(x) \cdot f(x,u) + \sigma_1(|x|) - \sigma_2(|h(x)|) - \sigma_3(|u|)\} \leq 0. \tag{2}$$

It is a consequence of our results that if there is an (even just) lower semicontinuous such solution (when "solution" is interpreted in a weak sense, for example in terms of viscosity or proximal subdifferentials), then there is also a smooth solution (usually, however, with different comparison functions $\sigma_i$'s). This is because the existence of a weak solution is already equivalent to IOSS, as shown in [22]. It is a routine observation that the above partial differential inequality can be posed in an equivalent way as a *Hamilton-Jacobi Inequality (HJI)*, in the special case of quadratic input "cost" $\sigma_3(r) = r^2$, and for systems $\dot{x} = f(x,u)$ which are affine in controls, i.e. systems of the form:

$$\dot{x} = g_0(x) + \sum_{i=1}^{m} u_i\, g_i(x) \tag{3}$$

(we are denoting by $u_i$ the $i$th component of $u$). Indeed, one need only replace the expression in (2) by its maximum value obtained at $u_i = (1/2)\nabla V(x) \cdot g_i(x)$, $i = 1, \ldots m$, thereby obtaining the following HJI:

$$\nabla V(x) \cdot g_0(x) + \frac{1}{4}\sum_{i=1}^{m}(\nabla V(x) \cdot g_i(x))^2 + \sigma_1(|x|) - \sigma_2(|h(x)|) \leq 0. \tag{4}$$

## 2 Definitions and statements of the main results

### 2.1 Systems of interest

We study a system whose dynamics depend on two types of inputs, which we call respectively *controls* and *disturbances*:

$$\dot{x}(t) = f(x(t), \mathbf{u}(t), \mathbf{w}(t)), \quad y(t) = h(x(t)). \tag{5}$$

Here, states evolve in $\mathbb{X} = \mathbb{R}^n$, controls are measurable, essentially bounded functions $\mathbf{u}$ on $\mathcal{I} = \mathbb{R}_{\geq 0}$ with values in $\mathbb{U} := \mathbb{R}^{m_u}$, and disturbances are measurable functions $\mathbf{w} : \mathcal{I} \to \Gamma$ with values in a compact metric space $\Gamma$ (which, unless otherwise specified, is always taken to be of the form $[-1,1]^{m_w}$); we will denote the set of all such functions by $\mathcal{M}_\Gamma$. In those cases when a different interval $\mathcal{I} \subset \mathbb{R}_{\geq 0}$ of definition for a control $\mathbf{u}$ is specified, we always apply the definitions to the extension of $\mathbf{u}$ to $\mathbb{R}_{\geq 0}$, using $\mathbf{u} \equiv 0$ on $\mathbb{R}_{\geq 0} \setminus \mathcal{I}$. The function $f : \mathbb{X} \times \mathbb{U} \times \Gamma \to \mathbb{X}$ is locally Lipschitz in $(x, u)$ uniformly on $w$, jointly continuous in $x$, $u$, and $w$, and such that $f(0, 0, w) = 0$ for any $w \in \Gamma$; and $h : \mathbb{X} \to \mathcal{Y} := \mathbb{R}^p$ is smooth ($C^1$) and vanishes at 0.

A function $\alpha : \mathbb{R}_{\geq 0} \to \mathbb{R}_{\geq 0}$ is *of class* $\mathcal{K}$ if it is continuous, positive definite, and strictly increasing, and is *of class* $\mathcal{K}_\infty$ if it is also unbounded. A function $\beta : \mathbb{R}_{\geq 0} \times \mathbb{R}_{\geq 0} \to \mathbb{R}_{\geq 0}$ is said to be *of class* $\mathcal{KL}$ if for each fixed $t \geq 0$, $\beta(\cdot, t)$ is of class $\mathcal{K}$, and for each fixed $s \geq 0$, $\beta(s, t)$ decreases to 0 as $t \to \infty$. Let $z(\cdot)$ be a measurable function.

The $L_\infty$ (essential supremum) norm of the restriction of $z$ to the interval $[t_1, t_2]$ is denoted by $\|z|_{[t_1, t_2]}\|$.

Given a state $\xi \in \mathbb{X}$, for each pair $(\mathbf{u}, \mathbf{w})$ denote by $x(t, \xi, \mathbf{u}, \mathbf{w})$ the unique maximal solution of the system (5), which is defined on some maximal interval $[0, t_{\max}(\xi, \mathbf{u}, \mathbf{w}))$. We will use the notation $y(t, \xi, \mathbf{u}, \mathbf{w}) := h(x(t, \xi, \mathbf{u}, \mathbf{w}))$, and, when unimportant or clear from the context, we will write $t_{\max}$ instead of $t_{\max}(\xi, \mathbf{u}, \mathbf{w})$, $x(t)$ instead of $x(t, \xi, \mathbf{u}, \mathbf{w})$ and $y(t)$ instead of $y(t, \xi, \mathbf{u}, \mathbf{w})$.

## 2.2 Notions of "Uniform Detectability" and dissipation functions

**Definition 2.1** A system of type (5) is said to be *uniformly input-output-to-state stable (UIOSS)* if, there exist functions $\beta \in \mathcal{KL}$ and $\gamma_1, \gamma_2 \in \mathcal{K}$ such that the estimate:

$$|x(t, \xi, \mathbf{u}, \mathbf{w})| \leq \max\left\{\beta(|\xi|, t), \gamma_1\left(\|\mathbf{u}|_{[0,t]}\|\right), \gamma_2\left(\|y|_{[0,t]}\|\right)\right\} \tag{6}$$

holds for any initial state $\xi \in \mathbb{X}$, control $\mathbf{u}$, disturbance $\mathbf{w}$, and time $t \in [0, t_{\max}(\xi, \mathbf{u}, \mathbf{w}))$.

**Definition 2.2** A smooth ($C^\infty$) function $V : \mathbb{X} \to \mathbb{R}_{\geq 0}$ is an *UIOSS-Lyapunov function* for system (5) if:

- there exist $\mathcal{K}_\infty$-functions $\alpha_1$, $\alpha_2$ such that

$$\alpha_1(|\xi|) \leq V(\xi) \leq \alpha_2(|\xi|) \tag{7}$$

  holds for all $\xi$ in $\mathbb{X}$, and

- there exists a $\mathcal{K}_\infty$-function $\alpha$ and $\mathcal{K}$-functions $\sigma_1$, $\sigma_2$ such that

$$\nabla V(\xi) \cdot f(\xi, u, w) \leq -\alpha(|\xi|) + \sigma_1(|u|) + \sigma_2(|h(\xi)|) \tag{8}$$

  for all $\xi$ in $\mathbb{X}$, all control values $u \in \mathbb{U}$, and all disturbance values $w \in \Gamma$. $\square$

Property (7) amounts to positive definiteness and properness of $V$; requiring the existence of an upper bound $\alpha_2$ is redundant, as it follows from the fact that $V$ is continuous and satisfies $V(0) = 0$. However, it is convenient to specify this bound explicitly, as it will be used in various estimates. Condition (8) is a *dissipation inequality* in the sense of [52].

**Remark 2.3** A smooth function $V : \mathbb{X} \to \mathbb{R}_{\geq 0}$, satisfying (7) on $\mathbb{X}$ with some $\alpha_1$, $\alpha_2$ of class $\mathcal{K}_\infty$, is an UIOSS-lyapunov function for a system (5) if and only if there exist functions $\alpha_3$ of class $\mathcal{K}_\infty$, and $\gamma$, and $\chi_1$ of class $\mathcal{K}$ such that

$$\nabla V(\xi) \cdot f(\xi, u, w) \leq -\alpha_3(|\xi|) + \gamma(|h(\xi)|), \tag{9}$$

for any $\xi \in \mathbb{X}$, $w \in \Gamma$, and $u \in \mathbb{U}$ such that $|\xi| \geq \chi_1(|u|)$.

Indeed, clearly (8) implies (9) with $\alpha_3(\cdot) := \alpha(\cdot)/2$, $\gamma \equiv \sigma_2$, and $\chi_1 := \alpha^{-1} \circ (2\sigma_1)$. To prove the other implication, assume now that (9) holds with some $\alpha_3 \in \mathcal{K}_\infty$ and $\gamma, \chi_1 \in \mathcal{K}$. Define $\sigma_1(\cdot) = \max\{0, \hat{\sigma}_1(\cdot)\}$, where

$$\hat{\sigma}_1(r) := \max\{\nabla V(\xi) \cdot f(\xi, u, w) + \alpha_3(\chi_1(|u|)) : |u| \leq r, |\xi| \leq \chi_1(r), w \in \Gamma\}.$$

Then $\sigma_1$ is continuous, $\sigma_1(0) = 0$, and one can assume that $\sigma_1$ is a $\mathcal{K}_\infty$-function (majorize it by one if it is not). We claim that (8) holds with $\alpha \equiv \alpha_3$ and $\sigma_2 \equiv \gamma$. Indeed, if $|\xi| \geq \chi_1(|u|)$, then (9) holds, from which (8) trivially follows. If $|\xi| < \chi_1(|u|)$, then, by definition of $\sigma_1$,

$$\sigma_1(|u|) \geq \nabla V(\xi) \cdot f(\xi, u, w) + \alpha_3(\chi_1(|u|))$$

for every $w$, which, in turn, implies (8). $\square$

A few particular cases of the UIOSS property have been studied in the literature. If the system (5) in consideration has no outputs and no disturbances, UIOSS reduces to the well-known ISS property, whose Lyapunov characterization was obtained in [42]. In case (5) is autonomous, UIOSS becomes OSS. This property was introduced in [43] where Lyapunov-type necessary and sufficient conditions were obtained. Finally, for systems with no disturbances, UIOSS is just IOSS. This property was introduced in [47], where it was conjectured that any IOSS control system admits a smooth IOSS-Lyapunov function. This conjecture will be proven here in a more general setting, for systems forced by both controls and disturbances. A few interesting applications of this Lyapunov characterization were also discussed in [47], one of them to be defined next.

## 2.3 Norm estimators

**Definition 2.4** A *state-norm estimator* (or *state-norm observer*) for a system $\Sigma$ of type (5) is a pair $(\Sigma_{n.o}, k(\cdot, \cdot))$, where $k : \mathbb{R}^\ell \times \mathcal{Y} \to \mathbb{R}$, and $\Sigma_{n.o}$ is a system

$$\dot{p} = g(p, u, y) \tag{10}$$

evolving in $\mathbb{R}^\ell$ and driven by the controls and outputs of $\Sigma$, such that the following conditions are satisfied:

- There exist $\mathcal{K}$-functions $\hat{\gamma}_1$ and $\hat{\gamma}_2$ and a $\mathcal{KL}$-function $\hat{\beta}$ such that for any initial state $\zeta \in \mathbb{R}^\ell$, all inputs $\mathbf{u}$ and $\mathbf{y}$, and any $t$ in the interval of definition of the solution $p(\cdot, \zeta, \mathbf{u}, \mathbf{y})$, the following inequality holds

$$|k(p(t, \zeta, \mathbf{u}, \mathbf{y}), y(t))| \leq \hat{\beta}(|\zeta|, t) + \hat{\gamma}_1\left(\|\mathbf{u}|_{[0,t]}\|\right) + \hat{\gamma}_2\left(\|\mathbf{y}|_{[0,t]}\|\right) \tag{11}$$

(in other words, the system (10) is IOS with respect to the inputs $u$ and $y$ and output $k$).

- There are functions $\rho \in \mathcal{K}$ and $\beta \in \mathcal{KL}$ so that, for any pair of initial states $\xi$ and $\zeta$ of systems (5) and (10) respectively, any control $\mathbf{u} : \mathbb{R}_{\geq 0} \to \mathbb{U}$, and any disturbance $\mathbf{w} \in \mathcal{M}_\Gamma$, we have

$$|x(t, \xi, \mathbf{u}, \mathbf{w})| \leq \beta(|\xi| + |\zeta|, t) + \rho(|k(p(t, \zeta, \mathbf{u}, \mathbf{y}_{\xi, \mathbf{u}, \mathbf{w}}), \mathbf{y}_{\xi, \mathbf{u}, \mathbf{w}}(t))|) \tag{12}$$

for all $t \in [0, t_{\max}(\xi, \mathbf{u}, \mathbf{w}))$. (Here $\mathbf{y}_{\xi, \mathbf{u}, \mathbf{w}}$ denotes the output trajectory of $\Sigma$, that is, $\mathbf{y}_{\xi, \mathbf{u}, \mathbf{w}}(t) = y(t, \xi, \mathbf{u}, \mathbf{w})$.) □

## 2.4 Statement of the main result

The main theorem to be proved in this paper, summarizing the equivalent characterizations of UIOSS, will be as follows:

**Theorem 1** *Let $\Sigma$ be a system of type (5). Then the following are equivalent:*

1. $\Sigma$ is UIOSS.

2. $\Sigma$ admits an UIOSS-Lyapunov function.

3. There is a state-norm estimator for $\Sigma$. ∎

The main contribution is in showing that 1 implies 2; the remaining implications are much easier.

## 2.5 Example: linear systems

A particular class of systems (5) is as follows. A *linear*, time-invariant system $\Sigma_{lin}$ with outputs is one for which $f$ and $h$ are linear, that is,

$$\begin{aligned} \dot{x} &= \mathbf{A}x + \mathbf{B}u \\ y &= \mathbf{C}x, \end{aligned} \tag{13}$$

where $\mathbf{A} \in \mathbb{R}^{n \times n}$, $\mathbf{B} \in \mathbb{R}^{n \times m}$, and $\mathbf{C} \in \mathbb{R}^{p \times n}$. (We assume that $m_w = 0$.)

Recall:

**Definition 2.5** A linear, time invariant system (13) is *detectable*, or *asymptotically observable*, if the following implication

$$\mathbf{C}x(t) \equiv 0 \implies x(t) \to 0 \tag{14}$$

holds for any trajectory $x(t)$ of (13), corresponding to the zero control $\mathbf{u} \equiv 0$. □

This is a totally routine linear systems theory fact, but we include the proof as a motivation for the nonlinear material to follow.

**Proposition 2.6** *If a linear system (13) is detectable, then it is IOSS.*

*Proof.* It is a well known fact (see, for example, [38]) that if a system (13) is detectable, then there exists a matrix $\mathbf{L} \in \mathbb{R}^{n \times p}$, such that the matrix $\mathbf{A} + \mathbf{LC}$ is Hurwitz, and, furthermore, the system

$$\dot{z}(t) = \mathbf{A}z(t) + \mathbf{B}u(t) + \mathbf{L}(\mathbf{C}z(t) - y(t)) \qquad (15)$$

referred to as an *observer* and driven by the controls and outputs of (13), has the property that if $x(t)$ and $z(t)$ are any solutions of (13) and (15) respectively, then $|x(t) - z(t)| \to 0$, and, in particular, if $x(0) = z(0)$, then $x(t) = z(t)$ for all nonnegative $t$. Fix an initial state $\xi$ and a control $\mathbf{u}$. Then the solution $x(t, \xi, \mathbf{u})$ of (13) is also the solution of (15) with $z(0) = \xi$, so that

$$x(t, \xi, \mathbf{u}) = e^{t(\mathbf{A}+\mathbf{LC})}\xi + \int_0^t e^{s(\mathbf{A}+\mathbf{LC})}[\mathbf{B}u(t-s) - \mathbf{L}y(t-s)]ds.$$

Choose two positive numbers $\delta'$ and $\delta$ so that $\Re\lambda \leq -\delta' < -\delta$ for every eigenvalue $\lambda$ of $\mathbf{A}+\mathbf{LC}$. Then there exists a polynomial $P(\cdot)$ and, consequently, a constant $K$ such that

$$\begin{aligned}
|x(t, \xi, \mathbf{u})| &\leq P(t)e^{-\delta' t}|\xi| + \int_0^t P(s)e^{-\delta' s}\left[\|\mathbf{B}\||u(t-s)| + \|\mathbf{L}\||y(t-s)|\right]ds \\
&\leq Ke^{-\delta t}|\xi| + K\frac{\|\mathbf{B}\|}{\delta}\|\mathbf{u}|_{[0,t]}\| + K\frac{\|\mathbf{L}\|}{\delta}\|y|_{[0,t]}\|.
\end{aligned}$$

Thus, the IOSS estimate (6) holds for (13) with the linear gains $\beta(r, t) := Ke^{-\delta t}r$ and $\gamma_1(r) = \gamma_2(r) := K\frac{\|\mathbf{B}\|}{\delta}r$. ∎

To find an IOSS-Lyapunov function for system (13), take any symmetric matrix $\mathbf{P} \in \mathbb{R}^{n \times n}$ such that

$$\mathbf{P}(\mathbf{A}+\mathbf{LC}) + (\mathbf{A}+\mathbf{LC})'\mathbf{P} = -I$$

(such a matrix $\mathbf{P}$ exists, because $\mathbf{A}+\mathbf{LC}$ is Hurwitz). Define

$$V(x) := x'\mathbf{P}x \qquad (16)$$

Notice that, since $\mathbf{P}(\mathbf{A}+\mathbf{LC}) = ((\mathbf{A}+\mathbf{LC})'\mathbf{P})'$, we have $x'\mathbf{P}(\mathbf{A}+\mathbf{LC})x = x'(\mathbf{A}+\mathbf{LC})'\mathbf{P}x$. Therefore

$$\begin{aligned}
\nabla V(x) \cdot f(x, u) &= 2x'\mathbf{P}(\mathbf{A}x + \mathbf{B}u) \\
&= 2x'\mathbf{P}((\mathbf{A}+\mathbf{LC})x + \mathbf{B}u - \mathbf{L}y) \\
&= 2x'\mathbf{P}(\mathbf{A}+\mathbf{LC})x + 2x'\mathbf{PB}u - 2x'\mathbf{PL}y \\
&= x'(\mathbf{P}(\mathbf{A}+\mathbf{LC}) + (\mathbf{A}+\mathbf{LC})'\mathbf{P})x + 2x'\mathbf{PB}u - 2x'\mathbf{PL}y \\
&\leq -|x|^2 + 2\|\mathbf{P}\|\|\mathbf{B}\||u||x| + 2\|\mathbf{P}\|\|\mathbf{L}\||y||x| \\
&\leq -|x|^2 + |x|^2/4 + 4\|\mathbf{P}\|^2\|\mathbf{B}\|^2|u|^2 + |x|^2/4 + 4\|\mathbf{P}\|^2\|\mathbf{L}\|^2|y|^2 \\
&\leq -|x|^2/2 + 4\|\mathbf{P}\|^2\|\mathbf{B}\|^2|u|^2 + 4\|\mathbf{P}\|^2\|\mathbf{L}\|^2|y|^2.
\end{aligned}$$

So, the UIOSS dissipation inequality holds for $V$ with gains defined by $\alpha(r) = r/2$, $\sigma_1(r) = 4\|\mathbf{P}\|^2\|\mathbf{B}\|^2r^2$, $\sigma_2(r) = 4\|\mathbf{P}\|^2\|\mathbf{L}\|^2r^2$.

## 2.6 Systems without controls

Let $\Omega$ be a compact metric space (which is always assumed to be of the form $[-1, 1]^m$ unless specified otherwise). Consider systems of the type

$$\dot{x} = f(x(t), \mathbf{d}(t)), \quad y(t) = h(x(t)), \qquad (17)$$

where $f : \mathbb{X} \times \Omega$ is locally Lipschitz in $x$ uniformly on $d$ and jointly continuous in $x$ and $d$, and $f(0, d) = 0$ for any $d \in \Omega$. The inputs are measurable functions $\mathbf{d} : \mathcal{I} \to \Omega$, and we use the term *disturbances* to refer to such $\Omega$-valued inputs. We will use $\mathcal{M}_\Omega$ to denote the collection of all such functions.

This system can be seen as a particular case of (5) that eschews controls and is driven only by disturbances. However, it will play an important role in our studies, therefore for convenience we will define the corresponding stability property and dissipation inequality for this system separately from the main definition 2.1.

**Definition 2.7** A system (17) is *uniformly output-to-state stable (UOSS)* if there exist some $\beta \in \mathcal{KL}$ and $\gamma_2 \in \mathcal{K}$ such that

$$|x(t, \xi, \mathbf{d})| \leq \max \left\{ \beta(|\xi|, t), \gamma_2 \left( \left\| y|_{[0,t]} \right\| \right) \right\} \tag{18}$$

for any disturbance $\mathbf{d}$, initial state $\xi \in \mathbb{X}$, and $t \in [0, t_{\max})$.

**Definition 2.8** A *UOSS-Lyapunov function* for system (17) is a smooth function $V : \mathbb{X} \to \mathbb{R}_{\geq 0}$ satisfying (7) and

$$\nabla V(x) \cdot f(x, d) \leq -\alpha_3(|x|) + \gamma(|h(x)|) \quad \forall\, x \in \mathbb{X}, \, \forall\, d \in \Omega, \tag{19}$$

with some class $\mathcal{K}_\infty$ functions $\alpha_i$ and a $\mathcal{K}$ function $\gamma$. For systems with no disturbances we simply say that $V$ is an OSS-Lyapunov function. □

### 2.6.1 "Modulo outputs" relative stability

Recall the classical notion of uniform global asymptotic stability for systems of type (17), ensuring that every solution of the system tends to the equilibrium and never goes too far from it. Suppose now that it does not matter how the system behaves when the information provided by the output is adequate, that is, the norm of the output dominates the norm of the current state. On the other hand, we want the system to decay nicely when the output does not help in determining how large the state is. This motivates the following "modulo output" definition of stability.

**Definition 2.9** A system of type (17) satisfies the GASMO *(global asymptotic stability modulo output)* property if there exist a function $\rho$ of class $\mathcal{K}_\infty$ and a function $\lambda$ of class $\mathcal{KL}$, such that, for all $\xi \in \mathbb{X}$, $\mathbf{d} \in \mathcal{M}_\Omega$, and any $T < t_{\max}(\xi, \mathbf{d})$, if

$$|x(t, \xi, \mathbf{d})| \geq \rho(|h(x(t, \xi, \mathbf{d}))|) \quad \forall\, 0 \leq t \leq T,$$

then the estimate

$$|x(t, \xi, \mathbf{d})| \leq \lambda(|\xi|, t) \quad \forall\, 0 \leq t \leq T. \tag{20}$$

holds. □

**Remark 2.10** If a system in consideration has no outputs, then the GASMO property becomes global asymptotic stability (GAS).

The following proposition provides an "ε-δ" characterization of the GASMO property.

**Proposition 2.11** A system of type (17) satisfies the GASMO property if and only if there exists a $\mathcal{K}_\infty$-function $\rho$ so that the following two properties hold:

1. For any $\varepsilon > 0$ and any $r > 0$, there exists some $T_{r,\varepsilon}$ such that for any $|\xi| \leq r$, any $\mathbf{d}$, and any $T \in [0, t_{\max}(\xi, \mathbf{d}))$ such that $T \geq T_{r,\varepsilon}$, if
$$|x(t, \xi, \mathbf{d})| \geq \rho(|y(t, \xi, \mathbf{d})|) \quad \text{for all } 0 \leq t \leq T,$$
then
$$|x(t, \xi, \mathbf{d})| < \varepsilon \quad \text{for all } t \in [T_{r,\varepsilon}, T].$$

2. There exists a $\mathcal{K}$-function $\vartheta$ such that for any $\xi \in \mathbb{X}$, any disturbance $\mathbf{d}$, and any $T < t_{\max}(\xi, \mathbf{d})$ such that
$$|x(t, \xi, \mathbf{d})| \geq \rho(|y(t, \xi, \mathbf{d})|) \quad \text{for all } 0 \leq t \leq T,$$
the following "bounded overshoot" estimate holds:
$$|x(T, \xi, \mathbf{d})| \leq \vartheta(|\xi|).$$

The necessity part is obvious. To prove the sufficiency, we need the following lemma, proved in section 3 of [26], although not explicitly stated in this form:

**Lemma 2.12** Let $\Phi(r, t) : (\mathbb{R}_{\geq 0})^2 \to \mathbb{R}_{\geq 0}$ be a map such that

1. for all $\varepsilon > 0$ and for all $R > 0$ there exists $T$ such that $\Phi(r, t) < \varepsilon$ for all $0 \leq r \leq R$ and for all $t \geq T$,

2. for all $\varepsilon > 0$ there exists $\delta > 0$ such that if $r \leq \delta$ then $\Phi(r, t) < \varepsilon$ for all $t > 0$.

Then $\Phi$ can be majorized by a $\mathcal{KL}$- function. □

*Proof of sufficiency for Proposition 2.11.* Consider the function:
$$\Phi(r, t) := \sup \{|x(t, \xi, \mathbf{d})| : |\xi| \leq r, \mathbf{d} \in \mathcal{M}_\Omega, |x(s, \xi, \mathbf{d})| \geq \rho(|y(s, \xi, \mathbf{d})|) \quad \forall s \in [0, t]\}.$$

Then the conditions 1) and 2) of Lemma 2.12 follow from the assumptions 1) and 2) of the proposition, so that one can majorize $\Phi$ by a $\mathcal{KL}$-function $\lambda$. ■

Suppose a system $\Sigma$ of type (17) satisfies the GASMO property with some $\mathcal{K}_\infty$ function $\rho$. By majorizing $\rho$ by another $\mathcal{K}_\infty$ function if necessary, we will assume that $\rho$ is smooth when restricted to $s > 0$ and also $\rho(s) > s$ for all positive $s$. We let
$$\mathcal{D} := \{\xi \in \mathbb{X}: |\xi| \leq \rho(|h(\xi)|)\},$$
$$\mathcal{E} := \mathbb{X} \setminus \mathcal{D},$$
and
$$\mathcal{E}_1 := \{\xi \in \mathbb{X} : |\xi| > 2\rho(|h(\xi)|)\}.$$

For each $\mathbf{d} \in \mathcal{M}_\Omega$ and $\xi \in \mathcal{E}$, define
$$\lambda_{\xi,\mathbf{d}} = \inf \{t \in [0, t_{\max}) : x(t, \xi, \mathbf{d}) \in \mathcal{D}\}, \tag{21}$$
with the convention $\lambda_{\xi,\mathbf{d}} = t_{\max}(\xi, \mathbf{d})$ if the trajectory never enters $\mathcal{D}$.

### 2.6.2 Integral variants

The UOSS property gives uniform estimates on states as a function of uniform bounds on outputs. There is a "finite energy output implies finite energy state" version as well:

**Definition 2.13** A system of type (17) is *integral to integral uniformly output to state stable (iiUOSS)* if there exist functions $\gamma$, $\kappa$ of class $\mathcal{K}$ and $\chi \in \mathcal{K}_\infty$ such that

$$\int_0^t \chi(|x(s,\xi,\mathbf{d})|)\, ds \;\leq\; \kappa(|\xi|) + \int_0^t \gamma(|h(x(s,\xi,\mathbf{d}))|)\, ds \tag{22}$$

for any initial state $\xi$, any disturbance $\mathbf{d} \in \mathcal{M}_\Omega$, and any time $t \in [0, t_{\max}(\xi, \mathbf{d}))$.  $\square$

Without loss of generality, $\gamma$ and $\kappa$ can be assumed to be of class $\mathcal{K}_\infty$.

**Definition 2.14** A system (5) is called *forward complete* if for every initial condition $\xi$, every input signal $\mathbf{u}$, and every disturbance $\mathbf{d}$ defined on $[0, +\infty)$, the corresponding trajectory $x(t, \xi, \mathbf{u}, \mathbf{d})$ is defined for all $t \geq 0$, i.e. $t_{\max}(\xi, \mathbf{u}, \mathbf{d}) = +\infty$.  $\square$

The following property, which is strictly weaker than forward completeness, was introduced in [30].

**Definition 2.15** A system (17) has the *unboundedness observability* property (UO) if

$$\limsup_{t \nearrow t_{\max}(\xi, \mathbf{d})} |y(t, \xi, \mathbf{d})| = +\infty \tag{23}$$

holds for each initial state $\xi$ and disturbance $\mathbf{d}$ with $t_{\max}(\xi, \mathbf{d}) < \infty$.  $\square$

The following useful characterization of UO was provided in [3].

**Proposition 2.16** A system (17) has the UO property if and only if there exist class $\mathcal{K}$ functions $\rho_1$, $\chi_1$, $\chi_2$, and a constant $c$, such that the following implication holds:

$$|h(x(t,\xi,\mathbf{d}))| \leq \rho_1(|x(t,\xi,\mathbf{d})|) \quad \forall t \in [0,T] \;\Rightarrow\; |x(t,\xi,\mathbf{d})| \leq \chi_1(t) + \chi_2(|\xi|) + c \quad \forall t \in [0,T], \tag{24}$$

for all $\xi \in \mathbb{X}$, $\mathbf{d} \in \mathcal{M}_\Omega$, and all $T \in [0, t_{\max}(\xi, \mathbf{d}))$.  $\square$

This proposition provides a uniform bound on all the states that can be reached by a UO system in given time from a given bounded set via a trajectory *not dominated by the output*. Notice that *for systems with disturbances (17), the UOSS property implies the UO property*.

### 2.6.3 Statement of the main result for the case of no inputs

**Theorem 2** *Let $\Sigma$ be a system of type (17). Then the following are equivalent:*

1. *$\Sigma$ is UOSS.*

2. *$\Sigma$ is GASMO.*

3. *$\Sigma$ is iiUOSS and UO.*

4. *$\Sigma$ admits an UOSS-Lyapunov function.*  ∎

## 2.7 Organization of the paper

Implications $2 \Rightarrow 3 \Rightarrow 1$ of Theorem 1 are proven in section 5. The most difficult to prove part of Theorem 1 is the implication $1 \Rightarrow 2$. The main technical result needed for this proof is implication $2 \Rightarrow 4$ of Theorem 2. This is proven in section 4. The construction of a UIOSS-Lyapunov function for an original system (5) is reduced, via a small gain argument, to the construction of a UOSS-Lyapunov function for a special system (17) related to the original system (5). This reduction is done in section 3.1, and section 3.2 completes the construction of UIOSS-Lyapunov functions.

Finally, implications $3 \Rightarrow 2$, $1 \Rightarrow 2$, and $4 \Rightarrow 3$ of Theorem 2 are proven in section 3.3, and $4 \Rightarrow 1$ follows from Theorem 1.

# 3 Reduction to the case of no controls

In this part we show how to reduce our main result to the particular case of systems with no controls.

## 3.1 Robust output to state stability

**Definition 3.1** System (5) is said to be *robustly output to state stable (ROSS)* if there exists a locally Lipschitz $\mathcal{K}_\infty$-function $\varphi$, called a *stability margin*, such that the system

$$\dot{x}(t) = g(x(t), \mathbf{d}(t)) := f(x(t), \mathbf{d}_u(t)\varphi(|x(t)|), \mathbf{w}(t)) \tag{25}$$

with disturbances $d := [d_u, w] \in \Omega := [-1, 1]^{m_u + m_w}$ and outputs $y = h(x)$ is UOSS.

Observe that the dynamics $g$ of system (25) are locally Lipschitz in $x$ uniformly in $d$, and also $g(0, d) = 0$ for all $d \in \Omega$.

**Lemma 3.2** If a system (5) is UIOSS, then it is ROSS.

The proof will follow from a few preliminary lemmas.

Let $\beta \in \mathcal{KL}$ and $\gamma_1, \gamma_2 \in \mathcal{K}_\infty$ be as in (6). Let $\alpha(r) = \beta(r, 0)$. Without loss of generality, we may assume that $\alpha$ is $\mathcal{K}_\infty$ and $\alpha(r) \geq r$ (so that $\alpha^{-1}(r) \leq r$).

Define $\varphi(r)$ to be a locally Lipschitz $\mathcal{K}_\infty$-function, which minorizes $\gamma_1^{-1}(\frac{1}{4}\alpha^{-1}(r))$ and can be extended as a Lipschitz function to a neighborhood of $[0, \infty)$. To prove the lemma we will show that $\varphi$ is a stability margin for (5).

**Proposition 3.3** Fix a $\xi \in \mathbb{X}$, a control $\mathbf{u}$, and a disturbance $\mathbf{w}$, and let $x(\cdot) := x(\cdot, \xi, \mathbf{u}, \mathbf{w})$ be the corresponding solution of the system (5). Let $T \in [0, t_{\max}(\xi, \mathbf{u}, \mathbf{w}))$. Then if $|\mathbf{u}(t)| \leq \varphi(|x(t)|)$ for almost all $t \in [0, T]$, the estimate

$$|x(t)| \leq \max\left\{\beta(|\xi|, t),\ \gamma_2(\|y|_{[0,t]}\|),\ \frac{|\xi|}{4}\right\} \tag{26}$$

holds for all $t \in [0, T]$.

*Proof.*

*Claim 1.* Suppose $T < t_{\max}(\xi, \mathbf{u}, \mathbf{w})$. If

$$|\mathbf{u}(t)| \leq \varphi(|x(t)|) \text{ for almost all } t \in [0, T), \tag{27}$$

then, for all $t \in [0, T)$,

$$|x(t)| \leq \max\left\{\alpha(|\xi|), 2\gamma_2(\|y|_{[0,t]}\|)\right\}. \tag{28}$$

*Proof of the Claim:* Suppose first that $\xi = 0$. In this case $x(t) \equiv 0$. Indeed, define $\mathbf{d}_u$ on $[0, t_{\max}(\xi, \mathbf{u}, \mathbf{w}))$ by

$$\mathbf{d}_u(t) = \begin{cases} 0, & \text{if } x(t) = 0 \\ \mathbf{u}(t)/\varphi(|x(t)|) & \text{if } x(t) \neq 0. \end{cases}$$

Then (27) implies that $\mathbf{d}_u \in \mathcal{M}_{[-1,1]^{m_u}}$, and that $x(\cdot)$ is the solution of (25) with $\xi = 0$, $\mathbf{w}$ as we picked, and $\mathbf{d}_u$ as we defined. Noticing that the constant function equal to $0$ is also a solution of this system, with the same initial state and the same disturbance, we conclude by uniqueness of solutions that $x(t) = 0$ for all nonnegative $t$, so that (28) trivially holds.

Suppose now that $\xi \neq 0$. Fix $\varepsilon$, such that $1 < \varepsilon < 2$. We will first show that, if $|\mathbf{u}(t)| \leq \varphi(|x(t)|)$ for almost all $t < T$, then the estimate

$$|x(t)| \leq \max\left\{\varepsilon\alpha(|\xi|), 2\gamma_2(\|y|_{[0,t]}\|)\right\} \tag{29}$$

holds for all $t \in [0, T)$. Indeed, notice that (29) is true as a strict inequality at $t = 0$, because $|\xi| \leq \alpha(|\xi|) < \varepsilon\alpha(|\xi|)$. If (29) fails at some $t \in [0, T)$, then there exists a

$$t_0 = \min\left\{t < T : x(t) = \max\left\{\varepsilon\alpha(|\xi|), 2\gamma_2(\|y|_{[0,t]}\|)\right\}\right\}.$$

Note that $t_0 > 0$ because at $t = 0$ we have a strict inequality in (29). So, (29) holds for all $t \in [0, t_0)$ and $x(t_0) = \max\left\{\varepsilon\alpha(|\xi|), 2\gamma_2(\|y|_{[0,t_0]}\|)\right\}$. Therefore for almost all $t \in [0, t_0)$ we have

$$\begin{aligned}
\gamma_1(\|\mathbf{u}|_{[0,t]}\|) &\leq \gamma_1(\|\varphi(|x(\cdot)|)|_{[0,t]}\|) \\
&\leq \max\left\{\frac{1}{4}\alpha^{-1}(\varepsilon\alpha(|\xi|)), \frac{1}{4}\alpha^{-1}(2\gamma_2(\|y|_{[0,t]}\|))\right\} \\
&\leq \max\left\{\frac{1}{4}\varepsilon\alpha(|\xi|), \frac{1}{4}2\gamma_2(\|y|_{[0,t]}\|)\right\} \\
&\leq \max\left\{\alpha(|\xi|), \gamma_2(\|y|_{[0,t]}\|)\right\}.
\end{aligned}$$

Then, since our system is UIOSS and $x(\cdot)$ is continuous, for all $t \in [0, t_0]$ we have

$$\begin{aligned}
|x(t)| &\leq \max\left\{\alpha(|\xi|), \gamma_1(\|\mathbf{u}|_{[0,t]}\|), \gamma_2(\|y|_{[0,t]}\|)\right\} \\
&= \max\left\{\alpha(|\xi|), \gamma_2(\|y|_{[0,t]}\|)\right\}.
\end{aligned}$$

On the other hand, $|x(t_0)| = \max\left\{\varepsilon\alpha(|\xi|), 2\gamma_2(\|y|_{[0,t_0]}\|)\right\}$ by definition of $t_0$. The contradiction proves the estimate (29). Letting $\varepsilon$ tend to $1$, we conclude that estimate (28) holds for all $t \in [0, T)$, completing the proof of Claim 1.

Hence, under the assumption of Claim 1, we have

$$\begin{aligned}
\gamma_1\left(\|\mathbf{u}|_{[0,t]}\|\right) &\leq \gamma_1\left(\|\varphi(|x(\cdot)|)|_{[0,t]}\|\right) \\
&\leq \max\left\{\frac{|\xi|}{4}, \frac{1}{4}\alpha^{-1}\left(2\gamma_2\left(\|y|_{[0,t]}\|\right)\right)\right\}
\end{aligned}$$

for all $t$ in $[0, t_{\max})$. So,

$$|x(t)| \leq \max\left\{\beta(|\xi|, t),\ \gamma_2\left(\|y|_{[0,t]}\|\right),\ \frac{|\xi|}{4},\ \frac{1}{4}\alpha^{-1}\left(2\gamma_2\left(\|y|_{[0,t]}\|\right)\right)\right\}$$

for all $t \in [0, t_{\max})$. Noticing that

$$\frac{1}{4}\alpha^{-1}\left(2\gamma_2\left(\|y|_{[0,t]}\|\right)\right) \leq \gamma_2\left(\|y|_{[0,t]}\|\right),$$

(because $\alpha^{-1}(r) \leq r$) we arrive at (26). ∎

**Lemma 3.4** Given any $\mathcal{KL}$-function $\hat{\beta}$, there exists a $\mathcal{KL}$-function $\beta$ and a $\mathcal{K}_\infty$-function $\nu$ such that for any $\tau > 0$, any continuous function $\mu : [0, \tau] \to \mathbb{R}_{\geq 0}$, and any nonnegative constant $C$, the following implication holds:

$$\forall t_1, t_2,\ 0 \leq t_1 < t_2 \leq \tau,\quad \mu(t_2) \leq \max\left\{\hat{\beta}(\mu(t_1), t_2 - t_1),\ \frac{\mu(t_1)}{2},\ C\right\} \tag{30}$$

implies

$$\mu(\tau) \leq \max\{\beta(\mu(0), \tau),\ \nu(C)\}. \tag{31}$$

*Proof.* By Proposition 7 in [39], there exist $\mu_1$ and $\mu_2 \in \mathcal{K}_\infty$ such that

$$\hat{\beta}(r, t) \leq \mu_1(\mu_2(r)e^{-t}),$$

so, by majorizing $\hat{\beta}$ as above if necessary, we can assume without loss of generality that $\hat{\beta}$ is continuous in its second variable, and $\hat{\beta}(r, 0) \geq r$ for all $r$. For any $r > 0$ define $T_r$ to be the first time when $\hat{\beta}(r, T_r) = r/2$. By replacing $\hat{\beta}$ with the $\mathcal{KL}$-function $\widetilde{\beta}$, defined by

$$\widetilde{\beta}(r, t) := \max\left\{\hat{\beta}(r, t),\ \hat{\beta}(r, 0)e^{-t}\right\},$$

we can assume without loss of generality that the series $\sum_{i=0}^{\infty} T_{\frac{r}{2^i}}$ diverges for every $r > 0$.

Define a function $\phi : \mathbb{R}_{\geq 0} \times \mathbb{R}_{\geq 0} \to \mathbb{R}_{\geq 0}$ as follows:

$$\phi(r, t) = \begin{cases} \hat{\beta}(r, t) & \text{for } t \in [0, T_r) \\ \hat{\beta}(\frac{r}{2^k}, t - \sum_{i=0}^{k-1} T_{\frac{r}{2^i}}) & \text{for } t \in \left[\sum_{i=0}^{k-1} T_{\frac{r}{2^i}}, \sum_{i=0}^{k} T_{\frac{r}{2^i}}\right), k = 1, 2, 3... \end{cases}.$$

Notice that the following two conditions hold for $\phi$:

1) For every $R, \varepsilon > 0$ there exists $\widetilde{t} > 0$ such that $\phi(r, t) \leq \varepsilon$ for all $r < R$ and $t > \widetilde{t}$.

Indeed, fix positive $R$ and $\varepsilon$ and find $k \in \mathbb{Z}$ such that $\hat{\beta}(R/2^k, 0) < \varepsilon$. Next, by continuity of $\hat{\beta}$ and by compactness of $[0, R]$ we can find a $\widetilde{t}$, such that $\sum_{i=0}^{k-1} T_{\frac{r}{2^i}} < \widetilde{t}$ for all positive $r < R$. Then, if $r < R$ and $t > \widetilde{t}$, then

$$\phi(r, t) = \hat{\beta}\left(\frac{r}{2^k}, t - \sum_{i=0}^{k-1} T_{\frac{r}{2^i}}\right) < \hat{\beta}\left(\frac{r}{2^k}, 0\right) < \varepsilon.$$

2) For all $\varepsilon > 0$ there exists $\delta > 0$ such that if $r \leq \delta$, then $\phi(r, t) < \varepsilon$ for all $t \geq 0$.

Indeed, for all $r$ and $t$, $\phi(r,t) \leq \hat{\beta}_0(r) := \hat{\beta}(r,0)$. For any positive $\varepsilon$, take $\delta = \delta(\varepsilon) := \hat{\beta}_0^{-1}(\varepsilon)$. Then, for all $t \geq 0$ we have
$$\phi(r,t) \leq \hat{\beta}\left(\hat{\beta}_0^{-1}(\varepsilon), 0\right) \leq \varepsilon.$$

Therefore, by Lemma 2.12, $\phi$ can be majorized by a $\mathcal{KL}$-function $\beta$.

Let $\nu(r) = \hat{\beta}(2r, 0)$.

Now pick any $\mu$, $C$ and $\tau$ satisfying (30). Define $T = \min\{t : \mu(t) \leq 2C\}$, and $T = \tau$ if $\mu(t) > 2C$ for all $t \geq 0$.

For any $t_1$ and $t_2$ in $[0,\tau]$ such that $0 \leq t_1 \leq t_2 \leq T$, we have $\mu(t_1) > 2C$, so that $\mu(t_1)/2 > C$, hence
$$\mu(t_2) \leq \max\left\{\hat{\beta}(\mu(t_1), t_2 - t_1), \mu(t_1)/2\right\}. \tag{32}$$

Suppose now that $\tau = T$. If $0 \leq \tau < T_{\mu(0)}$, then (32) with $t_1 = 0$, $t_2 = \tau$ yields
$$\mu(\tau) \leq \max\left\{\hat{\beta}(\mu(0), \tau), \frac{\mu(0)}{2}\right\} = \hat{\beta}(\mu(0), \tau),$$
where the equality follows from the definition of $T_{\mu(0)}$. Likewise, if
$$\tau \in \left[\sum_{i=0}^{k-1} T_{\frac{\mu(0)}{2^i}}, \sum_{i=0}^{k} T_{\frac{\mu(0)}{2^i}}\right),$$
then
$$\begin{aligned}\mu(\tau) &\leq \max\left\{\hat{\beta}\left(2^{-k}\mu(0), \tau - \sum_{i=0}^{k-1} T_{\frac{\mu(0)}{2^i}}\right), 2^{-(k+1)}\mu(0)\right\} \\ &= \hat{\beta}\left(2^{-k}\mu(0), \tau - \sum_{i=0}^{k-1} T_{\frac{\mu(0)}{2^i}}\right),\end{aligned}$$
where the inequality follows from (32) and the equality is implied by the definition of $T_{\frac{\mu(0)}{2^k}}$. Therefore we have
$$\mu(\tau) \leq \phi(\mu(0), \tau). \tag{33}$$

In case $\tau > T$, inequality (30) implies
$$\begin{aligned}\mu(\tau) &\leq \max\left\{\hat{\beta}(2C, \tau - T), \mu(T)/2, C\right\} \\ &= \max\left\{\hat{\beta}(2C, \tau - T), C, C\right\} \leq \hat{\beta}(2C, 0) = \nu(C). \tag{34}\end{aligned}$$

Combining (33) and (34) we obtain
$$\mu(\tau) \leq \max\{\phi(\mu(0), \tau), \nu(C)\} \leq \max\{\beta(\mu(0), \tau), \nu(C)\}.$$
∎

*Proof.* (of Lemma 3.2) We need to show that the system (25), corresponding to our system (5) with the stability margin $\varphi$ we have defined, is UOSS. Apply Lemma 3.4 to the $\mathcal{KL}$-function $\hat{\beta} := \beta$ to find appropriate functions $\beta_1 \in \mathcal{KL}$ and $\nu \in \mathcal{K}$. Assume given any initial state $\xi$ and disturbance $\mathbf{d} = [\mathbf{d}_u, \mathbf{w}]$, and let $x(t) := x(t, \xi, \mathbf{d}_u, \mathbf{w})$ be the corresponding solution. Fix any positive $t < t_{\max}(\xi, \mathbf{d})$, and define $\mathbf{u}$ by

$$\mathbf{u}(s) := \begin{cases} \varphi(|x(s)|)\mathbf{d}_u(s), & s \leq t \\ 0, & s > t. \end{cases}$$

Then, for all $s \leq t$ we have $x(s) = x(s, \xi, \mathbf{u}, \mathbf{w})$, where the latter is the solution of the original system (5) with control $\mathbf{u}$ and disturbance $\mathbf{w}$. Let $C = \gamma_2(\|y|_{[0,t]}\|)$. Then, for any $t_1$ and $t_2$ in $[0, t]$ we have $C \geq \gamma_2(\|y|_{[t_1, t_2]}\|)$. So, since $|\mathbf{u}(s)| \leq \varphi(|x(s)|)$ for all $s \in [0, t]$, Proposition 3.3 will imply that

$$|x(t_2)| \leq \max\left\{\beta(|x(t_1)|, t_2 - t_1), \frac{|x(t_1)|}{2}, C\right\}.$$

By the choice of $\beta_1$ and $\nu$ we have then,

$$|x(t)| \leq \max\left\{\beta_1(|\xi|, t), \nu(\gamma_2(\|y|_{[0,t]}\|))\right\},$$

proving the UOSS property for system (25) corresponding to the original UIOSS system. Thus, $\varphi$ is indeed a stability margin for the original system, and the proof of Lemma 3.2 is now complete. ∎

## 3.2 A UIOSS system admits a UIOSS-Lyapunov function

We show now how the main implication of Theorem 1 follows from Theorem 2.

**Lemma 3.5** (See Lemma 2.13 in [42]) Suppose a system $\Sigma$ of type (5) is ROSS. Let $V$ be a UOSS-Lyapunov function for the system (25) associated with $\Sigma$. Then $V$ is an UIOSS-Lyapunov function for $\Sigma$.

*Proof.* Let $\varphi$ be a stability margin for $\Sigma$. Since $V$ is a UOSS-Lyapunov function for (25), inequalities (7) and (19) hold with some $\alpha_1$, $\alpha_2$, $\alpha_3$ and $\gamma$. Pick a state $\xi \in \mathbb{X}$ and disturbance value $w \in \Gamma$. For any control value $u \in \mathbb{U}$ with $|u| \leq \varphi(|\xi|)$ we can find a $d_u \in [-1, 1]^{m_u}$ such that $u = d_u \varphi(|\xi|)$, so that by the dissipation inequality (19) for $V$ (applied with $d := [d_u, w]$) we have

$$\nabla V(\xi) \cdot f(\xi, u, w) = \nabla V(\xi) \cdot g(\xi, d) \leq -\alpha_3(|\xi|) + \gamma(|h(\xi)|),$$

proving (9) for $V$. So, the condition as in Remark 2.3 is satisfied for $V$ with $\chi_1 = \varphi^{-1}$, and $\alpha_i$ and $\gamma$ as before. Thus, $V$ is a UIOSS-Lyapunov function for $\Sigma$. ∎

By Theorem 2, the system (25) admits a UOSS-Lyapunov function $V$. Hence the following corollary follows.

**Corollary 3.6** If a system (5) is ROSS, then it admits a UIOSS-Lyapunov function. □

By Lemma 3.2, every UIOSS system is also ROSS, hence the implication $1 \Rightarrow 2$ of Theorem 1 follows.

## 3.3 UOSS and iiUOSS imply the GASMO property

**Lemma 3.7** A UOSS system of type (17) satisfies the GASMO property.

*Proof.* Assume that system (17) is UOSS. Without loss of generality, we may assume that $\gamma_2$ in Equation (18) is of class $\mathcal{K}_\infty$.

Let $\vartheta(s) = \beta(s, 0)$. Recall that we have assumed that $\vartheta(s) > s$ for all $s > 0$. Now let $\rho$ be any $\mathcal{K}_\infty$-function satisfying the inequality $\rho(s) > \vartheta(4\gamma_2(s))$ for all $s > 0$.

*Claim:* For any $\xi \in \mathbb{X}$, any $\mathbf{d} \in \mathcal{M}_\Omega$ and any $\tau \in [0, t_{\max}(\xi, \mathbf{d}))$, if

$$|x(t, \xi, \mathbf{d})| \geq \rho(|y(t, \xi, \mathbf{d})|) \text{ for all } 0 \leq t \leq \tau,$$

then

$$\gamma_2(|y(t, \xi, \mathbf{d})|) \leq |\xi|/2 \text{ for all } 0 \leq t \leq \tau,$$

and hence

$$|x(t, \xi, \mathbf{d})| \leq \beta(|\xi|, 0) = \vartheta(|\xi|) \text{ for all } 0 \leq t \leq \tau.$$

In particular, if $|x(t, \xi, \mathbf{d})| \geq \rho(|y(t, \xi, \mathbf{d})|)$ for all $t \in [0, t_{\max}(\xi, \mathbf{d}))$, then $t_{\max}(\xi, \mathbf{d}) = \infty$.

*Proof of the Claim:* If $\xi = 0$, the result is clear. Pick any $\xi \neq 0$, $\mathbf{d} \in \mathcal{M}_\Omega$ and assume that $|x(t, \xi, \mathbf{d})| \geq \rho(|y(t, \xi, \mathbf{d})|)$ for all $0 \leq t \leq \tau$ for some $\tau \in (0, t_{\max}(\xi, \mathbf{d}))$. Then, at $t = 0$,

$$\gamma_2(|y(0, \xi, \mathbf{d})|) \leq \gamma_2(\rho^{-1}(|\xi|)) \leq \gamma_2(\rho^{-1}(\vartheta(|\xi|))) < |\xi|/4.$$

Hence, $\gamma_2(|y(t, \xi, \mathbf{d})|) < |\xi|/4$ for all $t \in [0, \delta)$ for some $\delta > 0$. Let

$$t_1 = \inf\{t > 0 : \gamma_2(|y(t, \xi, \mathbf{d})|) \geq |\xi|/2\}.$$

Then $t_1 > 0$. Assume now that $t_1 \leq \tau$. Then

$$\gamma_2(|y(t_1, \xi, \mathbf{d})|) = |\xi|/2, \text{ and } \gamma_2(|y(t, \xi, \mathbf{d})|) < |\xi|/2,$$

for each $t \in [0, t_1)$, and hence for such $t$: $|x(t, \xi, d)| \leq \vartheta(|\xi|)$. Then, for each $0 \leq t \leq t_1$,

$$\gamma_2(|y(t, \xi, \mathbf{d})|) \leq \gamma_2(\rho^{-1}(|x(t, \xi, \mathbf{d})|)) \leq \gamma_2(\rho^{-1}(\vartheta(|\xi|))) < |\xi|/4.$$

By continuity, $\gamma_2(|y(t_1, \xi, \mathbf{d})|) \leq |\xi|/4$, contradicting the definition of $t_1$. This shows that it is impossible to have $t_1 \leq \tau$, and the proof of Claim is complete.

For each $r > 0$ let $T_r$ be any nonnegative number so that $\beta(r, t) < r/2$ for all $t \geq T_r$. Now, given any $r > 0$ and any $\varepsilon > 0$, for each $i = 1, 2, \ldots$, let $r_i := 2^{1-i}r$, and let $k(\varepsilon)$ be any positive integer so that $2^{-k(\varepsilon)}r < \varepsilon$ and define $T_{r,\varepsilon}$ as $T_{r_1} + T_{r_2} + \ldots + T_{r_{k(\varepsilon)}}$.

Pick any trajectory $x(t, \xi, \mathbf{d})$ as in the statement of Proposition 2.11, defined on an interval of the form $[0, T]$, with $T \geq T_{r,\varepsilon}$, with initial condition $|\xi| \leq r$ and disturbance $\mathbf{d} \in \mathcal{M}_\Omega$, satisfying $|x(t, \xi, \mathbf{d})| \geq \rho(|y(t, \xi, \mathbf{d})|)$ for all $t \in [0, T]$. Then, the above claim implies that $\gamma_2(|y(t, \xi, \mathbf{d})|) < |\xi|/2$ for all such $t$. Therefore, for any $t > T_{r_1} = T_r$,

$$\begin{aligned}|x(t, \xi, \mathbf{d})| &\leq \max\{\beta(|\xi|, t), |\xi|/2\} \\ &\leq \max\{\beta(r, t), r/2\} \leq r/2.\end{aligned}$$

Consider now the restriction of the trajectory to the interval $[T_{r_1}, T]$. This is the same as the trajectory that starts from the state $x(T_{r_1}, \xi, \mathbf{d})$, which has norm less than $r_1$, so by the same argument and the definition of $T_{r_2}$ we have that $|x(t, \xi, d)| \leq r/4$ for all $t \geq T_{r_2}$. Repeating on each interval $[T_{r_i}, T_{r_{i+1}}]$, we conclude that $|x(t, \xi, d)| < \varepsilon$ for all $T_{r,\varepsilon} \leq t \leq T$. ∎

**Lemma 3.8** Suppose a system of type (17) is iiUOSS and UO. Then it satisfies the GASMO property with $\rho(\cdot) := \max\left\{\chi^{-1}(2\gamma(\cdot)), \rho_1^{-1}(\cdot)\right\}$, where $\chi$ and $\gamma$ are as in the definition of iiUOSS and $\rho_1$ is as in Proposition 2.16.

To prove this lemma, we need the following elementary observation, which is a variant of what is usually referred to as "Barbălat's lemma":

**Proposition 3.9** Let $\mathcal{X} := \{x_\alpha, \ \alpha \in \mathcal{A}\}$ be a family of absolutely continuous curves in $\mathbb{X}$, each of which is defined on an interval $\mathcal{I}_\alpha$, either half-open ($\mathcal{I}_\alpha = [0, \lambda_\alpha)$) or closed ($\mathcal{I}_\alpha = [0, \lambda_\alpha]$). Suppose that

- $\mathcal{X}$ is closed with respect to shifts, that is, for all $\alpha \in \mathcal{A}$ and $T \in \mathcal{I}_\alpha$, there exists an $\alpha' \in \mathcal{A}$ such that $x_{\alpha'} \equiv x_{\alpha T}$, where $x_{\alpha T}$ is defined by $x_{\alpha T}(t) := x_\alpha(t+T)$, and $\lambda_{\alpha'} = \lambda_\alpha - T$.

- There exists a nonnegative, increasing function $\nu_3$, such that
$$|\dot{x}_\alpha(t)| \leq \nu_3(|x_\alpha(t)|) \ \forall \alpha \in \mathcal{A}, \quad \text{for almost all } t \in \mathcal{I}_\alpha$$

- There exist funcions $\kappa$ and $\chi$ of class $\mathcal{K}_\infty$ such that
$$\kappa(|x_\alpha(0)|) \geq \int_0^t \chi(|x_\alpha(s)|)\,ds \quad \forall \alpha \in \mathcal{A},\ t \in \mathcal{I}_\alpha.$$

Then for any two positive numbers $r$ and $\varepsilon$ there exists a $T_{r,\varepsilon}$, such that for all $\alpha \in \mathcal{A}$ and $t \in \mathcal{I}_\alpha$ the following holds:
$$t \geq T_{r,\varepsilon} \text{ and } |x_\alpha(0)| \leq r \Rightarrow |x_\alpha(t)| < \varepsilon.$$

*Proof.*

*Claim 1:* Given any $\varepsilon > 0$, there exists $\delta = \delta(\varepsilon)$ such that if $|x_\alpha(0)| \leq \delta$, then $|x_\alpha(t)| < \varepsilon$ for all $t \in I_\alpha$.

*Proof of Claim 1:* Fix a positive $\varepsilon$, and set
$$\delta(\varepsilon) = \min\left\{\frac{\varepsilon}{2}, \kappa^{-1}\left(\frac{\varepsilon\chi(\varepsilon/2)}{2\nu_3(\varepsilon)}\right)\right\}.$$

Pick any $\alpha \in \mathcal{A}$ such that $|x_\alpha(0)| \leq \delta$.

Suppose $|x_\alpha(\widetilde{t}_2)| \geq \varepsilon$ for some $\widetilde{t}_2 \in \mathcal{I}_\alpha$. Then there exist $t_1$ and $t_2$ with $t_1 < t_2 \leq \widetilde{t}_2$ such that $|x_\alpha(t_1)| = \varepsilon/2$ and $\varepsilon/2 < |x_\alpha(t)| < \varepsilon$ for all $t \in (t_1, t_2)$. Then
$$\frac{\varepsilon}{2} = |x_\alpha(t_2)| - |x_\alpha(t_1)| \leq |x_\alpha(t_2) - x_\alpha(t_1)| \leq \sup_{t_1 \leq t \leq t_2} |\dot{x}_\alpha(t)|(t_2 - t_1) \leq \nu_3(\varepsilon)(t_2 - t_1).$$

So,
$$\kappa(\delta) \geq \kappa(|\xi|) \geq \int_{t_1}^{t_2} \chi(|x_\alpha(s)|)\,ds > \frac{1}{2}\chi\left(\frac{\varepsilon}{2}\right)(t_2 - t_1) \geq \frac{\varepsilon\chi(\varepsilon/2)}{2\nu_3(\varepsilon)} \geq \kappa(\delta).$$

The obtained contradiction proves the claim.

*Claim 2:* Given positive numbers $r$ and $\delta$, there exists a time $\tau(r, \delta)$ such that if $|x_\alpha(0)| \leq r$ and $\tau(r, \delta) \in \mathcal{I}_\alpha$, then $\exists t_0 < \tau(r, \delta)$ such that $|x(t_0, \xi, \mathbf{d})| \leq \delta$.

*Proof of Claim 2:* Take $\tau(r,\delta) = \frac{2\kappa(r)}{\chi(\delta)}$. Then, if $\tau(r,\delta) \in I_\alpha$ and $|x_\alpha(t)| > \delta$ for all $t \in [0, \tau(r,\delta))$, then we have

$$\kappa(r) \geq \kappa(|x_\alpha(0)|) \geq \int_0^{\tau(r,\delta)} \chi(|x_\alpha(s)|)\, ds > \chi(\delta)\frac{\kappa(r)}{\chi(\delta)} = \kappa(r).$$

The obtained contradiction proves the claim.

Fix arbitrary positive $r$ and $\varepsilon$. By Claim 1, find $\delta(\varepsilon)$ such that if $|x_\alpha(0)| < \delta(\varepsilon)$ then $|x_\alpha(t)| \leq \varepsilon$ for all $t \in [0, \lambda_\alpha)$. Define $T_{r,\varepsilon} := \tau(r, \delta(\varepsilon))$, where $\tau(r, \delta(\varepsilon))$ is furnished by Claim 2. Then, if $|x_\alpha(0)| < r$ then, by Claim 2, there is a $t_0 < \tau(r, \delta(\varepsilon))$ with $|x_\alpha(t_0)| < \delta(\varepsilon)$. Consider now a function $x_{\alpha'}(\cdot) := x_\alpha(t_0 + \cdot)$ (it belongs to $\mathcal{X}$ by assumption). Since $|x_{\alpha'}| \leq \delta(\varepsilon)$, Claim 1 ensures that

$$|x_\alpha(t)| = |x_{\alpha'}(t - t_0)| \leq \varepsilon \quad \forall\, t \geq t_0 \geq T_{r,\varepsilon}.$$

This show that $T_{r,\varepsilon}$ satisfies the conclusion of the proposition. ∎

We now return to the proof of Lemma 3.8.

*Proof.* Recall that we have defined $\lambda_{\xi,\mathbf{d}} := \inf\{t \in [0, t_{\max}(\xi, \mathbf{d})) : |x(t, \xi, \mathbf{d})| \leq \rho(|y(t, \xi, \mathbf{d})|)\}$, and let $\lambda_{\xi,\mathbf{d}} = t_{\max}$ if $|x(t, \xi, d)| > \rho(|y(t, \xi, \mathbf{d})|)$ for all $t \in [0, t_{\max})$. Note that, given $\xi$ and $\mathbf{d}$, for all $t < \lambda_{\xi,\mathbf{d}}$ we have $\chi(|x(t, \xi, \mathbf{d})|) > 2\gamma(|h(x(t, \xi, \mathbf{d}))|)$ so that

$$\kappa(|\xi|) \geq \int_0^t \left(\chi(|x(s,\xi,\mathbf{d})|) - \gamma(|h(x(s,\xi,\mathbf{d}))|)\right) ds > \frac{1}{2}\int_0^t \chi(|x(s,\xi,\mathbf{d})|)\, ds. \tag{35}$$

Let $\nu_3$ be a $\mathcal{K}$-function such that $\max_{d \in \Omega} |f(x,d)| \leq \nu_3(|x|)$. Write $x_{\xi,\mathbf{d}}(\cdot) := x(\cdot, \xi, \mathbf{d})$. Notice that the family $\{x_{\xi,\mathbf{d}}(\cdot),\ \xi \in \mathbb{X},\ \mathbf{d} \in \mathcal{M}_\Omega\}$ with $\mathcal{I}_{\xi,\mathbf{d}} := [0, \lambda_{\xi,\mathbf{d}})$ satisfies all the assumptions of Proposition 3.9 (with "$\kappa$" $= 2\kappa$). Given any positive $r, \varepsilon$, Proposition 3.9 furnishes $T_{r,\varepsilon}$. This $T_{r,\varepsilon}$ obviously fits the first condition in the characterization of the GASMO property, provided by Proposition 2.11.

To find a function $\vartheta$ to ensure that the second part of Proposition 2.11 is satisfied, recall that, by Proposition 2.16, if a system (17) has the UO-property, then there exist class $\mathcal{K}_\infty$ functions $\rho_1$, $\mu_1$, $\mu_2$ and a constant $c > 0$ such that the following implication holds for all $\xi \in \mathbb{X}$, all $\mathbf{d} \in \mathcal{M}_\Omega$ and all $T \in [0, t_{\max}(\xi, \mathbf{d}))$:

$$|h(x(t,\xi,\mathbf{d}))| \leq \rho_1(|x(t,\xi,\mathbf{d})|)\ \forall t \in [0,T] \Rightarrow |x(t,\xi,\mathbf{d})| \leq \mu_1(t) + \mu_2(|\xi|) + c\ \forall t \in [0,T].$$

Therefore, if $|\xi| \leq r$ and $T_{r,r/2}$ is as defined above, then for all $t \in [0, \lambda_{\xi,\mathbf{d}})$ we have

$$|x(t,\xi,\mathbf{d})| \leq \mu_1(T_{r,r/2}) + \mu_2(r) + c, \quad \text{if } t < T_{r,r/2}$$

and

$$|x(t,\xi,\mathbf{d})| \leq r/2, \quad \text{if } t \geq T_{r,r/2}.$$

Thus, the following estimate holds for all such $t$:

$$|x(t,\xi,\mathbf{d})| \leq \widetilde{\vartheta}(|\xi|) := \max\left\{|\xi|/2,\ \mu_1(T_{|\xi|,|\xi|/2}) + \mu_2(|\xi|) + c\right\}.$$

Next, take a sequence $\{\varepsilon_k\}$, $k = 0, 1, 2...$, strictly decreasing to 0, with $\varepsilon_0 = 1$. For each $\varepsilon_k$, find $\delta_k = \delta(\varepsilon_k)$ as in the proof of Claim 1. Since $\delta_k \leq \varepsilon_k/2$, the sequence $\{\delta_k\}$ coverges to 0 as well. Find a function $\vartheta$ of class $\mathcal{K}$, such that:

1) $\vartheta(\delta_{k+1}) > \varepsilon_k \quad \forall k > 0$.

(this will ensure that $|x(t, \xi, \mathbf{d})| \leq \vartheta(|\xi|) \; \forall \xi$ with $|\xi| < \delta_0, \; \forall t \in [0, \lambda_{\xi, \mathbf{d}}]$.)

2) $\vartheta(s) \geq \widetilde{\vartheta}(s) \; \forall s > \delta_0$.

Then $\vartheta$ satisfies the second condition in the Proposition 2.11. This completes the proof. ∎

**Remark 3.10** The unboundedness observability assumption is crucial in proving the last lemma. The following example illustrates a disturbance-free iiOSS system which fails to be OSS (and, equivalently, fails to be GASMO).

Let $1_A(\cdot)$ denote the indicator function of a set $A$, and $\phi_\varepsilon$ be a $C^\infty$-bump function with support in $(-\varepsilon, \varepsilon)$:

$$\phi_\varepsilon(\xi) := \begin{cases} e^{-\frac{|\xi|^2}{\varepsilon^2 - |\xi|^2}}, & |\xi| < \varepsilon \\ 0, & |\xi| \geq \varepsilon. \end{cases} \tag{36}$$

Fix an arbitrary positive $\varepsilon < 0.25$ and consider a one dimensional autonomous system

$$\Sigma : \quad \dot{x} = f(x), \quad y = h(x)$$

where

$$\begin{aligned} f(x) &= x^3 \left[ 1_{(-\infty, -1]}(x)(1 - \phi_\varepsilon(x+1)) + 1_{[1,+\infty)}(x)(1 - \phi_\varepsilon(x-1)) \right] - \\ &\quad - x \left[ 1_{(-1,1)}(x)(1 - \phi_\varepsilon(x+1))(1 - \phi_\varepsilon(x-1)) \right], \end{aligned}$$

and $h$ is a smooth function such that $h(x) = x$ for all $x$ in $[-2, 2]$, and $h(x) = 0$ if $|x| \geq 3$.

*Claim:* The system $\Sigma$ is iiOSS.

*Proof:* Note that $\Sigma$ has a stable equilibrium at $x = 0$ and two unstable ones at $1$ and $-1$. If $|x| < 1$, then $\text{sign}(x) = -\text{sign}(f(x))$, so, if $|\xi| \leq 1$, then $|x(t, \xi)| \leq 1$ for any nonnegative $t$. Therefore, if $\xi \in [-1, 1]$, then for all $t \geq 0$ we have

$$\int_0^t |x(s, \xi)| \, ds = \int_0^t |h(x(s, \xi))| \, ds, \tag{37}$$

so, estimate (22) trivially follows for all $\xi \in [-1, 1]$ and all $t \in t_{\max}(\xi)$ with $\gamma = Id$ and any $\kappa \in \mathcal{K}$.

If $|\xi| \geq 1 + \varepsilon$, then $f(x) = x^3$, so that

$$x(t, \xi) = \frac{\text{sign}(\xi)}{\sqrt{\xi^{-2} - 2t}}.$$

Thus, in this case the solution $x(t, \xi)$ is defined for all nonnegative $t < t_{\max}(\xi) = \xi^{-2}/2$ and

$$\int_0^t |x(s, \xi)| \, ds \leq \int_0^{t_{\max}(\xi)} |x(s, \xi)| \, ds = \frac{1}{\xi} \leq \frac{1}{1 + \varepsilon}.$$

Let $\kappa$ be any $\mathcal{K}$-function such that $\kappa(1) \geq (1 + \varepsilon)^{-1}$. Suppose $1 < |\xi| < 1 + \varepsilon$. Let $\hat{t}$ be the time when $|x(\hat{t}, \xi)| = 1 + \varepsilon$. Then $t_{\max}(\xi) = \hat{t} + (1 + \varepsilon)^{-2}/2$. Also, $x(s, \xi) = h(x(s, \xi))$ for all $s \in [0, \hat{t}]$, so, in particular, equality (37) holds for all $t < \hat{t}$, which, again, trivially implies (22) with $\gamma = Id$ and any $\kappa \in \mathcal{K}$.

If $t > \hat{t}$, then

$$\begin{aligned}
\int_0^t |x(s,\xi)|\, ds &= \int_0^{\hat{t}} |x(s,\xi)|\, ds + \int_{\hat{t}}^t |x(s,\xi)|\, ds \\
&\leq \int_0^{\hat{t}} |x(s,\xi)|\, ds + \int_{\hat{t}}^{t_{\max}(\xi)} |x(s,\xi)|\, ds \\
&= \int_0^{\hat{t}} |h(x(s,\xi))|\, ds + \int_0^{t_{\max}(1+\varepsilon)} |x(s, 1+\varepsilon)|\, ds \\
&\leq \int_0^{\hat{t}} |h(x(s,\xi))|\, ds + \kappa(\xi) \\
&\leq \int_0^t |h(x(s,\xi))|\, ds + \kappa(\xi).
\end{aligned}$$

This shows that $\Sigma$ is iiOSS, as estimate (22) holds for $\Sigma$ with the $\kappa$ that we constructed and $\gamma = Id$.

*Claim 2.* System $\Sigma$ is not OSS.

Indeed, pick any initial state $\xi$ of large enough magnitude so that $h(\xi) = 0$. Then $h(x(t,\xi)) = 0$ for all $t < t_{\max}(\xi) = \xi^{-2}/2$. If $\Sigma$ were OSS, then there would exist some $\mathcal{KL}$-function $\beta$ such that $|x(t,\xi)| \leq \beta(|\xi|, t)$, but $|x(t,\xi)|$ tends to $\infty$ as $t \to t_{\max}(\xi)$, whereas $\beta(|\xi|, t) \leq \beta(|\xi|, 0)$. This contradiction proves the claim. □

We now prove implication $4 \Rightarrow 3$ of Theorem 2.

*Proof.* Suppose a system $\Sigma$ of type (17) is UOSS. We have already remarked that $\Sigma$ is UO, so we must show it is iiUOSS. By assumption there exists a smooth function $V$ satisfying (19) and (7) with some $\alpha_1$, $\alpha_2$, $\alpha_3$, and $\gamma$. Pick any $\xi$, $\mathbf{d}$ and $t \in [0, t_{\max}(\xi, \mathbf{d}))$. Integrating inequality (19) along the trajectory $x(\cdot, \xi, \mathbf{d})$ over $[0, t]$ we get

$$\begin{aligned}
\int_0^t \alpha_3(|x(t,\xi,\mathbf{d})|)dt &\leq V(x(0,\xi,\mathbf{d})) - V(x(t,\xi,\mathbf{d})) + \int_0^t \gamma(|h(x(t,\xi,\mathbf{d}))|)\, dt \\
&\leq \alpha_2(|\xi|) + \int_0^t \gamma(|h(x(t,\xi,\mathbf{d}))|)\, dt,
\end{aligned}$$

proving inequality (22) for system $\Sigma$, with $\chi = \alpha_3$ and $\kappa = \alpha_2$. ∎

With Lemma 3.7 in mind we conclude that the only step missing in establishing the Lyapunov characterization for UOSS is proving the implication $2 \Rightarrow 4$ in Theorem 2.

## 4 The case of no controls

### 4.1 Setup

Suppose a system $\Sigma$ of type (17) satisfies the GASMO property with some $\mathcal{K}_\infty$ function $\rho$. As discussed in 2.6.1, $\rho$ can be assumed to be smooth when restricted to $\mathbb{R}_{>0}$ and also $\rho(s) > s$ for all positive $s$. Recall the following notation, introduced in 2.6.1:

- $\mathcal{D} := \{\xi \in \mathbb{X}: |\xi| \leq \rho(|h(\xi)|)\}$,

- $\mathcal{E} := \mathbb{X} \setminus \mathcal{D}$, and
- $\mathcal{E}_1 := \{\xi \in \mathbb{X} : |\xi| > 2\rho(|h(\xi)|)\}$.

If $\mathcal{D} = \mathbb{X}$, then any proper, smooth and positive definite function $V : \mathbb{X} \to \mathbb{R}$ is a UOSS-Lyapunov function for (17). Indeed, because it is proper and finite, $V$ obviously satisfies (7) for some $\alpha_1$ and $\alpha_2$. Since $V$ is smooth, $|\nabla V(\xi)|$ is bounded above by a nondecreasing continuous function $\nu(|\xi|)$ and

$$\frac{d}{dt} V(x(t)) = \nabla V(x(t)) \cdot f(x(t), \mathbf{d}(t)) \leq \nu(|x(t)|)\nu_3(|x(t)|),$$

where $\nu_3(|\cdot|)$ is a $\mathcal{K}$-function majorizing $f(\cdot, v)$ for all $v \in \Omega$. Then, since $|x| \leq \rho(|h(x)|)$ for all $x \in \mathbb{X}$, we have

$$\frac{d}{dt} V(x(t)) \leq -\nu(|x(t)|)\nu_3(|x(t)|) + 2\nu(\rho(|h(x(t))|))\nu_3(\rho(|h(x(t))|)).$$

So, $V$ satisfies inequality

$$\nabla V(x) \cdot f(x, d) \leq -\alpha_3(|x|) + \gamma(|h(x)|) \quad \forall\, x \in \mathbb{X},\ \forall\, d \in \Omega,$$

(with $\alpha_3(\cdot) = \nu(\cdot)\nu_3(\cdot)$ and $\gamma(\cdot) = [2\nu \circ \rho(\cdot)][\nu_3 \circ \rho(\cdot)]$) which is the same as (8) for systems of type (17).

Suppose now that $\mathcal{D} \neq \mathbb{X}$. Recall that we have defined, for each $\xi \notin \mathcal{D}$ and $d \in \mathcal{M}_\Omega$, $\lambda_{\xi,\mathbf{d}} = \inf \{t \in [0, t_{\max}) : x(t, \xi, \mathbf{d}) \in \mathcal{D}\}$, with the convention $\lambda_{\xi,\mathbf{d}} = t_{\max}(\xi, \mathbf{d})$ if the trajectory never enters $\mathcal{D}$.

The GASMO property then implies

$$|x(t, \xi, \mathbf{d})| \leq \lambda(|\xi|, t), \qquad \forall\, \xi \in \mathcal{E},\ \forall\, \mathbf{d} \in \mathcal{M}_\Omega,\ \forall\, t \in [0, \lambda_{\xi,\mathbf{d}}) \tag{38}$$

for some $\lambda \in \mathcal{KL}$.

Note that, because of property (38), the system cannot have any equilibrium in $\mathcal{E}$, that is,

$$f(\xi, d) \neq 0$$

for every $\xi \in \mathcal{E}$ and every $d \in \Omega$. Moreover, replacing $\rho(s)$ by $c\rho(s)$ for some $c > 1$ if necessary, one may also assume that $f(\xi, d) \neq 0$ for all $\xi \in \partial \mathcal{D} \setminus \{0\}$, all $d \in \Omega$.

We introduce an auxiliary system $\widehat{\Sigma}$ which slows down the motions of the original one:

$$\dot{z} = \widehat{f}(z, d) = \frac{1}{1 + |f(z,d)|^2 + \kappa(z)} f(z, d) \tag{39}$$

where $\kappa$ is any smooth function $\mathbb{X} \to [0, \infty)$ with the property that

$$\kappa(\xi) \geq 2 \max_{d \in \Omega} |\nabla(\rho \circ |h|)(\xi) \cdot f(\xi, d)| \tag{40}$$

whenever $|h(\xi)| \geq 1$. (Recall that $\rho$ was assumed, without loss of generality, to be smooth for positive arguments.) For each disturbance $\hat{\mathbf{d}}$ (defined on $\mathbb{R}_{\geq 0}$) denote by

$$z(s, \xi, \hat{\mathbf{d}})$$

the value at time $s$ of the solution of the equation $\dot{z} = \widehat{f}(z, \hat{\mathbf{d}})$ with initial state $\xi$. Observe that, as $\widehat{f}$ is bounded, this solution exists for all nonnegative $s$.

*Claim 1 :* For each $\xi$ and each $\mathbf{d}$,

$$x(t, \xi, \mathbf{d}) = z(\sigma_{\xi,\mathbf{d}}(t), \xi, \mathbf{d} \circ \sigma_{\xi,\mathbf{d}}^{-1}) \qquad \forall t \in [0, t_{\max}(\xi, \mathbf{d})), \tag{41}$$

where $\sigma_{\xi,\mathbf{d}} : [0, t_{\max}(\xi, \mathbf{d})) \to \mathbb{R}_{\geq 0}$ is defined by

$$\sigma_{\xi,\mathbf{d}}(t) = \int_0^t \left[ 1 + |f(x(s, \xi, \mathbf{d}), \mathbf{d}(s))|^2 + \kappa(x(s, \xi, \mathbf{d})) \right] ds.$$

Moreover, $\sigma_{\xi,\mathbf{d}}(t) \to \infty$ as $t \to t_{\max}(\xi, \mathbf{d})$, so, we can define $\sigma_{\xi,\mathbf{d}}(t_{\max}(\xi, \mathbf{d})) := +\infty$ for convenience.

*Proof of Claim 1.* Indeed, writing $s = \sigma_{\xi,\mathbf{d}}(t)$ and computing the derivative of $x(\sigma_{\xi,\mathbf{d}}^{-1}(s), \xi, \mathbf{d})$ with respect to $s$, one has:

$$\begin{aligned}
f(x(t, \xi, \mathbf{d}), \mathbf{d}(t)) &= \frac{d}{dt} x(t, \xi, \mathbf{d}) = \frac{d}{dt} x(\sigma_{\xi,\mathbf{d}}^{-1} \circ \sigma_{\xi,\mathbf{d}}(t), \xi, \mathbf{d}) \\
&= \frac{d}{ds} x\left(\sigma_{\xi,\mathbf{d}}^{-1}(s), \xi, \mathbf{d}\right) \cdot \frac{d}{dt} \sigma_{\xi,\mathbf{d}}(t) \\
&= \frac{d}{ds} x(\sigma_{\xi,\mathbf{d}}^{-1}(s), \xi, \mathbf{d}) \left[ 1 + |f(x(t, \xi, \mathbf{d}), \mathbf{d}(t))|^2 + \kappa(x(t, \xi, \mathbf{d})) \right].
\end{aligned}$$

Therefore

$$\begin{aligned}
\frac{d}{ds} x(\sigma_{\xi,\mathbf{d}}^{-1}(s), \xi, \mathbf{d}) &= \frac{f(x(t, \xi, \mathbf{d}), \mathbf{d}(t))}{1 + |f(x(t, \xi, \mathbf{d}), \mathbf{d}(t))|^2 + \kappa(x(t, \xi, \mathbf{d}))} \\
&= \frac{f\left(x(\sigma_{\xi,\mathbf{d}}^{-1}(s), \xi, \mathbf{d}), \mathbf{d} \circ \sigma_{\xi,\mathbf{d}}^{-1}(s)\right)}{1 + |f(x(\sigma_{\xi,\mathbf{d}}^{-1}(s), \xi, \mathbf{d}), \mathbf{d} \circ \sigma_{\xi,\mathbf{d}}^{-1}(s))|^2 + \kappa(x(\sigma_{\xi,\mathbf{d}}^{-1}(s), \xi, \mathbf{d}))} \\
&= \widehat{f}\left(x\left(\sigma_{\xi,\mathbf{d}}^{-1}(s), \xi, d\right), \mathbf{d} \circ \sigma_{\xi,\mathbf{d}}^{-1}(s)\right),
\end{aligned}$$

for all $0 \leq s < \sigma_{\xi,\mathbf{d}}(t_{\max}(\xi, \mathbf{d}))$. Thus, the functions $z(s, \xi, \mathbf{d} \circ \sigma_{\xi,\mathbf{d}}^{-1})$ and $x(\sigma_{\xi,\mathbf{d}}^{-1}(s), \xi, \mathbf{d})$ satisfy the same differential equation (39) with initial state $\xi$ on $[0, \sigma_{\xi,\mathbf{d}}(t_{\max}(\xi, \mathbf{d})))$, therefore they coincide on $[0, \sigma_{\xi,\mathbf{d}}(t_{\max}(\xi, \mathbf{d})))$.

To show the limit property of $\sigma_{\xi,\mathbf{d}}$, suppose

$$\lim_{t \to t_{\max}(\xi, \mathbf{d})} \sigma_{\xi,\mathbf{d}}(t) = b < \infty$$

(note that the limit exists because $\sigma_{\xi,\mathbf{d}}(\cdot)$ is increasing). Let

$$K = \{z(s, \xi, \mathbf{d} \circ \sigma_{\xi,\mathbf{d}}^{-1}) : 0 \leq s < b\}.$$

Then $K$ is bounded, so $\bar{K}$ is compact, and by (41), $x(t, \xi, \mathbf{d}) \in K \subseteq \bar{K}$ for all $t \in [0, t_{\max}(\xi, \mathbf{d}))$, contradicting the maximality of $t_{\max}(\xi, \mathbf{d})$. ∎

For each initial state $\xi$ and each disturbance function $\mathbf{d}$, define

$$\theta_{\xi,\mathbf{d}} = \inf \{t \geq 0 : z(t, \xi, \mathbf{d}) \in \mathcal{D}\}, \tag{42}$$

where $\theta_{\xi,\mathbf{d}} = \infty$ if $z(t,\xi,\mathbf{d}) \notin \mathcal{D}$ for all $t \geq 0$. Note that $\theta_{\xi,\mathbf{d}} > 0$ for all $\mathbf{d} \in \mathcal{M}_\Omega$ and all $\xi \in \mathcal{E}$, because $\mathcal{E}$ is open. Observe also that if $\mathbf{d}_1 = \mathbf{d} \circ \sigma_{\xi,\mathbf{d}}$,

$$\theta_{\xi,\mathbf{d}} = \sigma_{\xi,\mathbf{d}_1}(\lambda_{\xi,\mathbf{d}_1}).$$

*Claim 2:* System $\hat{\Sigma}$ satisfies the GASMO property.

*Proof:* According to (38) and (41), we have, for every $\xi \in \mathcal{E}$ and each $\mathbf{d}$,

$$\begin{aligned} |z(t,\xi,\mathbf{d})| &= \left|x(\sigma_{\xi,\mathbf{d}_1}^{-1}(t),\xi,\mathbf{d}_1)\right| \\ &\leq \lambda(|\xi|, \sigma_{\xi,\mathbf{d}_1}^{-1}(t)) \leq \vartheta(|\xi|), \end{aligned}$$

for all $t \in [0, \theta_{\xi,\mathbf{d}})$, where $\vartheta(s) = \lambda(s,0)$, and $\mathbf{d}_1 = \mathbf{d} \circ \sigma_{\xi,\mathbf{d}}$. Let

$$M_r = 1 + \max_{d \in \Omega, |\xi| \leq \vartheta(r)} |f(\xi,d)|^2 + \max_{|\xi| \leq \vartheta(r)} \kappa(\xi).$$

Then, for any $\xi \in \mathcal{E}$ with $|\xi| \leq r$, it holds that

$$\sigma_{\xi,\mathbf{d}_1}(t) = \int_0^t \left[1 + |f(x(s,\xi,\mathbf{d}_1),\mathbf{d}_1(s))|^2 + \kappa(x(s,\xi,\mathbf{d}_1))\right] ds \leq M_r t$$

for all $t \in [0, \lambda_{\xi,\mathbf{d}_1})$, and hence, $\sigma_{\xi,\mathbf{d}_1}^{-1}(t) \geq \frac{t}{M_r}$ for all $|\xi| \leq r$, $t \in [0, \theta_{\xi,\mathbf{d}})$. Consequently, we have

$$|z(t,\xi,\mathbf{d})| \leq \hat{\lambda}(|\xi|,t) \qquad \forall t \in [0, \theta_{\xi,\mathbf{d}}), \tag{43}$$

where $\hat{\lambda}(s,t) = \lambda(s, \frac{t}{M_s})$ is clearly of class $\mathcal{KL}$. This shows that system $\hat{\Sigma}$ is GASMO. ∎

From now on, we let the function $\lambda$ of class $\mathcal{KL}$ be as in definition 2.9 for the system $\hat{\Sigma}$, that is, the following estimate holds for system (39):

$$|z(t,\xi,\mathbf{d})| \leq \lambda(|\xi|,t) \qquad \forall t \in [0, \theta_{\xi,\mathbf{d}}), \tag{44}$$

for all $\xi \in \mathcal{E}$ and all $\mathbf{d} \in \mathcal{M}_\Omega$.

According to Proposition 7 in [39], there exist $\mathcal{K}_\infty$-functions $\mu_1$ and $\mu_2$ such that

$$\lambda(r,t) \leq \mu_1(\mu_2(r)e^{-t}) \qquad \forall r, t \geq 0. \tag{45}$$

Define

$$\Xi(s) := \mu_1^{-1}(s).$$

The proof will now develop as follows. We first construct a continuous Lyapunov-like function $V_0$, defined on the set $\mathcal{E}_1$. Next $V_0$ is approximated by a Lipschitz continuous function (by the methods of nonsmooth analysis). The resulting function is then approximated by a smooth function $V_1$. Finally we extend $V_1$ to the rest of the state space, obtaining a Lyapunov-like function, smooth away from the origin, which is then approximated by a smooth Lyapunov function.

Before we start, we prove a technical lemma.

**Lemma 4.1** Suppose $\Sigma : \dot{z} = f(z,d)$, $y = h(z)$ is a system of type (17), and $p(\cdot)$ is a smooth function of class $\mathcal{K}_\infty$, such that the following conditions hold:

- $|f(\xi, d)| \leq 1$ for all $\xi \in \mathbb{X}$ and $d \in \Omega$,
- $|\nabla(p \circ |h|)(\xi) \cdot f(\xi, d)| \leq 1$ for all $d \in \Omega$ and all $\xi$ with $|h(\xi)| \geq 1$,
- $p(s) \geq s$ for all $s > 0$.

Pick any constant $a > 0$ and define $K_0 = p(1) + (2 + a)/a + 1$. Then for each $\xi \in \mathbb{X}$ such that

$$|\xi| \geq (1 + a)p(|h(\xi)|) \text{ and } |\xi| \geq K_0,$$

it holds that

$$|z(t, \xi, \mathbf{d})| > p(|h(z(t, \xi, \mathbf{d}))|)$$

for all $t \in [0, 1)$ and any $\mathbf{d} \in \mathcal{M}_\Omega$.

*Proof.* Fix $a$, $\xi$ and $\mathbf{d}$ as in the formulation of the lemma, and define

$$\theta := \min\{t : |z(t, \xi, \mathbf{d})| \leq p(|h(z(t, \xi, \mathbf{d}))|)\},$$

with the convention $\theta = +\infty$ if the inequality never holds for $t \geq 0$. Assume the lemma is false, so that $\theta < 1$.

Let $\eta = z(\theta, \xi, \mathbf{d})$, and let $\hat{\mathbf{d}}$ be the shift of $\mathbf{d}$ by $\theta$, that is $\hat{\mathbf{d}}(t) = \mathbf{d}(t + \theta)$. Since $|f(z, d)| \leq 1$ for all $z \in \mathbb{X}$ and all $d \in \Omega$, it holds that $|\eta| \geq |\xi| - \theta \geq K_0 - 1$. By the definitions of $\eta$ and $\theta$, one has

$$p(|h(\eta)|) = |\eta| \geq K_0 - 1, \tag{46}$$

so also $|h(\eta)| \geq p^{-1}(K_0 - 1) > 1$. Thus, $\left|h(z(s, \eta, \hat{\mathbf{d}}))\right| > 1$ for all $s$ near zero.

*Claim:* $\left|h(z(s, \eta, \hat{\mathbf{d}}))\right| > 1$ for all $s \in [-1, 0]$.

Assume the claim is false. Then there must exist some $-1 \leq s_0 < 0$ so that

$$s_0 = \max\left\{s \leq 0 : \left|h(z(s, \eta, \hat{\mathbf{d}}))\right| \leq 1\right\}.$$

We have that for each $s \in (s_0, 0]$, $\left|h(z(s, \eta, \hat{\mathbf{d}}))\right| > 1$.

Recall that $|\nabla(p \circ |h|)(z)f(z, d)| \leq 1$ for all $z$ with $|h(z)| \geq 1$ and all $d \in \Omega$. Thus

$$\left|\frac{d}{ds}p\left(\left|h(z(s, \eta, \hat{\mathbf{d}}))\right|\right)\right| \leq 1 \quad \forall\, s \in (s_0, 0].$$

This, in turn, implies that

$$p\left(\left|h(z(s_0, \eta, \hat{\mathbf{d}}))\right|\right) \geq p(|h(\eta)|) + s_0 \geq K_0 + s_0 - 1 > p(1),$$

and so, since $p$ is strictly increasing, $\left|h(z(s_0, \eta, \hat{\mathbf{d}}))\right| > 1$, thus contradicting the definition of $s_0$. This proves the claim.

It follows from the claim that $|h(z(s, \xi, \mathbf{d}))| > 1$ for all $s \in [0, \theta]$. Thus,

$$\begin{aligned} p(|h(\eta)|) &= |\eta| \geq |\xi| - \theta \\ &\geq (1 + a)p(|h(\xi)|) - \theta \geq (1 + a)p(|h(\eta)|) - (1 + a)\theta - \theta. \end{aligned}$$

(The last inequality used the fact that $\left|\frac{d}{ds}p(|h(z(s,\xi,\mathbf{d})|)\right| \leq 1$ for all $s \in [0,\theta]$). It follows that $p(|h(\eta)|) \leq \frac{2+a}{a}\theta$, so from (46) we know that

$$K_0 \leq 1 + p(|h(\eta)|) \leq 1 + \frac{2+a}{a}\theta,$$

contradicting the choice of $K_0$. This shows that it is impossible to have $\theta < 1$. ∎

## 4.2 Definitions and basic facts on relaxed controls

Recall that our disturbances $\mathbf{d}$ are measurable functions $\mathbb{R}_{\geq 0} \to \Omega = [-1,1]^m$.

Let $\mathcal{P}(\Omega)$ be the set of all Radon probability measures on $\Omega$. Bishop's theorem furnishes a weak norm on $\mathcal{P}(\Omega)$, whose corresponding metric topology coincides with the weak star topology on $\mathcal{P}(\Omega)$ (see [51], pages 40 and 267).

For any $T > 0$, we define $\mathcal{S}_T$ to be the set of all measurable functions from $[0,T]$ to $\mathcal{P}(\Omega)$, and $\mathcal{S}$ to be the set of all measurable functions from $\mathbb{R}_{\geq 0}$ to $\mathcal{P}(\Omega)$. We topologize $\mathcal{S}_T$ by weak convergence: $\{\nu_k(\cdot)\} \to \nu(\cdot)$ in $\mathcal{S}_T$ if and only if

$$\int_0^T \int_\Omega g(t,\omega) d[\nu_k(t)](\omega)\, dt \to \int_0^T \int_\Omega g(t,\omega) d[\nu(t)](\omega)\, dt$$

for all functions $g : [0,T] \times \Omega \to \mathbb{R}$ which are continuous in $\omega$, measurable in $t$ and such that

$$\max\{|g(t,\omega)|,\ \omega \in \Omega\}$$

is integrable on $[0,T]$. We say that $\{\nu_k\} \to \nu$ weakly in $\mathcal{S}$ if, for every $T > 0$, the sequence $\{\nu_k|_{[0,T]}\}$ of restrictions of $\nu_k$ to $[0,T]$ converges to $\nu|_{[0,T]}$ in $\mathcal{S}_T$.

Notice that, since every element of $\Omega$ can be identified with the $\delta$-measure, concentrated in it, $\Omega$ can be embedded into $\mathcal{P}(\Omega)$, $\mathcal{M}_\Omega$ into $\mathcal{S}$, and $\mathcal{M}_\Omega^T$ into $\mathcal{S}_T$ in the obvious way, where $\mathcal{M}_\Omega^T$ is the set of functions in $\mathcal{M}_\Omega$ restricted to $[0,T]$.

For each $\nu \in \mathcal{P}(\Omega)$, we denote

$$f(x,\nu) = \int_\Omega f(x,r)\, d\nu(r). \qquad (47)$$

Notice that for any relaxed control $\nu(\cdot)$, the function $f(\cdot,\nu(\cdot)) : (x,t) \to f(x,\nu(t))$ is Lipschitz in $x$ and measurable in $t$. Moreover, as for all $x \in \mathbb{X}$ and all $\nu \in \mathcal{P}(\Omega)$ we have a bound

$$|f(x,\nu)| \leq \max_{d \in \Omega}|f(x,d)|,$$

the solution of the "system" $\dot{x}(t) = f(x(t),\nu(t))$ exists for any initial condition $\xi$ and relaxed control $\nu$ on some maximal interval $[0, t_{\max}(\xi,\nu))$. Write $x(\cdot,\xi,\nu)$ to denote this solution. Just as in the case with ordinary controls, we will define

$$\lambda_{\xi,\nu} := \inf\{t \leq t_{\max} : x(s,\xi,\nu) \in \mathcal{D}\}.$$

The basic three facts we will be using in the next section are as follows:

*Fact 1.* For any $T > 0$, the space $\mathcal{S}_T$ is sequentially compact (See [51], Thm IV.2.1, page 272). Consequently, $\mathcal{S}$ is sequentially compact by a diagonalization argument.

*Fact 2.* For any $T > 0$, the set $\mathcal{M}_\Omega^T$ of ordinary controls on $[0, T]$ is dense in $\mathcal{S}_T$ (see [5], page 691, also [51]). Consequently, $\mathcal{M}_\Omega$ is dense in $\mathcal{S}$. The topology of $\mathcal{S}$ induces a topology on the subspace $\mathcal{M}_\Omega$. This topology is stronger than the topology of $L^p$ for any positive $p$: in fact, $\{\mathbf{d}_k(t)\} \to \mathbf{d}$ would imply

$$\int_0^T g(\mathbf{d}_k(t))dt \to \int_0^T g(\mathbf{d}(t))\,dt$$

for any continuous function $g : \mathbb{U} \to \mathbb{R}$ and any positive $T$.

*Fact 3.* For any $T > 0$, the mapping $(t, \xi, \nu(\cdot)) \mapsto z(t, \xi, \nu)$ is continuous on $[0, T] \times \mathbb{X} \times \mathcal{S}_T$ (see [1], Lemma 3.12).

Some relevant immediate consequences of Facts 1, 2, and 3 are as follows.

*Claim 1.* The function $(\xi, \nu) \to \lambda_{\xi,\nu}$ is lower semicontinuous on both $\nu$ and $\xi$. (This easily follows from the continuous dependence on the initial conditions and control, and from the fact that the set $\mathcal{D}$ is closed).

*Claim 2.* If system (17) is GASMO, then, for any initial value $\xi$, a relaxed disturbance $\nu \in \mathcal{S}$, and time $T < \lambda_{\xi,\nu}$ we have an estimate

$$x(T, \xi, \nu) \leq \lambda(|\xi|, T).$$

*Proof.* Pick $\xi \in \mathbb{X}$ and $\nu \in \mathcal{S}$. Let $\{\mathbf{d}_k\}$ be a sequence of ordinary disturbances, converging to $\nu$ in $\mathcal{S}$. Then, for large enough $k$ we have $t_{\max}(\xi, \mathbf{d}_k) > T$, and $x(\cdot, \xi, \mathbf{d}_k)$ converge to $x(\cdot, \xi, \nu)$ uniformly on $[0, T]$. Also, by Claim 1, $\lambda_{\xi, \mathbf{d}_k} \geq \lambda_{\xi, \nu}$. Since we assume the system to be GASMO, the estimate (20) holds for $\mathbf{d} := \mathbf{d}_k$ for all large enough $k$ and for all $t \leq T$, so, the Claim follows. ∎

With the previous claim in mind, we can say that the system (17) with $d \in \mathcal{S}$ is GASMO (this is a slight abuse of terminology because, strictly speaking, (17) with relaxed disturbances is not really a "system", as that would mean, by definition, that disturbances take values in a finite dimensional space). Consequently, the auxiliary system (39) is also GASMO for $\mathbf{d} \in \mathcal{S}$. Thus, we can assume that (44) holds for all $\xi \in \mathcal{E}$ and all $\mathbf{d} \in \mathcal{S}$.

### 4.3 Constructing a continuous Lyapunov-like function on $\mathcal{E}_1$

In the section 4.1 we introduced the system $\hat{\Sigma} : \dot{z} = \hat{f}(z, d) := \frac{f(z,d)}{1+|f(z,d)|^2+\kappa(z)}$, which slows down the motions of the original system $\Sigma$. Recall also that $\mathcal{D} := \{\xi : |\xi| \leq \rho(|h(\xi)|)\}$, and define the set

$$\mathcal{B} := \{\xi : \rho(|h(\xi)|) \leq |\xi| \leq 1.5\rho(|h(\xi)|)\}.$$

Let $f_0 : \mathbb{X} \to \mathbb{R}$ be defined by

$$f_0(\xi) = \max_{\mathrm{d} \in \Omega} \left|\widehat{f}(\xi, \mathrm{d})\right|.$$

Note that $f_0$ is locally Lipschitz, and recall that we have assumed with no loss of generality that $\Sigma$ has no equilibria on the set $\{x \in \mathbb{X} : |x| \geq \rho(|h(x)|)\}$ (otherwise replace $\rho(\cdot)$ by $c\rho(\cdot)$ where $c > 1$). In particular, $f_0(\xi) \neq 0$ for any $\xi \in \partial \mathcal{D} \setminus \{0\}$. Let $\phi : \mathbb{X} \setminus \{0\} \to [0, 1]$ be smooth and

$$\phi(x) = \begin{cases} 1 & x \in \mathcal{D} \\ 0 & x \in \mathbb{X} \setminus (\mathcal{D} \cup \mathcal{B}). \end{cases}$$

Now introduce another system, on the state space $\mathbb{X} \setminus \{0\}$,
$$\widetilde{\Sigma}: \quad \dot{z} = \widetilde{f}(z, d, v),$$
where disturbances $d$ are as before, and auxiliary controls $\mathbf{v}$ are measurable functions of time, taking values in $[-1, 1]^n$ (note that the dimension of the control set for $v$'s is the same as that of $\mathbb{X}$), where, for each $i = 1, 2, \ldots n$:
$$\widetilde{f}_i(z, d, v) := \widehat{f}_i(z, d) + 2\phi(z)f_0(z)v_i.$$

For each $T > 0$, let $\mathcal{W}_T$ denote the set of the auxiliary controls defined on $[0, T]$ (i.e., measurable functions $\mathbf{v} : [0, T] \to [-1, 1]^n$) equipped with the weak convergence topology; that is, "$\{\mathbf{v}_k\}$ converges weakly to $\mathbf{v}$ in $\mathcal{W}_T$" means
$$\int_0^T \varphi(s)\mathbf{v}_k(s)\,ds \to \int_0^T \varphi(s)\mathbf{v}(s)\,ds$$
for all functions $\varphi$ that are integrable over $[0, T]$. With the weak topology, $\mathcal{W}_T$ is sequentially compact (c.f. [38, Proposition 10.1.5]). Consequently, given any sequence $\{\mathbf{v}_k\}$ of controls defined on $[0, \infty)$, there exist some control $\mathbf{v}$ and a subsequence $\{\mathbf{v}_{k_j}\}$ such that $\mathbf{v}_k \to \mathbf{v}$ weakly on every interval $[0, T]$.

We denote the set of the auxiliary controls defined on $[0, \infty)$ by $\mathcal{W}$. Let $\{\mathbf{v}_k\} \subset \mathcal{W}$ and $\mathbf{v} \in \mathcal{W}$. We say that $\{\mathbf{v}_k\}$ weakly converges to $\mathbf{v}$ if, for every $T > 0$, the sequence of restrictions $\{\mathbf{v}_k|_{[0,T]}\}$ weakly converges to $\mathbf{v}|_{[0,T]}$ in $\mathcal{W}_T$.

Recall that $z(t, \xi, \mathbf{d})$ denotes the solution of $\widehat{\Sigma}$. Write $z(t, \xi, \mathbf{d}, \mathbf{v})$ for the value at time $t$ of the solution of $\widetilde{\Sigma}$ with initial state $\xi \neq 0$, disturbance $\mathbf{d} \in \mathcal{M}_\Omega$) and auxiliary control $\mathbf{v} \in \mathcal{W}$.

Observe that

- $\widetilde{\Sigma}$ is affine in $v$.

- Since for any $\xi \in \mathbb{X}$ and $d \in [-1, 1]^m$, $\left|\widehat{f}(\xi, d)\right| \leq 1$, we have $\left|\widetilde{f}(\xi, d, v)\right| \leq 3$ for any $\xi$, $d$ and $v$. In particular, this implies that $\widetilde{\Sigma}$ is forward complete.

- Suppose $\xi \notin \mathcal{D} \cup \mathcal{B}$ and pick $\mathbf{d} \in \mathcal{M}_\Omega$ and $\mathbf{v} \in \mathcal{W}$. Then there is some $t_0 > 0$ such that $z(t, \xi, \mathbf{d}) \notin \mathcal{D} \cup \mathcal{B}$ for all $t \in [0, t_0]$, and $z(t, \xi, \mathbf{d}, \mathbf{v}) \equiv z(t, \xi, \mathbf{d})$ on $[0, t_0]$.

To extend the definition of $z(t, \xi, \mathbf{d}, \mathbf{v})$ to the case when $\mathbf{d} \in \mathcal{S}$, we let, for a fixed $\mathbf{d} \in \mathcal{S}$ and $\mathbf{v} \in \mathcal{W}$,
$$g_{\mathbf{d},\mathbf{v}}(z, t) := \widehat{f}(z, \mathbf{d}(t)) + 2\phi(z)f_0(z)\mathbf{v}(t),$$
(where $f(z, \mathbf{d}(t))$ is as defined in (47) for $\nu := \mathbf{d}(t) \in \mathcal{P}(\Omega)$.) Then $g_{\mathbf{d},\mathbf{v}} : \mathbb{X} \setminus \{0\} \times \mathbb{R}_{\geq 0} \to \mathbb{X}$ is locally Lipschitz in its first variable, and $|g_{\mathbf{d},\mathbf{v}}(z, t)| \leq 3$ for all $(z, t) \in \mathbb{X} \setminus \{0\} \times \mathbb{R}_{\geq 0}$. Hence, the solution of
$$\begin{aligned}\dot{z}(t) &= g_{\mathbf{d},\mathbf{v}}(z, t) \\ z(0) &= \xi\end{aligned}$$
exists for all $\xi \in \mathbb{X} \setminus \{0\}$ and $t \geq 0$. We will denote it by $z(t, \xi, \mathbf{d}, \mathbf{v})$.

Observe that $\widetilde{f} : \mathbb{X} \setminus \{0\} \times \Omega \times [-1, 1]^n \to \mathbb{R}^n$ is continuous, and $\widetilde{f}(z, d, v)$ is locally Lipschitz in $z$ on $\mathbb{X} \setminus \{0\}$ uniformly on $(d, v) \in \Omega \times [-1, 1]^n$. The system $\widetilde{\Sigma}$ evolves in the state space $\mathbb{X} \setminus \{0\}$. As $\left|\widetilde{f}\right| \leq 3$ everywhere, trajectories are defined and unique for each initial value $\xi \in \mathbb{X} \setminus \{0\}$ and each pair of inputs $\mathbf{d}, \mathbf{v}$. Moreover, if $z(\cdot)$ is a maximal such trajectory, then either $z(t)$ is defined for all $t \geq 0$, or there is some $T > 0$ such that $\lim_{t \to T} z(t) = 0$. We prove next that this last case cannot happen.

**Lemma 4.2** For every ball $U$ around 0, there is a constant $c$ such that for any $\xi \in U$, $\xi \neq 0$, $\mathbf{d} \in \mathcal{S}$, $\mathbf{v} \in \mathcal{W}$, and $t \geq 0$, we have a lower bound

$$|z(t, \xi, \mathbf{d}, \mathbf{v})| > \frac{1}{2} |\xi|\, e^{-ct}. \tag{48}$$

*Proof.* Since $\hat{f}$ is locally Lipschitz in $x$ uniformly in $d$, and $\hat{f}(0, d) = 0$ for all $d$, we can find a positive constant $c$, such that $\left|\hat{f}(z, d)\right| \leq c|z|/3$ for all $z \in U$ and all $d \in \Omega$. Then also $|f_0(z)| \leq c|z|/3$; therefore

$$\left|\widetilde{f}(z, d, v)\right| \leq \left|\hat{f}(z, d)\right| + 2|f_0(z)| \leq c|z|$$

for all $v \in [-1, 1]^n$, all $d \in \Omega$ and, hence, all $d \in \mathcal{P}(\Omega)$. Now fix $\xi \in U \setminus \{0\}$, $\mathbf{d} \in \mathcal{S}$, and $\mathbf{v} \in \mathcal{W}$, and write $z(t) := z(t, \xi, \mathbf{d}, \mathbf{v})$. Note that the inequality (48) holds for $t = 0$, therefore it holds for all small enough $t > 0$. Suppose that (48) fails at some $t_2 > 0$, so that

$$|z(t_2)| \leq \frac{1}{2}|\xi|\, e^{-ct_2} < |\xi|. \tag{49}$$

Then there exists a $t_1 < t_2$ such that $|z(t_1)| = |\xi|$ and $|z(t)| \leq |\xi|$ for all $t \in [t_1, t_2]$. Let

$$w(t) := |z(t)|^2 / 2.$$

Then, for almost all $t \in [t_1, t_2]$

$$|\dot{w}(t)| = |z(t) \cdot \dot{z}(t)| \leq c\,|z(t)|^2 = 2cw(t).$$

In particular, this implies that $\dot{w}(t) + 2cw(t) \geq 0$. So, for all $t \in [t_1, t_2]$ we have

$$0 \leq e^{2ct}(\dot{w}(t) + 2cw(t)) = \frac{d(e^{2ct}w(t))}{dt},$$

implying that $e^{2ct}w(t) \geq e^{2ct_1}w(t_1)$ for all $t \in [t_1, t_2]$. Thus,

$$\frac{1}{2}|z(t_2)|^2 = w(t_2) \geq e^{2ct_1}w(t_1)e^{-2ct_2} = \frac{1}{2}|z(t_1)|^2\, e^{-2c(t_2 - t_1)} \geq \frac{1}{2}|\xi|^2\, e^{-2ct_2},$$

so that $|z(t_2)| \geq e^{-ct_2}|\xi|$, contradicting (49). ∎

**Corollary 4.3** For every $r > 0$ and $T > 0$ there is a $\sigma = \sigma(r, T) > 0$, such that for any $\mathbf{d} \in \mathcal{S}$, $\mathbf{v} \in \mathcal{W}$, $|\xi| \geq r$ and $t \leq T$ we have
$$|z(t, \xi, \mathbf{d}, \mathbf{v})| \geq \sigma.$$

□

**Lemma 4.4** For each $\mathbf{d} \in \mathcal{S}_T$ and each $\mathbf{v} \in \mathcal{W}_T$, if $\xi_k \to \xi$ in $\mathbb{X} \setminus \{0\}$, $\mathbf{d}_k \to \mathbf{d}$ in $\mathcal{S}_T$, and $\mathbf{v}_k \to \mathbf{v}$ in $\mathcal{W}_T$, $\{z(t, \xi_k, \mathbf{d}_k, \mathbf{v}_k)\}$ converges to $z(t, \xi, \mathbf{d}, \mathbf{v})$ uniformly on $[0, T]$.

*Proof.* Assume without loss of generality that $\xi_k \in U := B_{\frac{|\xi|}{2}}(\xi)$. Since $\left|\widetilde{f}\right| \leq 3$ and by Corollary 4.3, $\mathcal{R}_T(U) \subseteq B_{1.5|\xi|+3T}(0) \setminus B_\sigma(0)$, where $\sigma = \sigma(|\xi|/2, T)$ is as in Corollary 4.3. Let $M_1$ and $M_2$ be Lipschitz constants for $\hat{f}(\cdot, d)$ (uniformly for $d \in \Omega$) and $2f_0 \phi$ respectively on $B_{1.5|\xi|+3T}(0) \setminus B_\sigma(0)$. Write $z(t) := z(t, \xi, \mathbf{d}, \mathbf{v})$, and $z_k(t) := z(t, \xi_k, \mathbf{d}_k, \mathbf{v}_k)$. Now, for all $t \in [0, T]$ we have

$$
\begin{aligned}
|z(t) - z_k(t)| &= \left| \xi_k - \xi + \int_0^t \left( \widetilde{f}(z_k(t), \mathbf{d}_k(t), \mathbf{v}_k(t)) - \widetilde{f}(z(t), \mathbf{d}(t), \mathbf{v}(t)) \right) dt \right| \\
&\leq |\xi_k - \xi| + \left| \int_0^t \left( \widetilde{f}(z(t), \mathbf{d}_k(t), \mathbf{v}_k(t)) - \widetilde{f}(z(t), \mathbf{d}(t), \mathbf{v}(t)) \right) dt \right| \\
&\quad + \int_0^t \left| \widetilde{f}(z_k(t), \mathbf{d}_k(t), \mathbf{v}_k(t)) - \widetilde{f}(z(t), \mathbf{d}_k(t), \mathbf{v}_k(t)) \right| dt \\
&\leq |\xi_k - \xi| + \left| \int_0^t \left( \hat{f}(z(t), \mathbf{d}_k(t)) - \hat{f}(z(t), \mathbf{d}(t)) \right) dt \right| \quad (50) \\
&\quad + 2 \left| \int_0^t (f_0(z(t))\phi(z(t))\mathbf{v}_k(t) - f_0(z(t))\phi(z(t))\mathbf{v}(t)) \, dt \right| \quad (51) \\
&\quad + \int_0^t (M_1 + M_2) |z_k(t) - z(t)| \, dt.
\end{aligned}
$$

The integrals in (50) and (51) tend to 0 because of the convergence of $\{\mathbf{d}_k(\cdot)\}$ to $\mathbf{d}$ in $\mathcal{S}$ and the weak convergence of $\mathbf{v}_k$ to $\mathbf{v}$. So, for any $\varepsilon > 0$ we can find a $K$ such that, for all $k \geq K$,

$$
\begin{aligned}
|\xi_k - \xi| &+ \left| \int_0^T \left( \hat{f}(z(t), \mathbf{d}_k(t)) - \hat{f}(z(t), \mathbf{d}(t)) \right) dt \right| \\
&+ 2 \left| \int_0^t (f_0(z(t))\phi(z(t))\mathbf{v}_k(t) - f_0(z(t))\phi(z(t))\mathbf{v}(t)) \, dt \right| \leq \varepsilon e^{-(M_1+M_2)T}.
\end{aligned}
$$

Then, for all $k \geq K$ and $t \in [0, T]$ we have, by Gronwall inequality,

$$
|z(t) - z_k(t)| \leq \varepsilon e^{-(M_1+M_2)T} e^{(M_1+M_2)t} \leq \varepsilon
$$

As $\varepsilon$ was arbitrary, this proves uniform convergence. ∎

Also, since for any $\xi \in \partial D \setminus \{0\}$, $f_0(\xi) \neq 0$ and $f_0(\xi) \geq \left|\hat{f}(\xi, d)\right|$ for any $d \in \Omega$, we have the following controllability property on $\partial \mathcal{D} \setminus \{0\}$:

**Lemma 4.5** Let $\xi \neq 0$ be on $\partial \mathcal{D}$. Then, for each $\tau > 0$, there exists a neighborhood $U$ of $\xi$, such that for any $\eta \in U$ and any $\mathbf{d} \in \mathcal{M}_\Omega$, there is some control $\mathbf{v}$, and some $0 \leq t_1 \leq \tau$, such that $z(t_1, \eta, \mathbf{d}, \mathbf{v}) = \xi$ and $z(t, \eta, \mathbf{d}, \mathbf{v}) \in U$ for all $0 \leq t \leq t_1$.

*Proof.* Since $\phi(\xi) f_0(\xi) \neq 0$ and the function $\phi(\cdot) f_0(\cdot)$ is continuous, we can find a ball $U_1$ centered at $\xi$ and a constant $c_1$ such that $\phi(z) f_0(z) > c_1$ for every $z \in U_1$. Since $\phi(\xi) = 1$ and $\phi$ is continuous, one could also find a ball $U_2 \subseteq U_1$ centered at $\xi$, so that

$$
\left|\hat{f}(z, d)\right| \leq 1.5 \phi(z) f_0(z), \quad \forall z \in U_2, d \in \Omega.
$$

Fix $\tau > 0$ and let $B(\xi)$ be the ball of radius $\tau c_1/2$ centered at $\xi$. Define $U := B(\xi) \cap U_2$. Pick a point $\eta \in U$. Then $|\xi - \eta| < \tau c_1/2$, so,

$$\bar{v}_2 := \frac{2(\xi - \eta)}{\tau c_1}$$

has norm smaller than 1. Consider the "feedback law"

$$k(z, d) = \frac{1}{2}(1.5 k_1(z, d) + 0.5 \bar{v}_2),$$

where

$$k_1(z, d) := -\frac{\widehat{f}(z, d)}{1.5 \phi(z) f_0(z)}.$$

Notice that for all $z \in U$ and $d \in \Omega$ we have $|k_1(z, d)| \leq 1$ and

$$\begin{aligned}
\widetilde{f}(z, d, k(z, d)) &= \widehat{f}(z, d) + 2\phi(z) f_0(z) k(z, d) \\
&= \widehat{f}(z, d) + 1.5 \phi(z) f_0(z) k_1(z, d) + 0.5 \phi(z) f_0(z) \bar{v}_2 \\
&= 0.5 \phi(z) f_0(z) \bar{v}_2 = \frac{(\xi - \eta) \phi(z) f_0(z)}{\tau c_1}.
\end{aligned}$$

Thus, with any initial condition $\eta \in U$ and disturbance $\mathbf{d} \in \mathcal{M}_\Omega$, if the control $v(t) := k(z(t), d(t))$ is applied, then the trajectory of the system $\widetilde{\Sigma}$ will be the line segment, connecting $\xi$ and $\eta$, transversed with a velocity greater than $(\xi - \eta)/\tau$. So, there exists a $t_0 \leq \tau$ such that $z(t_0, \eta, \mathbf{d}, \mathbf{v}) = \xi$ and, since $U$ is convex, $z(t, \eta, \mathbf{d}, \mathbf{v}) \in U$ for all $t \geq t_0$. ∎

For each $\mathbf{d} \in \mathcal{S}$, let

$$\theta_{\mathbf{d}}(\xi, \mathbf{v}) = \inf \{t \geq 0 : z(t, \xi, \mathbf{d}, \mathbf{v}) \in \mathcal{D}\},$$

where as before, $\theta_{\mathbf{d}}(\xi, \mathbf{v}) = \infty$ if the trajectory never reaches $\mathcal{D}$.

**Lemma 4.6** The map $(\xi, \mathbf{v}, \mathbf{d}) \mapsto \theta_{\mathbf{d}}(\cdot, \cdot)$ is lower semicontinuous on $\mathcal{E} \times \mathcal{W} \times \mathcal{S}$.

*Proof.* Let $\{\xi_k\} \subset \mathcal{E}$, $\{\mathbf{v}_k\} \subset \mathcal{W}$ and $\{\mathbf{d}_k\} \subset \mathcal{S}$ be such that $\xi_k \to \xi$, $\mathbf{v}_k \to \mathbf{v}$, and $\mathbf{d}_k \to \mathbf{d}$ for some $\xi \in \mathcal{E}$, $\mathbf{v} \in \mathcal{W}$, and $\mathbf{d} \in \mathcal{S}$. We need to show that

$$\theta_{\mathbf{d}}(\xi, \mathbf{v}) \leq \liminf_{k \to \infty} \theta_{\mathbf{d}_k}(\xi_k, \mathbf{v}_k). \tag{52}$$

Let $\theta_k = \theta_{\mathbf{d}_k}(\xi_k, \mathbf{v}_k)$. Without loss of generality, we may assume that

$$\liminf_{k \to \infty} \theta_k = \theta_0 < \infty.$$

Passing to a subsequence if necessary, we assume that $\theta_k \to \theta_0$. Thus, there exists some $K$ such that $\theta_k \leq \theta_0 + 1$ for all $k \geq K$. Since $\{z(t, \xi_k, \mathbf{d}_k, \mathbf{v}_k)\}$ converges to $z(t, \xi, \mathbf{d}, \mathbf{v})$ uniformly on $[0, \theta_0 + 1]$, it follows that

$$z(\theta_0, \xi, \mathbf{d}, \mathbf{v}) = \lim_{k \to \infty} z(\theta_k, \xi_k, \mathbf{d}_k, \mathbf{v}_k).$$

Since $\mathcal{D}$ is closed and $z(\theta_k, \xi_k, \mathbf{d}_k, \mathbf{v}_k) \in \mathcal{D}$ for each $k$, we know that $z(\theta_0, \xi, \mathbf{d}, \mathbf{v}) \in \mathcal{D}$, and hence, $\theta_{\mathbf{d}}(\xi, \mathbf{v}) \leq \theta_0$. ∎

Define, for $\xi \in \mathcal{E}, \mathbf{d} \in \mathcal{S}$ and $\mathbf{v} \in \mathcal{W}$,

$$V_{\mathbf{v}}(\xi, \mathbf{d}) := \int_0^{\theta_{\mathbf{d}}(\xi,\mathbf{v})} \Xi(|z(t,\xi,\mathbf{d},\mathbf{v})|)\, dt,$$

and, for $\xi \in \mathcal{E}$ and $\mathbf{d} \in \mathcal{S}$,

$$\widetilde{V}_0(\xi, \mathbf{d}) := \inf_{\mathbf{v} \in \mathcal{W}} V_{\mathbf{v}}(\xi, \mathbf{d}).$$

Note that for some $\mathbf{v}$ and $\mathbf{d}$, $V_{\mathbf{v}}(\xi, \mathbf{d})$ may take $\infty$ as its value, but $\widetilde{V}_0(\xi, \mathbf{d})$ is always finite, since

$$\widetilde{V}_0(\xi, \mathbf{d}) \leq V_O(\xi, \mathbf{d}) \qquad (53)$$

where $O(\cdot)$ is the control identically equal to 0. Recall that, by definition of $\Xi$ we have then, for all $\xi \in \mathcal{E}_1$ and $\mathbf{d} \in \mathcal{S}$,

$$\begin{aligned} V_O(\xi, \mathbf{d}) &= \int_0^{\theta_{\mathbf{d}}(\xi,O)} \Xi(|z(t,\xi,\mathbf{d},O)|)\, dt \\ &= \int_0^{\theta_{\mathbf{d}}(\xi,O)} \Xi(|z(t,\xi,\mathbf{d})|)\, dt \leq \int_0^{\infty} \Xi(\mu_1(\mu_2(|\xi|)e^{-t}))\, dt \leq \mu_2(|\xi|), \end{aligned} \qquad (54)$$

where $\mu_1$ and $\mu_2$ are as in equation (45).

**Lemma 4.7** *The function* $V_{(\cdot)}(\cdot,\cdot) : \mathcal{E} \times \mathcal{S} \times \mathcal{W} \to \mathbb{R}_{\geq 0}$ *is lower semicontinuous.*

*Proof.* Let $(\xi, \mathbf{d}, \mathbf{v}) \in \mathcal{E} \times \mathcal{S} \times \mathcal{W}$, and let $\{\xi_k\} \to \xi$, $\{\mathbf{d}_k\} \to \mathbf{d}$ and $\{\mathbf{v}_k\} \to \mathbf{v}$, where $\xi_k \in \mathcal{E}$ for all $k$.

*Case 1.* $V_{\mathbf{v}}(\xi, \mathbf{d}) < \infty$. In this case, for any $\varepsilon > 0$, there exists some $0 < T < \theta_{\mathbf{d}}(\xi, \mathbf{v})$ such that

$$V_{\mathbf{v}}(\xi, \mathbf{d}) = \int_0^{\theta_{\mathbf{d}}(\xi,\mathbf{v})} \Xi(|z(t,\xi,\mathbf{d},\mathbf{v})|)\, dt \leq \int_0^{T} \Xi(|z(t,\xi,\mathbf{d},\mathbf{v})|)\, dt + \varepsilon.$$

Without loss of generality we can assume that all $\xi_k$ are within the unit distance from $\xi$. Recall that the reachable set from the unit ball around $\xi$ is bounded. Since $z(t, \xi_k, \mathbf{d}_k, \mathbf{v}_k)$ converges to $z(t, \xi, \mathbf{d}, \mathbf{v})$ uniformly on $[0, T]$, and $\Xi(\cdot)$ is uniformly continuous on compacts, there exists some $K > 0$ such that

$$|\Xi(|z(t,\xi,\mathbf{d},\mathbf{v})|) - \Xi(|z(t,\xi_k,\mathbf{d}_k,\mathbf{v}_k)|)| < \frac{\varepsilon}{1+T} \quad \forall k > K,\ \forall t \in [0,T].$$

This implies that

$$\int_0^T \Xi(|z(t,\xi_k,\mathbf{d}_k,\mathbf{v}_k)|)\, dt \geq \int_0^T \Xi(|z(t,\xi,\mathbf{d},\mathbf{v})|)\, dt - \varepsilon$$

for all $k \geq K$. By Lemma 4.6, there exists some $K_1 \geq K$ such that $\theta_{\mathbf{d}_k}(\xi_k, \mathbf{v}_k) > T$ for all $k \geq K_1$. Thus, for all $k \geq K_1$,

$$\begin{aligned} V_{\mathbf{v}_k}(\xi_k, \mathbf{d}_k) &\geq \int_0^T \Xi(|z(t,\xi_k,\mathbf{d}_k,\mathbf{v}_k)|)\, dt \\ &\geq \int_0^T \Xi(|z(t,\xi,\mathbf{d},\mathbf{v})|)\, dt - \varepsilon \geq V_{\mathbf{v}}(\xi, \mathbf{d}) - 2\varepsilon. \end{aligned}$$

As $\varepsilon$ was arbitrary, we conclude that

$$V_{\mathbf{v}}(\xi, \mathbf{d}) \leq \liminf V_{\mathbf{v}_k}(\xi_k, \mathbf{d}_k).$$

*Case 2.* $V_{\mathbf{v}}(\xi, \mathbf{d}) = \infty$.

In this case, $\theta_{\mathbf{d}}(\xi, \mathbf{v}) = \infty$. Fix an integer $k \geq 0$. There exists some $T_k$ such that

$$\int_0^{T_k} \Xi(|z(t, \xi, \mathbf{d}, \mathbf{v})|)\, dt \geq k.$$

Repeating the same argument used above, one sees that

$$\int_0^{T_k} \Xi(|z(t, \xi_l, \mathbf{d}_l, \mathbf{v}_l)|)\, dt \geq k - 1$$

for all $l \geq L_0$ for some $L_0$. By Lemma 4.6, there is some $L_1 \geq L_0$ such that $\theta_{\mathbf{d}_l}(\xi_l, \mathbf{v}_l) \geq T_k$ for all $l \geq L_1$. Consequently, for all $l \geq L_1$,

$$V_{\mathbf{v}_l}(\xi_l, \mathbf{d}_l) \geq \int_0^{T_k} \Xi(|z(t, \xi_l, \mathbf{d}_l, \mathbf{v}_l)|)\, dt \geq k - 1.$$

Since $k > 0$ can be picked arbitrarily, it follows that

$$\liminf V_{\mathbf{v}_l}(\xi_l, \mathbf{d}_l) = \infty.$$

In both cases, we have shown that $\liminf V_{\mathbf{v}_l}(\xi_l, \mathbf{d}_l) \geq V_{\mathbf{v}}(\xi, \mathbf{d})$. The lower semicontinuity property follows readily. ∎

**Lemma 4.8** For every $\xi \in \mathcal{E}$, $\mathbf{d} \in \mathcal{S}$, there exists a control $\bar{\mathbf{v}}$ such that

$$V_{\bar{\mathbf{v}}}(\xi, \mathbf{d}) = \widetilde{V}_0(\xi, \mathbf{d}). \tag{55}$$

*Proof.* Let $\mathbf{v}_k$ be a sequence of controls such that $V_{\mathbf{v}_k}(\xi, \mathbf{d}) \searrow \widetilde{V}_0(\xi, \mathbf{d})$. Without loss of generality we are assuming that all these controls are defined for all positive $t$ (by letting them equal to 0 where they are not defined). Extract from $\{\mathbf{v}_k\}$ a subsequence $\{\mathbf{v}_{k_l}\}$ converging weakly to some limit $\bar{v}$ in $\mathcal{W}$. Without relabeling, we assume that $\mathbf{v}_k \to \bar{\mathbf{v}}$. By Lemma 4.7,

$$V_{\bar{\mathbf{v}}}(\xi, \mathbf{d}) \leq \lim_{k \to \infty} V_{\mathbf{v}_k}(\xi, \mathbf{d}) = \widetilde{V}_0(\xi, \mathbf{d}).$$

Combining this with the fact that $\widetilde{V}_0(\xi, \mathbf{d}) \leq V_{\mathbf{v}}(\xi, \mathbf{d})$ for all $\mathbf{v} \in \mathcal{W}$, one proves (55). ∎

**Corollary 4.9** For any $\xi \in \mathcal{E}, \mathbf{d} \in \mathcal{S}, \mathbf{v} \in \mathcal{W}$, and $0 \leq T < \theta_{\mathbf{d}}(\xi, \mathbf{v})$, it holds that

$$\widetilde{V}_0(\xi, \mathbf{d}) \leq \int_0^T \Xi(|z(s, \xi, \mathbf{d}, \mathbf{v})|)\, ds + \widetilde{V}_0(z(T, \xi, \mathbf{d}, \mathbf{v}), \mathbf{d}_T), \tag{56}$$

where $\mathbf{d}_T(t) = \mathbf{d}(t + T)$ for all $t \geq 0$.

*Proof.* Suppose the assertion is not true. Then there exist $\xi \in \mathcal{E}$, $T > 0$, $\mathbf{v} \in \mathcal{W}$, and $\mathbf{d} \in \mathcal{S}$ such that (56) fails. By Lemma 4.8, one can find a control $\mathbf{v}_1$ such that $\widetilde{V}_0(z(T, \xi, \mathbf{d}, \mathbf{v}), \mathbf{d}_T) = V_{\mathbf{v}_1}(z(T, \xi, \mathbf{d}, \mathbf{v}), \mathbf{d}_T)$. Define $\bar{\mathbf{v}}$ to be the concatenation of $\mathbf{v}$ and $\mathbf{v}_1$. Then, letting $\theta := \theta_{\mathbf{d}_T}(z(T, \xi, \mathbf{d}, \mathbf{v}), \mathbf{v}_1)$ and noticing that $\theta_{\mathbf{d}}(\xi, \bar{\mathbf{v}}) = \theta + T$, we get, by our assumption,

$$\begin{aligned}
\widetilde{V}_0(\xi, \mathbf{d}) &> \int_0^T \Xi(|z(s, \xi, \mathbf{d}, \mathbf{v})|) \, ds + \widetilde{V}_0(z(T, \xi, \mathbf{d}, \mathbf{v}), \mathbf{d}_T) \\
&= \int_0^T \Xi(|z(s, \xi, \mathbf{d}, \mathbf{v})|) \, ds + \int_0^\theta \Xi(|z(t, z(T, \xi, \mathbf{d}, \mathbf{v}), \mathbf{d}_T, \mathbf{v}_1)|) \, dt \\
&= \int_0^{\theta+T} \Xi(|z(s, \xi, \mathbf{d}, \bar{\mathbf{v}})|) \, ds,
\end{aligned}$$

which contradicts with the minimality of $\widetilde{V}_0(\xi, \mathbf{d})$. ∎

**Lemma 4.10** For each $\xi \in \partial \mathcal{D}$, the following holds:

$$\lim_{\eta \to \xi} \widetilde{V}_0(\eta, \mathbf{d}) = 0, \tag{57}$$

uniformly in $\mathbf{d} \in \mathcal{S}$, that is, for any $\varepsilon > 0$, there is a neighborhood $U$ of $\xi$, such that $\widetilde{V}_0(\eta, \mathbf{d}) < \varepsilon$ for all $\eta \in U \cap \mathcal{E}$, all $\mathbf{d} \in \mathcal{S}$.

*Proof.* If $\xi = 0$, the result follows from (53) and (54) . Suppose now that $\xi \neq 0$.

Let $\varepsilon > 0$ be given. Let $\tau = \frac{\varepsilon}{\Xi(|\xi|+1)}$. Find a neighborhood $U$ of $\xi$ as in Lemma 4.5. Shrinking $U$ if necessary, we assume that $|\xi - \eta| \leq 1$ for all $\eta \in U$.

Suppose $\eta \in U \cap \mathcal{E}$ and $\mathbf{d} \in \mathcal{M}_\Omega$. By the controllability property, there is some control $v$ such that $z(t_1, \eta, \mathbf{d}, \mathbf{v}) = \xi \in \partial \mathcal{D}$ for some $t_1 \in [0, \tau]$, and that $z(t, \eta, \mathbf{d}, \mathbf{v}) \in U$ for all $t \in [0, t_1]$. Thus,

$$V_{\mathbf{v}}(\eta, \mathbf{d}) \leq \int_0^{t_1} \Xi(|z(s, \eta, \mathbf{d}, \mathbf{v})|) \, ds \leq \tau \cdot \Xi(|\xi| + 1) \leq \varepsilon,$$

from which it follows that $\widetilde{V}_0(\eta, \mathbf{d}) \leq \varepsilon$.

The above shows that $\widetilde{V}_0(\eta, \mathbf{d}) \leq \varepsilon$ for all non-relaxed $\mathbf{d} \in \mathcal{M}_\Omega$, $\eta \in U \cap \mathcal{E}$. Pick a relaxed $\mathbf{d} \in \mathcal{S}$. Then there exists a sequence $\{\mathbf{d}_k\} \subset \mathcal{M}_\Omega$ such that $\mathbf{d}_k \to \mathbf{d}$ in the topology of relaxed controls. Let $\eta \in U \cap \mathcal{E}$. Let $\mathbf{v}_k$ be such that $\widetilde{V}_0(\eta, \mathbf{d}_k) = V_{\mathbf{v}_k}(\eta, \mathbf{d}_k)$. By weak sequential compactness of $\mathcal{W}$, one may assume, after taking a subsequence, that $\mathbf{v}_k \to \bar{\mathbf{v}}$ for some $\bar{\mathbf{v}}$. By Lemma 4.7,

$$V_{\bar{\mathbf{v}}}(\eta, \mathbf{d}) \leq \liminf_{k \to \infty} V_{\mathbf{v}_k}(\eta, \mathbf{d}_k) \leq \varepsilon,$$

and consequently, $\widetilde{V}_0(\eta, \mathbf{d}) \leq \varepsilon$. This shows that $\widetilde{V}_0(\eta, \mathbf{d}) \leq \varepsilon$ for all $\eta \in U \cap \mathcal{E}$, all $d \in \mathcal{S}$. ∎

To prove the continuity of $\widetilde{V}_0$, we also need the following result.

**Lemma 4.11** Suppose for some $\xi \in \mathcal{E}$, $\mathbf{d} \in \mathcal{S}$, and $\mathbf{v} \in \mathcal{W}$, $V_{\mathbf{v}}(\xi, \mathbf{d}) < \infty$. Then there exists some $\xi_0 \in \partial \mathcal{D}$ such that

$$\lim_{t \to \theta_{\mathbf{d}}(\xi, \mathbf{v})} z(t, \xi, \mathbf{d}, \mathbf{v}) = \xi_0. \tag{58}$$

*Proof.* Suppose $V_{\mathbf{v}}(\xi,\nu) < \infty$. This means that

$$\int_0^{\theta_\nu(\xi,\mathbf{v})} \Xi(|z(s,\xi,\nu,\mathbf{v})|)\,ds < \infty. \tag{59}$$

If $\theta_{\mathbf{d}}(\xi,\mathbf{v}) < \infty$, then (58) follows from the continuity of $z(\cdot,\xi,\nu,\mathbf{v})$ with $\xi_0 = z(\theta_{\mathbf{d}}(\xi,\mathbf{v}))$.

Suppose now that $\theta_{\mathbf{d}}(\xi,\mathbf{v}) = \infty$. Since the integral in (59) converges, $\int_t^\infty \Xi(|z(s,\xi,\mathbf{d},\mathbf{v})|)\,ds$ decreases to 0 as $t$ tends to $\infty$. Consider the family of functions $\{x_t(\cdot),\ t>0\}$, defined by $x_t(s) := z(t+s,\xi,\mathbf{d},\mathbf{v})$, $\mathcal{I}_t = [0,\infty)$. By Lemma 4.3, the trajectory $z(s,\xi,\mathbf{d},\mathbf{v})$ does not reach the origin in finite time, hence, there exists a positive, strictly decreasing function $\varphi$ such that $\varphi(s) < |z(s,\xi,\mathbf{d},\mathbf{v})|$ for all $s > 0$. Find a $\mathcal{K}_\infty$-function $\kappa$ such that

$$\kappa(\varphi(t)) > \int_t^\infty \Xi(|z(s,\xi,\mathbf{d},\mathbf{v})|)\,ds = \int_0^\infty \Xi(|x_t(s)|)\,ds.$$

Then the family $\{x_t(\cdot),\ t>0\}$ satisfies all the conditions of Proposition 3.9 (with $\chi := \Xi$). Take $r := |\xi|$. Then, for any $\varepsilon > 0$, $|z(t,\xi,\mathbf{d},\mathbf{v})| < \varepsilon$ for all $t > T_{r,\varepsilon}$. So, the conclusion of the lemma follows. ∎

**Proposition 4.12** *The function $\widetilde{V}_0 : \mathcal{E} \times \mathcal{S} \to \mathbb{R}$ is continuous.*

*Proof.* Fix $\xi \in \mathcal{E}$, $\mathbf{d} \in \mathcal{S}$. Suppose $\xi_k \to \xi$, $\mathbf{d}_k \to \mathbf{d}$, where $\xi_k \in \mathcal{E}$. Let $\{k_j\}$ be a subsequence of $\{k\}$ such that

$$\lim_{j\to\infty} \widetilde{V}_0(\xi_{k_j}, \mathbf{d}_{k_j}) = \liminf_{k\to\infty} \widetilde{V}_0(\xi_k, \mathbf{d}_k). \tag{60}$$

For each $k$, let $\mathbf{v}_k$ be such that $\widetilde{V}_0(\xi_k, \mathbf{d}_k) = V_{\mathbf{v}_k}(\xi_k, \mathbf{d}_k)$. Notice that $\lim_{k\to\infty} V_{\mathbf{v}_k}(\xi_k, \mathbf{d}_k)$ exists, because of (60).

By sequential compactness of $\mathcal{W}$, there exists a subsequence of $\{\mathbf{v}_{k_j}\}$ converging to some $\bar{\mathbf{v}} \in \mathcal{W}$. Without relabeling, we assume that $\mathbf{v}_{k_j} \to \bar{\mathbf{v}}$. It then follows from Lemma 4.7 that

$$V_{\bar{\mathbf{v}}}(\xi,\mathbf{d}) \leq \lim V_{\mathbf{v}_{k_j}}(\xi_{k_j},\mathbf{d}_{k_j}) = \liminf \widetilde{V}_0(\xi_k,\mathbf{d}_k).$$

Consequently, $\widetilde{V}_0(\xi,\mathbf{d}) \leq \liminf \widetilde{V}_0(\xi_k,\mathbf{d}_k)$. To complete the proof, we will show that

$$\widetilde{V}_0(\xi,\mathbf{d}) \geq \limsup_{k\to\infty} \widetilde{V}_0(\xi_k,\mathbf{d}_k). \tag{61}$$

Let $\mathbf{v}$ be a control such that $\widetilde{V}_0(\xi,\mathbf{d}) = V_{\mathbf{v}}(\xi,\mathbf{d})$. Let $\varepsilon > 0$ be given. By Lemma 4.11, there is some $\xi_0 \in \partial \mathcal{D}$ such that (58) holds. By Lemma 4.10, there is a neighborhood $U$ of $\xi_0$ such that

$$\widetilde{V}_0(\eta,\nu) < \varepsilon/4 \ \text{ for all } \eta \in U \cap \mathcal{E},\ \text{all } \nu \in \mathcal{S}. \tag{62}$$

Let $0 < T < \theta_{\mathbf{d}}(\xi,\mathbf{v})$ be such that $z(T,\xi,\mathbf{d},\mathbf{v}) \in U$. Then, since $\{z(t,\xi_k,\mathbf{d}_k,\mathbf{v})\}$ converges to $z(t,\xi,\mathbf{d},\mathbf{v})$ uniformly on $[0,T]$, it follows that $z(T,\xi_k,\mathbf{d}_k,\mathbf{v}) \in U$ for $k \geq K_1$ for some $K_1$. By Lemma 4.6, one may assume that $T < \theta_{\mathbf{d}_k}(\xi_k,\mathbf{v})$ for all $k \geq K_1$. Consequently, $\eta = z(T,\xi_k,\mathbf{d}_k,\mathbf{v})$ is also in $\mathcal{E}$, so, applying (62) with $\nu = (\mathbf{d}_k)_T$, we have

$$\widetilde{V}_0(z(T,\xi_k,\mathbf{d}_k,\mathbf{v}),(\mathbf{d}_k)_T) < \varepsilon/4, \qquad \forall\, k \geq K_1,$$

where $(\mathbf{d}_k)_T(t) = \mathbf{d}_k(T+t)$. Using the uniform convergence property of $\{z(t, \xi_k, \mathbf{d}_k, \mathbf{v})\}$, it follows that there is some compact set $\mathcal{K}$ such that $z(t, \xi_k, \mathbf{d}_k, \mathbf{v}) \in \mathcal{K}$ for all $k$, all $t \in [0, T]$. Using also the uniform continuity of $\Xi(\cdot)$ on compacts, one sees that there is some $K_2 \geq K_1$ such that

$$\int_0^T \Xi(|z(s, \xi_k, \mathbf{d}_k, \mathbf{v})|)\, ds \leq \int_0^T \Xi(|z(s, \xi, \mathbf{d}, \mathbf{v})|)\, ds + \varepsilon/2$$

for all $k \geq K_2$. Thus, Equation (56) implies

$$\begin{aligned}\widetilde{V}_0(\xi_k, \mathbf{d}_k) &\leq \int_0^T \Xi(|z(s, \xi_k, \mathbf{d}_k, \mathbf{v})|)\, ds + \widetilde{V}_0(z(T, \xi_k, \mathbf{d}_k, \mathbf{v}), (\mathbf{d}_k)_T) \\ &\leq V_{\mathbf{v}}(\xi, \mathbf{d}) + \varepsilon/2 + \varepsilon/2 = \widetilde{V}_0(\xi, \mathbf{d}) + \varepsilon\end{aligned}$$

for all $k \geq K_2$. From this it follows that

$$\widetilde{V}_0(\xi, \mathbf{d}) \geq \limsup \widetilde{V}_0(\xi_k, \mathbf{d}_k) - \varepsilon.$$

Letting $\varepsilon \to 0$, one proves (61). ∎

For each $\xi \in \mathcal{E}$, define

$$V_0(\xi) := \sup_{\mathbf{d} \in \mathcal{M}_\Omega} \widetilde{V}_0(\xi, \mathbf{d}). \tag{63}$$

Note that the supremum is finite for each $\xi \in \mathcal{E}$; in fact, by (53) and (54), we have the upper bound

$$V_0(\xi) \leq \mu_2(|\xi|), \qquad \forall\, \xi \in \mathcal{E}, \tag{64}$$

where $\mu_2$ is as in inequality (45). Also observe that the same function $V_0$ results if the supremum in (63) is taken over the set $\mathcal{S}$. This follows from the fact that $\mathcal{M}_\Omega$ is dense in $\mathcal{S}$ and the continuity property of $\widetilde{V}_0$.

By sequential compactness of $\mathcal{S}$ and by continuity of $\widetilde{V}_0$, we get the following:

**Corollary 4.13** For any $\xi$ in $\mathcal{E}$ there exists a (possibly relaxed) disturbance $\mathbf{d}$ such that $V_0(\xi) = \widetilde{V}_0(\xi, \mathbf{d})$, that is, $V_0(\xi) = \max_{\bar{\mathbf{d}} \in \mathcal{S}} \widetilde{V}_0(\xi, \bar{\mathbf{d}})$. Moreover, $V_0 : \mathcal{E} \to \mathbb{R}$ is continuous.

*Proof.* Fix $\xi \in \mathcal{E}$. By definition of $V_0$, there exists a sequence of disturbances $\{\mathbf{d}_k(\cdot)\}$ such that $\widetilde{V}_0(\xi, \mathbf{d}_k) \nearrow V_0(\xi)$. By sequential compactness of $\mathcal{S}$, we can extract a subsequence $\{\mathbf{d}_{k_i}\}$, converging in $\mathcal{S}$ to some $\mathbf{d}$. By continuity of $\widetilde{V}_0$,

$$V_0(\xi) = \lim_{i \to \infty} \widetilde{V}_0(\xi, \mathbf{d}_{k_i}) = \widetilde{V}_0(\xi, \mathbf{d}),$$

proving the first statement of the corollary.

Now fix $\xi \in \mathcal{E}$. Take any sequence $\{\xi_k\} \in \mathcal{E}$, converging to $\xi$ and such that the sequence $\{V_0(\xi_k)\}$ converges. Let $\mathbf{d}_k(\cdot)$ and $\mathbf{d} \in \mathcal{S}$ be maximizing disturbances for $\xi_k$ and $\xi$ respectively, that is,

$$V_0(\xi_k) = \widetilde{V}_0(\xi_k, \mathbf{d}_k) \text{ and } V_0(\xi) = \widetilde{V}_0(\xi, \mathbf{d}).$$

Extracting a subsequence if necessary, let $\hat{\mathbf{d}}$ be a limit of $\{\mathbf{d}_k\}$ in $\mathcal{S}$. Then, by continuity of $\widetilde{V}_0$ and by definition of $V_0$, we have

$$V_0(\xi) \geq \widetilde{V}_0(\xi, \hat{\mathbf{d}}) = \lim_{k \to \infty} \widetilde{V}_0(\xi_k, \mathbf{d}_k) = \lim_{k \to \infty} V_0(\xi_k).$$

Consequently, if $\{\xi_i\} \subset \mathcal{E}$ is any sequence, converging to $\xi$, then

$$V_0(\xi) \geq \limsup_{i \to \infty} V_0(\xi_i),$$

proving upper-semicontinuity of $V_0$.

On the other hand, again by continuity of $\widetilde{V_0}$, we have

$$\liminf_{k \to \infty} V_0(\xi_k) \geq \lim_{k \to \infty} \widetilde{V_0}(\xi_k, \mathbf{d}) = V_0(\xi),$$

showing the lower-semicontinuity of $V_0$. Thus we conclude that $V_0$ is continuous on $\mathcal{E}$. ∎

**Lemma 4.14** There exists a $\mathcal{K}_\infty$-function $\underline{\alpha}$ such that

$$\underline{\alpha}(|\xi|) \leq V_0(\xi)$$

for all $\xi \in \mathcal{E}_1$. □

Notice that if $|\xi| \geq 1.6\rho(|h(\xi)|)$ and $\xi \neq 0$, then $\xi \in \mathcal{E}$, because $|\xi| > |\xi|/1.6 \geq \rho(|h(\xi)|)$.
*Proof.* (of Lemma 4.14) Recall that if $\xi \notin \mathcal{D} \cup \mathcal{B}$, then $\widetilde{f}(\xi, d, v) = \hat{f}(\xi, d)$ for any $d$ and $v$, so that

$$\left|\nabla(1.5\rho \circ |h|)(\xi) \cdot \widetilde{f}(\xi, d, v)\right| = \frac{1.5\left|\nabla(\rho \circ |h|)(\xi) \cdot f(\xi, d)\right|}{1 + |f(\xi, d)|^2 + \kappa(\xi)} \leq 1,$$

where the last inequality follows from (40). Therefore the assumptions of Lemma 4.1 are satisfied with $p := 1.5\rho$ and $f(\xi, d) := \widetilde{f}(\xi, d, v) = \hat{f}(\xi, d)$. By Lemma 4.1, we can find a constant $K_0$ such that if $|\xi| > K_0$ and $|\xi| \geq 1.6\rho(|h(\xi)|)$, then $z(t, \xi, \mathbf{d}, v) \notin \mathcal{D} \cup \mathcal{B}$ for all $t \in [0, 1)$. In particular, for such a $\xi$ we will have $\theta_\mathbf{d}(\xi, \mathbf{v}) > 1$ and $|z(t, \xi, \mathbf{d}, \mathbf{v})| > |\xi| - 1$ for all positive $t < 1$. Hence, the inequality

$$\Xi(|z(t, \xi, \mathbf{d}, \mathbf{v})|) > \Xi(|\xi| - 1), \quad \forall t \in [0, 1)$$

holds for any $\mathbf{d} \in \mathcal{M}_\Omega$, $\mathbf{v} \in \mathcal{W}$, and any $\xi$ such that

$$|\xi| > \max\{K_0 + 1,\ 1.6\rho(|h(\xi)|)\}. \tag{65}$$

Therefore, for any $\xi$ as in (65) and any $\mathbf{d} \in \mathcal{M}_\Omega$, the following estimate holds:

$$V_0(\xi) \geq \widetilde{V_0}(\xi, \mathbf{d}) \geq \int_0^1 \Xi(|z(t, \xi, \mathbf{d}, \mathbf{v})|)dt \geq \Xi(|\xi| - 1).$$

Next, notice that $V_0$ is strictly positive on $\mathcal{E}$. Indeed, $\left|\widetilde{f}(\xi, d, v)\right| \leq 3$ for any $\xi \in \mathcal{E}$, $d \in \Omega$, and $v \in [-1, 1]^n$, so that

$$|z(s, \xi, \mathbf{d}, \mathbf{v}) - \xi| \leq \frac{1}{2}\mathrm{dist}(\xi, \mathcal{D}) \quad \forall\, s \leq \frac{\frac{1}{2}\mathrm{dist}(\xi, \mathcal{D})}{3},\ \mathbf{d} \in \mathcal{M}_\Omega,\ \mathbf{v} \in \mathcal{W}.$$

Therefore $\theta_\mathbf{d}(\xi, \mathbf{v}) \geq \frac{1}{6}\mathrm{dist}(\xi, \mathcal{D})$ and $|z(s, \xi, \mathbf{d}, \mathbf{v})| \geq \xi - \frac{1}{2}\mathrm{dist}(\xi, \mathcal{D})$ for all $s \leq \frac{1}{6}\mathrm{dist}(\xi, \mathcal{D})$, all $\mathbf{d} \in \mathcal{M}_\Omega$, and $\mathbf{v} \in \mathcal{W}$. So, we have

$$V_\mathbf{v}(\xi, \mathbf{d}) \geq \int_0^{\mathrm{dist}(\xi, \mathcal{D})/6} \Xi\left(|z(s, \xi, \mathbf{d}, \mathbf{v})|\right)\, ds \geq \frac{\mathrm{dist}(\xi, \mathcal{D})}{6}\Xi(|\xi| - \mathrm{dist}(\xi, \mathcal{D})/2).$$

This shows that
$$\inf_{\mathbf{d}\in\mathcal{M}_\Omega}\inf_{v\in\mathcal{W}} V_{\mathbf{v}}(\xi,\mathbf{d}) \geq \frac{\mathrm{dist}(\xi,\mathcal{D})}{6}\Xi\left(|\xi| - \frac{\mathrm{dist}(\xi,\mathcal{D})}{2}\right). \tag{66}$$

Thus also the $\sup_{\mathbf{d}\in\mathcal{M}_\Omega}\inf_{v\in\mathcal{W}} V_{\mathbf{v}}(\xi,\mathbf{d})$ satisfies (66), and hence the same inequality holds for $V_0$.

Since $V_0$ is lower semicontinuous, it attains its minimum on any compact set. For each positive $l$ define
$$r_l := \frac{1}{l}\max\{K_0+1,\ 1.6\rho(|h(\xi)|)\} \quad\text{and}\quad m_l = \inf\{V_0(z) : z \in \mathcal{E}_1,\ r_l \leq |z| \leq r_1\}.$$

Since the sequence $\{m_l\}$ is non-increasing and positive, and $\Xi$ is of class $\mathcal{K}_\infty$, we can find a $\mathcal{K}_\infty$-function $\underline{\alpha}$ such that
$$\underline{\alpha}(s) < m_l \quad \forall\, s \in [r_l, r_{l-1}], \forall\, l > 1$$
and
$$\underline{\alpha}(s) < \Xi(s-1) \forall\, s \geq \max\{K_0+1,\ 1.6\rho(|h(\xi)|)\}.$$

By construction, $\underline{\alpha}$ will be a lower bound for $V_0$ on $\mathcal{E}_1$. ∎

Combining Lemma 4.14 with (64), we get the following, using $\bar{\alpha} := \mu_2$:
$$\underline{\alpha}(|\xi|) \leq V_0(\xi) \leq \bar{\alpha}(|\xi|), \qquad \forall\, \xi \in \mathcal{E}_1. \tag{67}$$

The following lemma and corollary summarize the dissipation properties for $V_0$.

**Lemma 4.15** For any $\xi \in \mathbb{X} \setminus (\mathcal{D} \cup \mathcal{B})$ and $\mathbf{d} \in \mathcal{M}_\Omega$, and any $t_0$ such that $z(t,\xi,\mathbf{d}) \notin \mathcal{D} \cup \mathcal{B}$ for all $t \in [0,t_0]$, the following dissipation inequality holds:
$$V_0(z(t,\xi,\mathbf{d})) - V_0(\xi) \leq -\int_0^{t_0} \Xi(|z(t,\xi,\mathbf{d})|)\, dt.$$

*Proof.* Fix $\xi \in \mathbb{X} \setminus (\mathcal{D} \cup \mathcal{B})$, $\mathbf{d} \in \mathcal{M}_\Omega$, and any positive $t_0$ as in the formulation of the lemma. Let $\varepsilon > 0$ be given. Find $\mathbf{d}_1$ such that
$$V_0(z(t_0,\xi,\mathbf{d})) - \widetilde{V}_0(z(t_0,\xi,\mathbf{d}),\mathbf{d}_1) < \varepsilon,$$
and let $\mathbf{v}_1 \in \mathcal{W}$ be a control such that $\widetilde{V}_0(z(t_0,\xi,\mathbf{d}),\mathbf{d}_1) = V_{\mathbf{v}_1}(z(t_0,\xi,\mathbf{d}),\mathbf{d}_1)$. Let $\widetilde{\mathbf{d}}$ be defined by
$$\widetilde{\mathbf{d}}(t) = \begin{cases} \mathbf{d}(t) & \text{if } 0 \leq t \leq t_0, \\ \mathbf{d}_1(t-t_0) & \text{if } t > t_0. \end{cases}$$

Then $z(t, z(t_0,\xi,\mathbf{d}), \mathbf{d}_1, \mathbf{v}_1) = z(t+t_0, \xi, \widetilde{\mathbf{d}}, \widetilde{\mathbf{v}})$ for all $t \geq 0$. By assumption,
$$z(t,\xi,\mathbf{d}) \notin \mathcal{D} \cup \mathcal{B} \qquad \forall\, t \in [0,t_0],$$
therefore we have
$$z(t,\xi,\mathbf{d}) = z(t,\xi,\mathbf{d},\mathbf{v}) \qquad \forall\, t \in [0,t_0], \quad \forall\, \mathbf{v} \in \mathcal{W}. \tag{68}$$

Notice also that (68) implies that $\theta_{\widetilde{\mathbf{d}}}(\xi, \mathbf{v}) > t_0$ for all $\mathbf{v} \in \mathcal{W}$ and

$$\begin{aligned}
\widetilde{V}_0(\xi, \widetilde{\mathbf{d}}) &= \min_{\mathbf{v} \in \mathcal{W}} \int_{t_0}^{\theta_{\widetilde{\mathbf{d}}}(\xi, \mathbf{v})} \Xi\left(\left|z(s, \xi, \widetilde{\mathbf{d}}, \mathbf{v})\right|\right) ds \\
&= \int_0^{t_0} \Xi\left(\left|z(s, \xi, \widetilde{\mathbf{d}})\right|\right) ds + \min_{\mathbf{v} \in \mathcal{W}} \int_{t_0}^{\theta_{\widetilde{\mathbf{d}}}(\xi, \mathbf{v})} \Xi\left(\left|z(s, \xi, \widetilde{\mathbf{d}}, \mathbf{v})\right|\right) ds \\
&= \int_0^{t_0} \Xi(|z(s, \xi, \mathbf{d})|) ds + \min_{\mathbf{v} \in \mathcal{W}} \int_0^{\theta_{\mathbf{d}_1}(z(t_0, \xi, \mathbf{d}), \mathbf{v})} \Xi(|z(s, z(t_0, \xi, \mathbf{d}), \mathbf{d}_1, \mathbf{v})|) ds \\
&= \int_0^{t_0} \Xi(|z(s, \xi, \mathbf{d})|) ds + \widetilde{V}_0(z(t_0, \xi, \mathbf{d}), \mathbf{d}_1).
\end{aligned}$$

Consequently, one has

$$\widetilde{V}_0(z(t_0, \xi, \mathbf{d}), \mathbf{d}_1) = \widetilde{V}_0(\xi, \widetilde{\mathbf{d}}) - \int_0^{t_0} \Xi\left(|z(s, \xi, \mathbf{d})|\right) ds.$$

Thus,

$$\begin{aligned}
V_0(z(t_0, \xi, \mathbf{d})) &\leq \widetilde{V}_0(z(t_0, \xi, \mathbf{d}), \mathbf{d}_1) + \varepsilon \\
&= \widetilde{V}_0(\xi, \widetilde{\mathbf{d}}) - \int_0^{t_0} \Xi(|z(s, \xi, \mathbf{d})|) ds + \varepsilon \\
&\leq V_0(\xi) - \int_0^{t_0} \Xi(|z(s, \xi, \mathbf{d})|) ds + \varepsilon.
\end{aligned}$$

Letting $\varepsilon \to 0$, we get the desired inequality. ∎

Thus, we have proven that the UOSS dissipation inequality holds for $V_0$ along the trajectories of the slower system $\hat{\Sigma}$ which are entirely contained in $\mathbb{X} \setminus (\mathcal{D} \cup \mathcal{B})$. It follows immediately that the same estimate holds along the trajectories of the original system $\Sigma$:

**Corollary 4.16** For any $\xi \in \mathbb{X} \setminus (\mathcal{D} \cup \mathcal{B})$ and $\mathbf{d} \in \mathcal{M}_\Omega$, and any $t_0$ such that $x(t, \xi, \mathbf{d}) \notin \mathcal{D} \cup \mathcal{B}$ for all $t \in [0, t_0]$, the following dissipation inequality holds:

$$V_0(x(t, \xi, \mathbf{d})) - V_0(\xi) \leq -\int_0^{t_0} \Xi(|x(t, \xi, \mathbf{d})|) \, dt. \quad (69)$$

*Proof.* Pick an initial state $\xi \in \mathbb{X} \setminus (\mathcal{D} \cup \mathcal{B})$, a disturbance $\mathbf{d} \in \mathcal{M}_\Omega$, and an appropriate $t_0$. Then

$$\begin{aligned}
V_0(x(t_0, &\xi, \mathbf{d})) - V_0(x(\xi)) \\
&= V_0(z(\sigma_{\xi, \mathbf{d}}(t_0), \xi, \mathbf{d} \circ \sigma_{\xi, \mathbf{d}}^{-1})) - V_0(z(\sigma_{\xi, \mathbf{d}}(t_1), \xi, \mathbf{d} \circ \sigma_{\xi, \mathbf{d}}^{-1})) \\
&\leq -\int_0^{\sigma_{\xi, \mathbf{d}}(t_0)} \Xi(|z(s, \xi, \mathbf{d} \circ \sigma_{\xi, \mathbf{d}}^{-1})|) \, ds \\
&= -\int_0^{t_0} \Xi(|z(\sigma_{\xi, \mathbf{d}}(t), \xi, \mathbf{d} \circ \sigma_{\xi, \mathbf{d}}^{-1})|) \, d\sigma_{\xi, \mathbf{d}}(t) \\
&= -\int_0^{t_0} \Xi(|x(t, \xi, \mathbf{d})|) \frac{d}{dt} \sigma_{\xi, \mathbf{d}}(t) \, dt \\
&= -\int_0^{t_0} \Xi(|x(t, \xi, \mathbf{d})|)[1 + |f(x(t, \xi, \mathbf{d}), \mathbf{d}(t))|^2 + \kappa(x(s, \xi, \mathbf{d}))] \, dt \\
&\leq -\int_0^{t_0} \Xi(|x(t, \xi, \mathbf{d})|) \, dt.
\end{aligned}$$

∎

## 4.4 Some definitions and facts from nonsmooth analysis

**Definition 4.17** A vector $\zeta \in \mathbb{R}^n$ is a *proximal subgradient* (respectively, *proximal supergradient*) of the function $V : \mathbb{R}^n \to (-\infty, +\infty]$ at $x$, if there exists some positive $\sigma$ such that, for all $x'$ in some neighborhood of $x$,

$$V(x') \geq V(x) + \zeta \cdot (x' - x) - \sigma |x' - x|^2 \tag{70}$$

$$(\text{correspondingly,} \quad V(x') \leq V(x) + \zeta \cdot (x' - x) + \sigma |x' - x|^2). \tag{71}$$

The (possibly empty) set of all proximal subgradients (respectively, supergradients) of $V$ at $x$ is called the *proximal subdifferential* and is denoted $\partial_P V(x)$ (respectively, *proximal superdifferential*, denoted $\partial^P V(x)$). Note that the definitions imply that if the function $V$ is differentiable at $x$, then both the subdifferential and superdifferential sets must be subsets of the $\{\nabla V(x)\}$. □

**Lemma 4.18** Let $\Sigma$ be a system of type (17), $\xi$ be a vector in $\mathbb{X}$, $\bar{d} \in \Omega$ and $V : \mathbb{X} \to \mathbb{R}$. Then, if there exist some continuous $\alpha_\xi : \mathbb{R}_{\geq 0} \to \mathbb{R}$ and $\varepsilon > 0$ such that the following inequality holds for all $\tau < \varepsilon$:

$$V(x(\tau, \xi, \mathbf{d})) - V(\xi) \leq \int_0^\tau \alpha_\xi(t) dt \tag{72}$$

(where $\mathbf{d}$ is the constant disturbance equal to $\bar{d}$), then for any $\zeta \in \partial_P V(\xi)$ the proximal form of inequality (72) holds:

$$\zeta \cdot f(\xi, \bar{d}) \leq \alpha_\xi(0). \tag{73}$$

*Proof.* It follows from (70) and (72) that, for all $\tau$ close enough to 0 we have

$$\int_0^\tau \alpha_\xi(t) dt \geq V(x(\tau, \xi, \mathbf{d})) - V(\xi) \geq \zeta \cdot (x(\tau, \xi, \mathbf{d}) - \xi) - \sigma |x(\tau, \xi, \mathbf{d}) - \xi|^2.$$

Dividing by $\tau$ and passing to the limit as $\tau$ tends to 0, we get (73). ■

## 4.5 Smoothing out a continuous Lyapunov function

The next result shows how to approximate a continuous function $V$ by a locally Lipschitz one in a weak $C^1$ sense. The function $V$ is assumed to be bounded below, or up to a translation by a constant, nonnegative.

**Lemma 4.19** Let $\Sigma : \dot{x} = f(x, d)$ be a system, with $x \in \mathbb{X} = \mathbb{R}^n$, and $d \in \Omega$, a compact metric space, so that $f(x, d)$ is locally Lipschitz in $x$ uniformly on $d$ and jointly continuous in $x$ and $d$. Assume that we are given:

- an open subset $\mathcal{O}$ of $\mathbb{X}$;

- a continuous, nonnegative function $V : \mathcal{O} \to \mathbb{R}$ satisfying

$$\zeta \cdot f(x, d) \leq \Theta(x, d) \quad \forall x \in \mathcal{O}, \ \zeta \in \partial_P V(x), \ d \in \Omega \tag{74}$$

with some continuous function $\Theta : \mathcal{O} \times \Omega \to \mathbb{R}$;

- two positive, continuous functions $\Upsilon_1$ and $\Upsilon_2$ on $\mathcal{O}$.

Then there exists a function $\widetilde{V} : \mathcal{O} \to \mathbb{R}$, locally Lipschitz on $\mathcal{O}$, such that

$$0 \leq V(x) - \widetilde{V}(x) \leq \Upsilon_1(x) \quad \forall x \in \mathcal{O}, \tag{75}$$

and

$$L_{f_d}\widetilde{V}(x) \leq \Theta(x,d) + \Upsilon_2(x) \quad \text{for all } d \in \Omega \text{ and almost all } x \in \mathcal{O}, \tag{76}$$

where $f_d$ is the vector field defined by $f_d = f(\cdot, d)$. $\square$

Note that, by Rademacher's theorem, the directional derivative $L_{f_d}\widetilde{V}(x) = \nabla V(x) \cdot f(x,d)$ is defined for almost all $x$, because $\widetilde{V}$ is locally Lipschitz. The proof will follow closely along the lines of the proof of the similar result for CLF's, found in [8] or [36].

We will first prove a "local" version of the result. For any $K$ which is a compact subset of $\mathcal{O}$ and $r > 0$, we introduce the following notations:

- $\bar{B}_r(K) := \{x \in \mathbb{X} : \exists \xi \in K \text{ with } |x - \xi| \leq r\}$ – the closed $r$-fattening of $K$, for $r > 0$.

- $\beta(K) := \sup_{x \in K} V(x)$,

- $m_K := \frac{1}{4} \min \{\Upsilon_2(x) : x \in K\}$,

- $\ell_K :=$ a Lipschitz constant for $f$ with respect to $x$ in $K$, that is:

$$|f(x,d) - f(x',d)| \leq \ell_K |x - x'|, \quad \forall d \in \Omega, \forall x \in K,$$

- $\omega_K(\cdot) :=$ the modulus of continuity of $V$ on $K$, that is:

$$\omega_K(\delta) := \sup \{V(x) - V(x') : |x - x'| \leq \delta; x, x' \in K\},$$

- $\pi_K(\cdot) :=$ the modulus of continuity of $\Theta$ on $K \times \Omega$, that is:

$$\pi_K(\delta) := \sup \{\Theta(x,d) - \Theta(x',d') : |x - x'| + \text{dist}(d,d') \leq \delta; x, x' \in K, d, d' \in \Omega\}.$$

To approximate the given continuous function $V$ by a locally Lipschitz one, we would like to use the notion of "Iosida-Moreau inf-convolution", well known in convex analysis. Fix a parameter $\alpha \in (0, 1]$. Suppose for the moment that $V$ is defined on the whole $\mathbb{X}$. Define

$$V_\alpha(x) := \min_{y \in \mathbb{X}} \left[ V(y) + \frac{1}{2\alpha^2} |y - x|^2 \right].$$

For each fixed $x$, the set of points $y$ where the minimum is attained is nonempty, because $V$ is bounded below. Denote one of them by $y_\alpha(x)$.

Fix a compact $K$ and $\alpha \in (0, 1]$, let

$$K_\alpha := \bar{B}_{\alpha\sqrt{2\beta(K)}}(K).$$

The following four claims summarize some of the useful properties of $V_\alpha$, proven in [8], [36] (or see the primary sources such as [9]):

*Claim 1:* For all $x \in K$,
$$|y_\alpha(x) - x|^2 \leq \min\left\{2\alpha^2 \beta(K),\ 2\alpha^2 \omega_{K_\alpha}\left(\alpha\sqrt{2\beta(K)}\right)\right\}.$$

*Proof of Claim 1:* By definition of $V_\alpha$ and $\beta(K)$, we have
$$\frac{1}{2\alpha^2}|y_\alpha(x) - x|^2 \leq V(x) - V(y_\alpha(x)) \leq V(x) \leq \beta(K), \tag{77}$$

so that $|y_\alpha(x) - x|^2 \leq 2\alpha^2 \beta(K)$. On the other hand, the first inequality in (77) implies also that
$$\begin{aligned}
|y_\alpha(x) - x|^2 &\leq 2\alpha^2(V(x) - V(y_\alpha(x))) \\
&\leq 2\alpha^2 \omega_{K_\alpha}(|y_\alpha(x) - x|) \leq 2\alpha^2 \omega_{K_\alpha}\left(\alpha\sqrt{2\beta(K)}\right),
\end{aligned}$$

proving the claim.

Let
$$\zeta_\alpha(x) = \frac{x - y_\alpha(x)}{\alpha^2}.$$

*Claim 2:* For any $x \in \mathbb{X}$,
$$\zeta_\alpha(x) \in \partial_P V(y_\alpha(x)), \quad \text{and} \tag{78}$$
$$\zeta_\alpha(x) \in \partial^P V_\alpha(x), \tag{79}$$

*Claim 3:* For any $x \in K$,
$$V_\alpha(x) \leq V(x) \leq V_\alpha(x) + \omega_{K_\alpha}\left(\alpha\sqrt{2\beta(K)}\right).$$

*Claim 4:* $V_\alpha$ is locally Lipschitz.

Now recall that, in our setting, $V$ is defined on an open subset $\mathcal{O}$ of $\mathbb{X}$, so, we can't minimize the expression in the definition of $V_\alpha$ over the whole state space. However, for any compact subset $K$ of $\mathcal{O}$ we can choose $\alpha$ small enough so that
$$K_\alpha \subset \mathcal{O},$$

hence, for any $x$ in $K$ we could define $V_\alpha$ minimizing over $y \in \mathcal{O}$, and the same function $V_\alpha$ results on $K$.

**Lemma 4.20** Assume that $V$ satisfies (74), a compact $K$ is fixed, and $\alpha \in (0,1]$ is chosen to satisfy

- $K_\alpha = \bar{B}_{\alpha\sqrt{2\beta(K)}}(K) \subset \mathcal{O}$,
- $\pi_{K_\alpha}\left(\alpha\sqrt{2\beta(K)}\right) \leq m_K$, and
- $2\ell_{K_\alpha}\omega_{K_\alpha}\left(\alpha\sqrt{2\beta(K)}\right) \leq m_K$.

Then the function
$$V_\alpha(x) := \inf_{y \in \mathcal{O}} \left[ V(y) + \frac{1}{2\alpha^2} |y - x|^2 \right] \tag{80}$$

will possess the following property:

for all $x \in K$, all $d \in \Omega$, and all $\zeta \in \partial_P V_\alpha(x)$,
$$\zeta \cdot f(x, d) \leq \Theta(x, d) + 2m_K.$$

*Proof.* Fix any $x \in K$. The choice of $\alpha$ insures that the infimum in (80) is a minimum, and it is achieved at some $y_\alpha(x) \in K_\alpha$. The definition of $\zeta_\alpha(x)$ and Claim 1 imply that

$$|\zeta_\alpha(x)||y_\alpha(x) - x| = \frac{|y_\alpha(x) - x|^2}{\alpha^2} \leq 2\omega_{K_\alpha} \left( \alpha \sqrt{2\beta(K)} \right). \tag{81}$$

Using again the fact that $y_\alpha(x) \in K_\alpha$,

$$|f(x, d) - f(y_\alpha(x), d)| \leq \ell_{K_\alpha} |x - y_\alpha(x)|.$$

Combining the last inequality with (81) we obtain

$$|\zeta_\alpha(x)||f(x, d) - f(y_\alpha(x), d)| \leq 2\ell_{K_\alpha} \omega_{K_\alpha} \left( \alpha \sqrt{2\beta(K)} \right). \tag{82}$$

Now, by Claim 2, $\zeta_\alpha(x) \in \partial_P V(y_\alpha(x))$. Hence,

$$\begin{aligned}
\zeta_\alpha(x) \cdot f(y_\alpha(x), d) &\leq \Theta(y_\alpha(x), d) \\
&\leq \Theta(x, d) + \pi_{K_\alpha}(|x - y_\alpha(x)|) \\
&\leq \Theta(x, d) + \pi_{K_\alpha}\left( \alpha \sqrt{2\beta(K)} \right) \\
&\leq \Theta(x, d) + m_K.
\end{aligned}$$

So, by the last inequality, (82), and the choice of $\alpha$, we have

$$\begin{aligned}
\zeta_\alpha(x) \cdot f(x, d) &= \zeta_\alpha(x) \cdot f(y_\alpha(x), d) + \zeta_\alpha(x) \cdot (f(x, d) - f(y_\alpha(x), d)) \\
&\leq \Theta(x, d) + m_K + m_K.
\end{aligned} \tag{83}$$

Next, we show that $\partial_P V_\alpha(x) \subseteq \{\zeta_\alpha\}$ for all $x \in K$. Pick any $\zeta \in \partial_P V_\alpha(x)$. By definition of the proximal subgradient, the inequality

$$\zeta \cdot (y - x) \leq V_\alpha(y) - V_\alpha(x) + o(|y - x|) \tag{84}$$

holds for all $y$ near $x$. Since $\zeta_\alpha(x)$ is a proximal supergradient of $V_\alpha$ at $x$, we also have

$$-\zeta_\alpha(x) \cdot (y - x) \leq -V_\alpha(y) + V_\alpha(x) + o(|y - x|) \tag{85}$$

Adding (84) and (85), we get

$$(\zeta - \zeta_\alpha(x)) \cdot (y - x) \leq o(|y - x|) \tag{86}$$

for all $y$ sufficiently close to $x$. Substituting $y = x + h(\zeta - \zeta_\alpha(x))$ in (86) and letting $h$ tend to 0, we arrive at $\zeta = \zeta_\alpha(x)$. ∎

Now we are ready to prove the main lemma of the section.

*Proof.* (of Lemma 4.19) For every $x \in \mathcal{O}$ find an $r_x > 0$ small enough so that $\bar{B}_{r_x}(x) \subset \mathcal{O}$. Then the collection of open balls $\{B_{r_x}(x), x \in \mathcal{O}\}$ forms an open covering of $\mathcal{O}$. By paracompactness of $\mathcal{O}$ we can find a locally finite refinement $\{B_i, i \in \mathbb{N}\}$ of $\{B_{r_x}(x), x \in \mathcal{O}\}$ (c.f. [6, Lemma 4.1]). Moreover, since $\cup_{x \in \mathcal{O}} B_{r_x} = \mathcal{O}$, we also have $\cup_i B_i = \mathcal{O}$. Let $\{\varphi_i, i \in \mathbb{N}\}$ be a partition of unity, subordinate to $\{B_i\}$. For each index $i$, let

$$\mathcal{J}_i = \{j \in \mathbb{N} \,:\, B_i \cap B_j \neq \emptyset\}.$$

Notice that $\mathcal{J}_i$ is finite for all $i$, because of local finiteness of the covering $\{B_i\}$. For every $i \in \mathbb{N}$ define

$$M_i := \sup |\nabla \varphi_i(x)| |f(x, d)|,$$

where the supremum is taken over all $x \in \cup_{j \in \mathcal{J}_i} \bar{B}_j$ and all $d \in \Omega$; and

$$N_i = \max_{j \in \mathcal{J}_i} \operatorname{card}(\mathcal{J}_j),$$

that is $N_i$ denotes the maximum cardinality of $\mathcal{J}_j$ for $j$ such that $B_j$ intersects $B_i$. Next, for each compact $K = \bar{B}_i$, $i \in \mathbb{N}$ choose an $\alpha_i$ as in the formulation of Lemma 4.20, satisfying the following two additional conditions:

$$N_i \omega_{K_{\alpha_i}}\left(\alpha_i \sqrt{2\beta(\bar{B}_i)}\right) M_i < \frac{1}{2} \inf_{x \in \bar{B}_i} \Upsilon_2(x) \tag{87}$$

and

$$\omega_{K_{\alpha_i}}\left(\alpha_i \sqrt{2\beta(\bar{B}_i)}\right) < \frac{1}{2} \inf_{x \in \bar{B}_i} \Upsilon_1(x). \tag{88}$$

(where $K_{\alpha_i}$ denotes $B_{\alpha_i \sqrt{2\beta(\bar{B}_i)}}(\bar{B}_i)$ ).

Next, for each $i \in \mathbb{N}$, define $V_{\alpha_i}$ on $\bar{B}_i$ as in (80) (this can be done, because, by the choice of $\alpha_i$, $K_{\alpha_i} \subseteq \mathcal{O}$). Since $V_{\alpha_i}$ is locally Lipschitz, it is differentiable almost everywhere by Rademacher's theorem. Lemma 4.20 implies that for almost all $x \in \bar{B}_i$ and all $d \in \Omega$

$$L_{f_d} V_{\alpha_i}(x) \leq \Theta(x, d) + \Upsilon_2(x)/2.$$

Define

$$\widetilde{V} := \sum_{i=1}^{+\infty} V_{\alpha_i} \varphi_i.$$

Strictly speaking, this does not make sense, because $V_{\alpha_i}$ is not defined outside $\bar{B}_i$, but that does not matter, because $\varphi_i$ vanishes outside $B_i$ anyway.

Since each $V_{\alpha_i}(x) \leq V(x)$ for all $x \in B_i$ and $\varphi_i$'s add up to 1, the definition of $\widetilde{V}$ shows that $\widetilde{V}(x) \leq V(x)$ for all $x \in \mathcal{O}$. It is also clear from the definition of $\widetilde{V}$ that $\widetilde{V}$ is locally Lipschitz on $\mathcal{O}$.

We claim that for almost all $x \in \mathcal{O}$ and all $d \in \Omega$ the following holds:

1. $0 \leq V(x) - \widetilde{V}(x) \leq \Upsilon_1(x)$;

2. $\nabla \widetilde{V}(x) \cdot f(x, d) \leq \Theta(x, d) + \Upsilon_2(x)$.

Take any $x \in \mathcal{O}$ and find $i \in \mathbb{N}$ such that $x \in B_i$. Define also $\mathcal{J}_x := \{j \in \mathbb{N} : x \in B_j\}$. Note that $\mathcal{J}_x \subseteq \mathcal{J}_i$ and that $\widetilde{V}(x) := \sum_{j \in \mathcal{J}_x} V_{\alpha_j}(x) \varphi_j(x)$. Then Claim 3 and the choice of $\alpha_j$ (condition (88)) imply

$$V(x) - V_{\alpha_j}(x) \leq \omega_{K_{\alpha_j}}\left(\alpha_j \sqrt{2\beta(\bar{B}_j)}\right) \leq \Upsilon_1(x)/2$$

for all $j \in \mathcal{J}_x$. Thus,

$$\begin{aligned}
V(x) - \widetilde{V}(x) &= V(x) - \sum_{j=1}^{+\infty} V_{\alpha_j} \varphi_j \\
&= V(x) - \sum_{j \in \mathcal{J}_i} V_{\alpha_j} \varphi_j \\
&\leq V(x) - \min_{j \in \mathcal{J}_x}\{V_{\alpha_j}(x)\} \left(\sum_{j \in \mathcal{J}_x} \varphi_j(x)\right) \\
&= V(x) - \min_{j \in \mathcal{J}_x}\{V_{\alpha_j}(x)\} \\
&\leq \frac{\Upsilon_1(x)}{2},
\end{aligned}$$

proving the first statement.

To prove the second statement, write

$$\begin{aligned}
L_{f_d}\widetilde{V}(x) &= \sum_{j \in \mathcal{J}_x} L_{f_d} V_{\alpha_j}(x)\, \varphi_j(x) + \sum_{j \in \mathcal{J}_x} V_{\alpha_j}(x)\, L_{f_d}\varphi_j(x) \\
&= \sum_{j \in \mathcal{J}_x} L_{f_d} V_{\alpha_j}(x)\, \varphi_j(x) + \sum_{j \in \mathcal{J}_x} V_{\alpha_j}(x)\, L_{f_d}\varphi_j(x) - V(x)\sum_{j \in \mathcal{J}_x} L_{f_d}\varphi_j(x) \quad (89) \\
&\leq \max_{j \in \mathcal{J}_x} L_{f_d} V_{\alpha_j}(x) \sum_{j \in \mathcal{J}_x} \varphi_j(x) + \sum_{j \in \mathcal{J}_x} \left(V_{\alpha_j}(x) - V(x)\right) L_{f_d}\varphi_j(x) \\
&\leq \Theta(x,d) + \Upsilon_2(x)/2 + \sum_{j \in \mathcal{J}_i} |V_{\alpha_j}(x) - V(x)|\, |L_{f_d}\varphi_j(x)| \\
&\leq \Theta(x,d) + \Upsilon_2(x)/2 + \sum_{j \in \mathcal{J}_i} M_j \omega_{K_{\alpha_j}}\left(\alpha_j\sqrt{2\beta(\bar{B}_j)}\right) \quad (90) \\
&\leq \Theta(x,d) + \Upsilon_2(x)/2 + \sum_{j \in \mathcal{J}_i} \frac{\Upsilon_2(x)}{2N_j} \quad (91) \\
&\leq \Theta(x,d) + \Upsilon_2(x)/2 + \Upsilon_2(x)/2 \quad (92) \\
&\leq \Theta(x,d) + \Upsilon_2(x)
\end{aligned}$$

where the equality (89) follows from the fact that $\sum_{j \in \mathcal{J}_x} L_{f_d}\varphi_j(x) = 0$, inequality (90) follows from Lemma 4.20, and Claim 3; inequality (91) follows from the choice of $\alpha_i$ (condition (87)); and (92) is implied by the fact that

$$\text{card}\mathcal{J}_i \leq N_j \quad \forall j \in \mathcal{J}_i.$$

This completes the proof of the lemma. ∎

The lemma we have just proved will provide the "continuous ⇒ locally Lipschitz away from zero" step in the smoothing process. To obtain a smooth Lyapunov function, we will use the following simple smoothing result. The proof is given, for the special case when $\alpha$ does not depend on $d$, in [26], but the general case ($\alpha$ depends on $d$) is proved in exactly the same manner, so we omit the proof here.

**Lemma 4.21** Let $\mathcal{O}$ be an open subset of $\mathbb{R}^n$, and let $\Omega$ be a compact subset of $\mathbb{R}^l$, and assume given:

- a locally Lipschitz function $\Phi : \mathcal{O} \to \mathbb{R}$;

- a continuous map $f : \mathbb{R}^n \times \Omega \to \mathbb{R}^n$, $(x, d) \mapsto f(x, d)$ which is locally Lipschitz on $x$ uniformly on $d$;

- a continuous function $\alpha : \mathcal{O} \times \Omega \to \mathbb{R}$ and continuous functions $\mu, \nu : \mathcal{O} \to \mathbb{R}_{>0}$ such that for each $d \in \Omega$,
$$L_{f_d}\Phi(\xi) \leq \alpha(\xi, d), \quad \text{a.e. } \xi \in \mathcal{O}, \tag{93}$$
where $f_d$ is the vector field defined by $f_d = f(\cdot, d)$ (recall that $\nabla \Phi$ is defined a.e., since $\Phi$ is locally Lipschitz by Rademacher's Theorem). Then there exists a smooth function $\Psi : \mathcal{O} \to \mathbb{R}$ such that
$$|\Phi(\xi) - \Psi(\xi)| < \mu(\xi), \quad \forall \xi \in \mathcal{O}$$
and for each $d \in \Omega$,
$$L_{f_d}\Psi(\xi) \leq \alpha(\xi, d) + \nu(\xi), \quad \forall \xi \in \mathcal{O}.$$
□

The next result immediately follows by Lemma 4.21.

**Corollary 4.22** Under the assumptions of Lemma 4.19 there also exists a smooth function $\hat{V}$ on $\mathcal{O}$, satisfying inequalities
$$\left|V(x) - \hat{V}(x)\right| \leq \Upsilon_1(x) \quad \forall x \in \mathcal{O} \tag{94}$$
and
$$L_{f_d}\hat{V}(x) \leq \Theta(x, d) + \Upsilon_2(x) \quad \text{for all } d \in \Omega \text{ and all } x \in \mathcal{O}. \tag{95}$$

*Proof.* Suppose that $\mathcal{O}$, $V$, $\Upsilon_1$, and $\Upsilon_2$ are given and the assumptions of Lemma 4.19 hold. Replacing $\Upsilon_1$ by $\Upsilon_1/2$ and $\Upsilon_2$ by $\Upsilon_2/2$, and applying Lemma 4.19, we can find a locally Lipschitz function $\widetilde{V}$, defined on $\mathcal{O}$ and satisfying (75) and (95) with $\Upsilon_1/2$ and $\Upsilon_2/2$ instead of $\Upsilon_1$ and $\Upsilon_2$. Next, Lemma 4.21, applied with the same $\mathcal{O}$ and with $\Phi := \widetilde{V}$, $\alpha := \Theta + \Upsilon_2/2$, $\mu := \Upsilon_1/2$, and $\nu := \Upsilon_2/4$ furnishes a smooth function $\hat{V} := \Psi$ as needed. ∎

In section 4.3 we have constructed a function $V_0$, satisfying inequalities (67) on $\mathcal{E}_1$ and (69) along all trajectories of $\Sigma$, contained in $\mathcal{E}_1$. Therefore, by Lemma 4.18, the proximal inequality (73) holds for $V_0$ at any interior point of $\mathcal{E}_1$. Then Corollary 4.22, applied with $\mathcal{O} := \text{int } \mathcal{E}_1$, $\Upsilon_1(\cdot) := \underline{\alpha}(|\cdot|)/2$, $\Upsilon_2(\cdot) := \Xi(|\cdot|)/2$, $\Theta(x, d) := \Xi(|x|)$, provides a smooth function
$$V_1 : \text{int } \mathcal{E}_1 \to \mathbb{R}_{>0}; \quad V_1 := \hat{V}_0$$

satisfying the following two conditions for all $\xi \in \text{int } \mathcal{E}_1$:

$$\underline{\alpha}(|\xi|)/2 \leq V_1(\xi) \leq \bar{\alpha}(|\xi|) + \underline{\alpha}(|\xi|)/2$$

(we will replace $\underline{\alpha}(|\cdot|)/2$ by $\underline{\alpha}(\cdot)$ and $\bar{\alpha}(|\cdot|) + \underline{\alpha}(|\cdot|)/2$ by $\bar{\alpha}(|\cdot|)$ from now on to avoid cluttering the notation.)

$$L_{f_d} V_1(\xi) \leq -\Xi_1(|\xi|) \quad \forall d \in \Omega, \tag{96}$$

where $\Xi_1(\cdot) \equiv \Xi(\cdot)/2$.

## 4.6 Extending to the rest of $\mathbb{X}$ and smoothing at the origin

To construct a UOSS-Lyapunov-like function defined on the whole $\mathbb{X}$, we must "patch" $V_1$ with some smooth, proper, and positive definite function such that the dissipation inequality still holds.

**Lemma 4.23** Suppose $\Sigma$ is a system of type (17), and $\rho$ is a function of class $\mathcal{K}_\infty$. Define $\mathcal{E}_1 := \{x \in \mathbb{X} : |x| > 2\rho(|h(x)|)\}$ and suppose that $V_1 : \mathcal{E}_1 \to \mathbb{R}_{\geq 0}$ is a smooth function satisfying, with some suitable $\mathcal{K}_\infty$ functions, the inequality

$$\underline{\alpha}(|\xi|) \leq V_1(\xi) \leq \bar{\alpha}(|\xi|) \tag{97}$$

and inequality (96) on $\mathcal{E}_1$. Then there exist a Lyapunov-like function $V_2$ for $\Sigma$, smooth away from the origin, and a class $\mathcal{K}_\infty$ function $\Phi$, such that

$$V_2(x) = \Phi \circ V_1(x) \quad \forall x \text{ such that } |x| > 3\rho(|h(x)|),$$

and the following dissipation inequality holds with some $\check{\alpha}_3 \in \mathcal{K}_\infty$, $\check{\gamma} \in \mathcal{K}$

$$\nabla V_2(\xi) \cdot f(\xi, d) \leq -\check{\alpha}_3(|\xi|) + \check{\gamma}(|h(\xi)|) \quad \forall \xi \neq 0, \quad \forall d \in \Omega. \tag{98}$$

*Proof.* Let

$$\mathcal{E}_2 := \{\xi : |\xi| > 3\rho(|h(\xi)|)\}$$

Since the sets $\{\xi : |\xi| \geq 3\rho(|h(\xi)|)\}$ and $\{\xi : |\xi| \leq 2\rho(|h(\xi)|)\}$ are disjoint and closed in the topology of $\mathbb{X} \setminus \{0\}$, one can find a smooth function $\phi : \mathbb{X} \setminus \{0\} \to [0,1]$ with the property that

$$\phi(\xi) = \begin{cases} 1, & \text{if } |\xi| \geq 3\rho(|h(\xi)|) \\ 0, & \text{if } |\xi| \leq 2\rho(|h(\xi)|) \end{cases}$$

and $\phi$ is non-zero elsewhere. It is easy to see that $|\nabla \phi(x)|$ is bounded above by a $\mathcal{K}$ function $\nu_2$ of $|x|$ outside the unit ball centered at 0. One can also find a smooth, strictly increasing function $\nu_1 : [0,1] \to \mathbb{R}_{\geq 0}$, such that $\nu_1(0) = 0$ and $|\nabla \phi(x)| \leq \frac{1}{\nu_1(|x|)}$ for all $x$ such that $0 < |x| \leq 1$.

Let $\nu_3$ be a $\mathcal{K}$-function such that $\max_{d \in \Omega} |f(\xi, d)| \leq \nu_3(|\xi|)$. Take any smooth function $\pi_1 : [0, \bar{\alpha}(1)] \to \mathbb{R}_{\geq 0}$ with $\pi_1'(s) > 0$ for all $s \in (0, \bar{\alpha}(1))$, such that

$$\frac{\pi_1(\bar{\alpha}(r))}{\nu_1(r)} < s(r), \quad \frac{\pi_1(r^2)}{\nu_1(r)} < s(r)$$

for some $\mathcal{K}$-function $s$ and all $0 < r \leq 1$. Let $\pi_2$ be any $\mathcal{K}$-function such that $\pi_2(r) \leq \pi_1'(r)$ for all nonnegative $r \leq 1$.

Let $\Phi(r) = \int_0^r \pi_2(r_1) dr_1$. Then $\Phi(r) \leq \pi_1(r)$ for all $r \leq 1$, so that $\frac{\Phi(\bar{\alpha}(r))}{\nu_1(r)} < s(r)$ and $\frac{\Phi(r^2)}{\nu_1(r)} < s(r)$ for all $r \in (0, 1]$.

Now let
$$V_2(\xi) = \phi(\xi)\Phi(V_1(\xi)) + (1-\phi(\xi))\Phi(|\xi|^2). \tag{99}$$

If $|\xi| > 3\rho(|h(\xi)|)$, then $V_2 \equiv \Phi \circ V_1$ in a neighborhood of $\xi$, so that, for all $d \in \Omega$,
$$\nabla V_2(\xi) \cdot f(\xi, d) = \Phi'(V_1(\xi))[\nabla V_1(\xi) \cdot f(\xi, d)] \leq -\pi_2(\underline{\alpha}(|\xi|))\Xi_1(|\xi|). \tag{100}$$

On the other hand, if $|\xi| \leq 3\rho(|h(\xi)|)$, then

$\nabla V_2(\xi) \cdot f(\xi, d)$
$= [\nabla \phi(\xi) \cdot f(\xi, d)] \Phi(V_1(\xi)) + \phi(\xi) \Phi'(V_1(\xi)) [\nabla V_1(\xi) \cdot f(\xi, d)]$
$\quad - [\nabla \phi(\xi) \cdot f(\xi, d)] \Phi(|\xi|^2) + (1-\phi(\xi)) 2\Phi'(|\xi|^2) [\xi \cdot f(\xi, d)]$
$\leq |\nabla \phi(\xi)| |f(\xi, d)| \Phi(V_1(\xi)) + |\nabla \phi(\xi)| |f(\xi, d)| \Phi(|\xi|^2) + 2|(1-\phi(\xi))| \Phi'(|\xi|^2) |\xi| |f(\xi, d)|$
$\leq |\nabla \phi(\xi)| |f(\xi, d)| \Phi(\bar{\alpha}(|\xi|)) + |\nabla \phi(\xi)| |f(\xi, d)| \Phi(|\xi|^2) + 2\Phi'(|\xi|^2) |\xi| |f(\xi, d)|,$

where the first inequality follows from the fact that $\phi(\xi) \Phi'(V_1(\xi)) [\nabla V_1(\xi) \cdot f(\xi, d)] \leq 0$. Next, the definition of $\Phi$ provides the following bounds for the three terms in the right-hand side of the last inequality:

$$\Phi(\bar{\alpha}(|\xi|)) |\nabla \phi(\xi)| |f(\xi, d)| \leq \max\{s(|\xi|), \Phi(\bar{\alpha}(|\xi|))\nu_2(|\xi|)\} \nu_3(|\xi|),$$

$$2\Phi'(|\xi|^2)[\xi \cdot f(\xi, d)] \leq 2\pi_2(|\xi|^2) |\xi| \nu_3(|\xi|),$$

$$\Phi(|\xi|^2) |\nabla \phi(\xi)| |f(\xi, d)| \leq \max\{s(|\xi|), \Phi(|\xi|^2)\nu_2(|\xi|)\} \nu_3(|\xi|).$$

Define $\check{\alpha}_3$, $\check{\gamma}$, $\check{\alpha}_1$ and $\check{\alpha}_2$ by
$$\check{\alpha}_3(r) := \pi_2(\underline{\alpha}(r))\Xi_1(r),$$

$$\check{\gamma}(r) := 2\pi_2((3\rho(r))^2)3\rho(r)\nu_3(3\rho(r))$$
$$\quad + \max\{s(3\rho(r)), \Phi((3\rho(r))^2)\nu_2(3\rho(r))\} \nu_3(3\rho(r)) + \check{\alpha}_3(3\rho(r)),$$

$$\check{\alpha}_1(r) = \min\{\Phi(r^2), \Phi \circ \underline{\alpha}(r)\},$$

and
$$\check{\alpha}_2(r) = \max\{\Phi(r^2), \Phi \circ \bar{\alpha}(r)\}.$$

Then inequalities
$$\check{\alpha}_1(|x|) \leq V_2(x) \leq \check{\alpha}_2(|x|) \quad \forall \, x \neq 0 \tag{101}$$

and (98) hold for $V_2$. Define $V_2(0) := 0$. Then inequality (7) holds for $V_2$ with $\alpha_1 := \check{\alpha}_1$, $\alpha_2 := \check{\alpha}_2$ on $\mathbb{X}$, and in particular implies that $V_2$ is continuous as 0. Then, $V_2$ is a UOSS Lyapunov-like function for $\Sigma$, smooth away from the origin. ∎

Recall that $V_2(\xi) \equiv \Phi(V_1(\xi))$ for all $\xi$ with $|\xi| > 3\rho(|h(\xi)|)$. Therefore inequality (96) implies that
$$\nabla V_2 \cdot f(\xi, d) \leq -\alpha_3(|\xi|) \text{ for all } |\xi| > 3\rho(|h(\xi)|), \text{ and all } d \in \Omega. \tag{102}$$

**Proposition 4.24** Suppose that a system $\Sigma$ of type (17) admits a continuous UOSS Lyapunov-like function $V_2$, smooth away from 0 and satisfying inequalities (101), (98), and (102). Then $\Sigma$ admits a UOSS-Lyapunov function. □

The basic idea used to obtain a Lyapunov function, smooth on the whole $\mathbb{X}$, is composing $V_2$ with some appropriately chosen $\mathcal{K}_\infty$-function $\beta$. This technique was previously utilized in [26]. We will need the following generalization of Lemma 4.3 from [26], where this function $\beta$ is constructed. In our setting we also need the derivative of $\beta$ to be of class $\mathcal{K}_\infty$, which was not required in [26]. The proof is only a slight modification of the proof of the mentioned lemma.

**Lemma 4.25** Assume that $V : \mathbb{R}^n \longrightarrow \mathbb{R}_{\geq 0}$ is $C^0$, positive definite, and the restriction $V|_{\mathbb{R}^n \setminus \{0\}}$ is $C^\infty$.

Then, given any $m \in \mathbb{N}$ there exists a $\mathcal{K}_\infty$-function $\beta_m$, smooth on $(0, \infty)$, satisfying the following conditions:

- $\beta_m^{(i)}(t) \to 0$ as $t \to 0^+$ for each $i = 0, 1, \ldots$,

- $\beta_m^{(i)} \in \mathcal{K}_\infty \;\forall i \leq m$,

- $W_m := \beta \circ V$ is a $C^\infty$ function on all of $\mathbb{R}^n$.

□

We now return to the proof of Proposition 4.24.

*Proof.* Take the function $V_2$, and apply Lemma 4.25 to get a function $\beta_1$, with derivative in class $\mathcal{K}$, such that
$$V_3 := \beta_1 \circ V_2$$
is smooth on $\mathbb{X}$ ($\mathcal{A} = \{0\}$ in our case). Let $\alpha_1(\cdot) := \beta_1(\check{\alpha}_1(\cdot)/2)$, $\alpha_2(\cdot) := \beta_1(\check{\alpha}_2(\cdot) + \check{\alpha}_1(\cdot)/2)$. Then $\alpha_i \in \mathcal{K}_\infty$ and,
$$\alpha_1(|x|) \leq V_3(x) \leq \alpha_2(|x|).$$

Furthermore, for any $x \in \mathbb{X} \setminus \{0\}$ we have
$$\begin{aligned}
L_{f_d} V_3(x) &= \beta_1'(V_2(x)) \left[ L_{f_d} V_2 \right](x) \\
&\leq \beta_1'(V_2(x)) \left( -\check{\alpha}_3(|x|)/2 + \check{\gamma}(|h(x)|) \right) \\
&\leq -\beta_1'(\check{\alpha}_1(|x|)/2)\check{\alpha}_3(|x|)/2 + \beta_1'(\check{\alpha}_2(|x|) + \check{\alpha}_1(|x|)/2)\check{\gamma}(|h(x)|).
\end{aligned}$$

Recall that, because of (102), the $\check{\gamma}$ term in the last estimate can be dropped if $|x| > 3\rho(|h(x)|)$, so that we can write
$$L_{f_d} V_3(x) \leq -\beta_1'(\check{\alpha}_1(|x|)/2)\check{\alpha}_3(|x|)/2 + \beta_1'(\check{\alpha}_2(3\rho(|h(x)|)) + \check{\alpha}_1(3\rho(|h(x)|))/2)\check{\gamma}(|h(x)|)$$

Thus $V_3$ is a smooth UOSS-Lyapunov function, satisfying the dissipation inequality with
$$\alpha_3(\cdot) := \beta_1'(\check{\alpha}_1(\cdot)/2)\,\check{\alpha}_3(\cdot)/2$$
and
$$\gamma(\cdot) := \beta_1'(\check{\alpha}_2(3\rho(\cdot)) + \check{\alpha}_1(3\rho(\cdot)))/2)\,\check{\gamma}(\cdot).$$
This completes the construction. ■

# 5 Norm-observers

As conjectured in [47] and proved in this presentation, every IOSS system admits an IOSS-Lyapunov function. One of the main motivations for the notion of IOSS, and for deriving Lyapunov characterizations, is the fact that a Lyapunov function enables us to get insights into the behavior of control systems. In particular, it may be useful to have an estimate of how far the system is from the equilibrium at any given time, and in some situations this "norm-estimate" is sufficient for the design of a stabilizer. We next provide a construction for a norm-observer in the most general case – for systems of type (5), assuming that we have a smooth UIOSS-Lyapunov function at our disposal.

## 5.1 Exponential decay Lyapunov functions

**Definition 5.1** Let $\Sigma$ be a system of type (5). A $C^1$ function $V : \mathbb{X} \to \mathbb{R}_{\geq 0}$ is an *exponential decay UIOSS-Lyapunov function* for $\Sigma$ if it satisfies (7) with some $\alpha_1$ and $\alpha_2$, and the following version of inequality (8):

$$\nabla V(x) \cdot f(x,u,w) \leq -V(x) + \sigma_1(|u|) + \sigma_2(|h(x)|) \quad \forall x \in \mathbb{X},\ u \in \mathbb{U},\ w \in \Gamma \qquad (103)$$

holds with some $\sigma_1$ and $\sigma_2 \in \mathcal{K}$.

**Lemma 5.2** Suppose $V$ is an UIOSS-Lyapunov function for a system $\Sigma$ of type (5), satisfying inequality (8). Then there exists a $\mathcal{K}_\infty$-function $\rho$ such that a function $W := \rho \circ V$ is an exponential decay UIOSS-Lyapunov function for $\Sigma$.

*Proof.* Assume that system (5) admits an UIOSS-Lyapunov function with $\alpha_i$ ($i=1,2$) as in (7) and with $\alpha, \sigma_1, \sigma_2$ as in (8). Replacing $\alpha$ by $\alpha \circ \alpha_2^{-1}$, we have:

$$\frac{d}{dt}V(x(t,\xi,\mathbf{u})) \leq -\alpha(V(x(t,\xi,\mathbf{u}))) + \sigma_1(|\mathbf{u}(t)|) + \sigma_2(|y(t,\xi,\mathbf{u})|) \qquad (104)$$

for almost all $t \in [0,\ t_{\max}(\xi,u))$. According to Lemma 12 in [33], there exists some function $\rho \in \mathcal{K}_\infty$ which can be extended as a $C^1$ function to a neighborhood of $[0,\infty)$ such that $\rho'(r)\frac{\alpha(r)}{2} \geq \rho(r)$ for all $r \geq 0$. Consider the function $W(\xi) := \rho(V(\xi))$. Observe that $W$ is again proper and positive definite. Along any trajectory $x(t) := x(t,\xi,\mathbf{u})$ (with $y(t) := y(t,\xi,\mathbf{u})$), at any point where (8) holds, one has that $\frac{d}{dt}W(x(t)) = \rho'(V(x(t)))\frac{d}{dt}V(x(t))$ is upper bounded by:

$$-\rho'(V(x(t)))\frac{\alpha(V(x(t)))}{2} + \rho'(V(x(t)))\left(-\frac{\alpha(V(x(t)))}{2} + \sigma_1(|\mathbf{u}(t)|) + \sigma_2(|y(t)|)\right)$$

which in turn is bounded by

$$-\rho(V(x(t)) + \rho'(V(x(t)))\left(-\frac{\alpha(V(x(t)))}{2} + \sigma_1(|\mathbf{u}(t)|) + \sigma_2(|y(t)|)\right). \qquad (105)$$

Observe that when $V(x(t)) \geq \alpha^{-1}(2\sigma_1(|\mathbf{u}(t)|) + 2\sigma_2(|y(t)|))$ it holds that:

$$\rho'(V(x(t)))\left(-\frac{\alpha(V(x(t)))}{2} + \sigma_1(|\mathbf{u}(t)|) + \sigma_2(|y(t)|)\right) \leq 0, \qquad (106)$$

while if instead $V(x(t)) \leq \alpha^{-1}(2\sigma_1(|\mathbf{u}(t)|) + 2\sigma_2(|y(t)|))$, then:

$$\rho'(V(x(t)))\left(\sigma_1(|\mathbf{u}(t)|) + \sigma_2(|y(t)|)\right) \leq \hat{\sigma_1}(|\mathbf{u}(t)|) + \hat{\sigma_2}(|y(t)|) \tag{107}$$

for some $\mathcal{K}_\infty$-functions $\hat{\sigma}_1$ and $\hat{\sigma}_2$ (using here the fact that $\rho'(s)$ is a continuous function). Combining (106) and (107), one concludes from the estimate (105) on $\frac{d}{dt}W(x(t))$ that

$$\frac{d}{dt}W(x(t)) \leq -W(x(t)) + \hat{\sigma_1}(|\mathbf{u}(t)|) + \hat{\sigma_2}(|y(t)|)$$

for almost all $t \in [0, t_{\max})$. ∎

## 5.2 Construction of a norm observer

**Proposition 5.3** Suppose that a system $\Sigma$ admits an exponential decay UIOSS Lyapunov function $V$. Then the pair $(\Sigma_{n.o}, k)$, where

$$\Sigma_{n.o}: \quad \dot{p} = -p + \sigma_1(|u|) + \sigma_2(|y|), \tag{108}$$

with $\sigma_1$ and $\sigma_2$ as in (103) and $k(\cdot, \cdot)$ defined by $k(s, r) = s$, is a norm-estimator for $\Sigma$.

*Proof.* Assume without loss of generality that the function $\alpha_2$ in the definition of $V$ satisfies $r \leq \alpha_2(r)$ for all nonnegative $r$.

The system (108) is ISS with respect to $u$ and $y$, since it can be seen as an asymptotically stable linear system driven by the input $(\sigma_1(|u|), \sigma_2(|y|))$, so, inequality (11) obviously holds. Pick any initial states $\xi$, $\zeta$ of $\Sigma$ and (108) respectively, any control $\mathbf{u}$, and any disturbance $\mathbf{w}$. Consider the resulting trajectory $(x(t), p(t))$ of the composite system. Property (103) implies that

$$\frac{d}{dt}(V(x(t)) - p(t)) \leq -(V(x(t)) - p(t)) \tag{109}$$

for almost all $t \in [0, t_{\max}(\xi, \mathbf{u}, \mathbf{w}))$. Thus

$$V(x(t)) \leq p(t) + e^{-t}(V(\xi) - \zeta) \leq |p(t)| + 2e^{-t}\alpha_2(|\xi| + |\zeta|)$$

(using $r \leq \alpha_2(r)$). This can be written as (12) with $\rho := \alpha_1^{-1}(2(\cdot))$ and $\beta(s, t) := \alpha_1^{-1}(4e^{-t}\alpha_2(s))$. ∎

Implication $2 \Rightarrow 3$ of Theorem 1 now follows from Lemma 5.2 and Proposition 5.3.

We now turn to the proof of implication $3 \Rightarrow 1$ of Theorem 1.

*Proof.* Assume that $(\Sigma_{n.o}, k)$ is some norm-estimator for $\Sigma$. Choose any initial state $\xi$ for $\Sigma$, any input $\mathbf{u}$, disturbance $\mathbf{w}$, and the special initial state $\zeta = 0$ for $\Sigma_{n.o}$. Then inequality (11) becomes

$$|k(p(t, 0, \mathbf{u}, \mathbf{y}_{\xi,\mathbf{u},\mathbf{w}}), y(t, \xi, \mathbf{u}, \mathbf{w}))| \leq \hat{\gamma}_1\left(\|\mathbf{u}|_{[0,t]}\|\right) + \hat{\gamma}_2\left(\|\mathbf{y}_{\xi,\mathbf{u},\mathbf{w}}|_{[0,t]}\|\right) \tag{110}$$

for all $t \in [0, t_{\max})$, with some class-$\mathcal{K}$ functions $\hat{\gamma}_1$ and $\hat{\gamma}_2$ (the $\mathcal{KL}$-term vanishes because $\zeta = 0$). Then, combining (110) with the estimate (12) we get

$$\begin{aligned}|x(t, \xi, \mathbf{u}, \mathbf{w})| &\leq \beta(|\xi|, t) + \rho(|k(p(t, 0, \mathbf{u}, \mathbf{y}_{\xi,\mathbf{u},\mathbf{w}}), y(t, \xi, \mathbf{u}, \mathbf{w}))|) \\ &\leq \max\left\{2\beta(|\xi|, t), 4\rho(\hat{\gamma}_1(\|\mathbf{u}|_{[0,t]}\|)), 4\rho(\hat{\gamma}_1(\|y_{\xi,\mathbf{u},\mathbf{w}}|_{[0,t]}\|))\right\}.\end{aligned}$$

This proves the UIOSS property for $\Sigma$. ∎

# 6 Pointers for future research

## 6.1 Integral variants of UIOSS

The UIOSS property gives uniform estimates of the states in terms of the uniform bounds on outputs and essential bounds on controls. A natural question to ask is what property will result if instead of the uniform (or essential) bounds we use other "finite energy" concepts, such as, for example, $L^\gamma$-type norms (defined by $\|g|_{[\sigma,\tau]}\|_\gamma = \int_\sigma^\tau \gamma(|g(t)|)dt$) of inputs and/or outputs, where the "$\gamma$'s" for inputs and outputs are some appropriately chosen functions of class $\mathcal{K}_\infty$, which depend on the system. For systems without controls, the iiUOSS property provides a "finite energy outputs $\Rightarrow$ finite energy state" characterization, which is "almost" equivalent to UOSS (See Theorem 2). Searching for the *uniform* estimate of the states in terms of $L^\gamma$ norms of inputs and outputs leads to the following definition:

**Definition 6.1** A system of type (5) is *uniformly integral input-output to state stable*, if there exist functions $\alpha_x \in \mathcal{K}_\infty$, $\beta \in \mathcal{KL}$, $\gamma_1$, and $\gamma_2 \in \mathcal{K}$ such that

$$\alpha_x\left(|x(t,\xi,\mathbf{w},\mathbf{u})|\right) \leq \beta(|\xi|,t) + \int_0^t \left(\gamma_1(|\mathbf{u}(s)|) + \gamma_2\left(|y(s,\xi,\mathbf{w},\mathbf{u})|\right)\right) ds \tag{111}$$

for all $\xi \in \mathbb{X}$, all $\mathbf{w}$ and $\mathbf{u}$, and all $t \in [0, t_{\max}(\xi,\mathbf{u},\mathbf{w}))$.

This general definition may be adjusted in obvious manners to all the particular cases of system (5).

**Remark 6.2** It is easy to see that estimate (111) is equivalent to the following estimate (with different bounding functions, of course)

$$|x(t,\xi,\mathbf{w},\mathbf{u})| \leq \max\left\{\beta(|\xi|,t), \gamma\left(\int_0^t (\gamma_1(|\mathbf{u}(s)|))ds\right), \gamma\left(\int_0^t \gamma_2(|y(s,\xi,\mathbf{w},\mathbf{u})|)ds\right)\right\}, \tag{112}$$

□

In the particular case of systems without outputs and disturbances, a Lyapunov characterization of the UiIOSS property, reduced to integral input to state stability (iISS), was obtained in [4]. By repeating the proof of the implication $1 \Rightarrow 2$ of Theorem 1 in [4] one can show that a system of type (5) will be UiIOSS if it admits a smooth, proper Lyapunov function $V : \mathbb{X} \to R_{\geq 0}$, satisfying inequality (8) with some $\sigma_1$ and $\sigma_2$ of class $\mathcal{K}$, and a *positive definite* function $\alpha$. Whether or not this sufficient condition is also necessary for UiIOSS is not known. Notice, however, that this condition is weaker than the corresponding property for UIOSS, as the dissipation condition for a UIOSS-Lyapunov function requires $\alpha$ to be of class $\mathcal{K}_\infty$. Thus, any UIOSS system will also be UiIOSS. The converse implication is not true, as demonstrated in the following example.

**Example 6.3** The construction is similar to the one used in Remark 3.10, so, we recall that $\phi_\varepsilon(\cdot)$ denotes a $C^\infty$-bump function as in (36), and $1_A(\cdot)$ is the indicator function of a set $A$.

Choose $\varepsilon_f = 0.1$ and any $\varepsilon_h < (1-\varepsilon_f)e^{-1}$, and consider the autonomous system

$$\Sigma_1: \quad \dot{x} = f(x); \quad y = h(x)$$

with

$$f_1(x) = x\left[1_{(-\infty,-1]}(x)(1-\phi_{\varepsilon_f}(x+1)) + 1_{[1,+\infty)}(x)(1-\phi_{\varepsilon_f}(x-1))\right]$$
$$-x\left[1_{(-1,1)}(x)(1-\phi_{\varepsilon_f}(x+1))(1-\phi_{\varepsilon_f}(x-1))\right],$$

and

$$h_1(x) = 1 - \phi_{\varepsilon_h}(x).$$

evolving in $\mathbb{X} = \mathbb{R}$. We claim that $\Sigma_1$ is iOSS but not OSS.

Consider also an autonomous system $\Sigma_2$ on $\mathbb{R}$ with

$$f_2(x) = x; \qquad h_2(x) \equiv 1 \tag{113}$$

This system cannot serve as a counterexample, because $h_2(0) \neq 0$. However, its behavior away from 0 is identical to that of $\Sigma_1$, so that considering it will simplify the presentation.

Let $x_i(t,\xi)$, $i = 1, 2$, denote the solutions of $\Sigma_i$, starting at $\xi$, and let $y_i(t,\xi)$ denote the corresponding output trajectories. It is easy to see that both $\Sigma_1$ and $\Sigma_2$ are forward complete.

Since the dynamics of $\Sigma_1$ and $\Sigma_2$ are odd functions and outputs are even (but only the magnitudes of outputs are involved in the estimates), we only need to consider trajectories starting from the positive initial states, as the same estimates will work for trajectories contained in the other half-line.

Since $f_1(x) < f_2(x)$ for all $x \in (0, 1+\varepsilon_f)$ and $f_1(x) = f_2(x)$ for $x \geq 1+\varepsilon_f$, we have, for all $t > 0$, $x_1(t,\xi) < x_2(t,\xi)$ for any $\xi \in (0, 1+\varepsilon_f)$, and $x_1(t,\xi) = x_2(t,\xi) = \xi e^t$ if $|\xi| \geq 1+\varepsilon_f$. In particular, this shows that $\Sigma_1$ is not OSS, because the trajectory diverges to $+\infty$, but $h_1(\xi e^t)$ is bounded.

However, observe that when $t \leq |\xi|$,

$$|x_2(t,\xi)| \leq |\xi| e^{|\xi|} \leq e^{2|\xi|} |\xi| e^{-\frac{t}{1+|\xi|}},$$

whereas when $t > |\xi|$, we have

$$|x_2(t,\xi)| \leq te^t.$$

Letting

$$\beta(r,t) := re^{2r-\frac{t}{1+r}}; \qquad \gamma_2 := \mathrm{Id}; \qquad \gamma(t) := te^t,$$

and noticing that

$$\int_0^t \gamma_2(|y_2(s,\xi)|)ds = t,$$

we conclude $x_2(t,\xi)$ satisfies estimate (112).

Now observe that

- If $\xi \geq 1$, then $1 \leq x_1(t,\xi) \leq x_2(t,\xi)$ for all $t \geq 0$, so that $y_1(t,\xi) = y_2(t,\xi) = 1$ and $x_1(t,\xi)$ satisfies (112).

- If $\xi \in [0, 1-\varepsilon_f]$, then
$$x_1(t,\xi) = \xi e^{-t} \leq \beta(|\xi|, t).$$

- Finally, if $\xi \in (1-\varepsilon_f, 1)$, then
  - for all $t \geq 0$ it holds that $|x_1(t,\xi)| < 1$;

- for all $t \in [0,1]$ it holds that $|x_1(t,\xi)| > (1-\varepsilon_f)e^{-1} > \varepsilon_h$, so that $y_1(t,\xi) = 1 = y_2(t,\xi)$.

Therefore $x_1(t,\xi)$ satisfies (112) for all $t \leq 1$ and

$$\int_0^t \gamma_2(|y_2(s,\xi)|)ds \geq 1 \geq x_1(t,\xi) \quad \forall\, t \geq 1.$$

This shows that $x_1(t,\xi)$ satisfies (112) for all $\xi$ and $t \geq 0$, thus, $\Sigma_1$ is, indeed, iOSS.

It could also be of interest to consider yet another property, providing a uniform estimate for the state in terms of the uniform norm of the output and $L^\gamma$ norm of the input.

**Definition 6.4** A system of type (5) satisfies the U(iI)OSS property if there exist functions $\alpha_x \in \mathcal{K}_\infty$, $\beta \in \mathcal{KL}$, $\gamma_1$, and $\gamma_2 \in \mathcal{K}$ such that

$$\alpha_x\left(|x(t,\xi,\mathbf{w},\mathbf{u})|\right) \leq \beta(|\xi|,t) + \int_0^t \gamma_1(|\mathbf{u}(s)|)ds + \gamma_2\left(\|y_{\xi,\mathbf{w},\mathbf{u}}|_{[0,t]}\|\right) \tag{114}$$

for all $\xi \in \mathbb{X}$, all $\mathbf{w}$ and $\mathbf{u}$, and all $t \in [0, t_{\max}(\xi, \mathbf{u}, \mathbf{w}))$.

Deriving a Lyapunov characterization for this property may be a challenge. It would be logical to expect that a good candidate for a U(iI)OSS-Lyapunov function can be a proper function $V: \mathbb{X} \to \mathbb{R}_{\geq 0}$, satisfying inequality (8) with $\alpha$ being either positive definite or class $\mathcal{K}_\infty$. However, we can see right away that neither of these two possibilities is the right guess. Indeed, notice that both OSS and iISS are particular cases of the U(iI)OSS property. If a U(iI)OSS-Lyapunov function would satisfy (8) with $\alpha$ of class $\mathcal{K}_\infty$ (as required by OSS), it would follow that every disturbance-free U(iI)OSS system without outputs is ISS, which is not true (see [4] for an example of an iISS system, which is not ISS). On the other hand, if a U(iI)OSS-Lyapunov function satisfied (8) with $\alpha$ positive definite (as required by iISS), it would imply that having such a dissipation function is sufficient for OSS, which is not so, because this would mean that every iOSS system is OSS.

## 6.2 Incremental input-output to state stability

As mentioned in the introduction, the detectability property for nonlinear systems is not equivalent to zero-detectability. In seaching for a correct notion for nonlinear detectability one could think of the following generalization of UIOSS:

**Definition 6.5** A system (5) is *incrementally uniformly input-output to state stable* ($\Delta$UIOSS) if there exists some $\beta \in \mathcal{KL}$ and $\gamma_1, \gamma_2 \in \mathcal{K}$ such that, for every two initial states $\xi_1$ and $\xi_2$, any two controls $\mathbf{u}_1$ and $\mathbf{u}_2$, and any disturbance $\mathbf{w}$,

$$\begin{aligned}|x(t,\xi_1,\mathbf{u}_1,\mathbf{w}) - x(t,\xi_2,\mathbf{u}_2,\mathbf{w})| &\leq \max\{\beta(|\xi_1-\xi_2|,t),\\ &\gamma_1\left(\|(\mathbf{u}_1-\mathbf{u}_2)|_{[0,t]}\|\right), \gamma_2\left(\|(\mathbf{y}_{\xi_1,\mathbf{u}_1,\mathbf{w}} - \mathbf{y}_{\xi_2,\mathbf{u}_2,\mathbf{w}})|_{[0,t]}\|\right)\}\end{aligned} \tag{115}$$

for all $t$ in the common domain of definition.

Deriving a right Lyapunov characterization for this property may lead to a construction of a full-order observer.